\documentclass[10pt]{article}

\usepackage{latexsym}
\usepackage{amsmath,amsfonts}
\usepackage[dvips]{epsfig}
\usepackage[margin=1in]{geometry}
\usepackage{amssymb}
\usepackage{appendix}
\usepackage{color}

\newcommand{\q}{\quad}
\newcommand{\p}{\partial}
\newcommand{\e}{\epsilon}

\newcommand{\dv}{{\rm div}}
\newcommand{\0}{\|_{L^2}}

\newcommand{\2}{\|_{H^2}}

\newcommand{\bl}{\color{black}}

\newcommand{\ef}{\eqref}
\DeclareMathOperator*{\esssup}{ess\,sup}

\newcommand{\be}{\begin{equation}}
\newcommand{\ee}{\end{equation}}

\newtheorem{thm}{Theorem}[section]
\newtheorem{coro}[thm]{Corollary}
\newtheorem{lem}[thm]{Lemma}
\newtheorem{prop}[thm]{Proposition}
\newtheorem{rmk}[thm]{Remark}
\newtheorem{defn}[thm]{Definition}

\numberwithin{equation}{section}
\begin{document}

\begin{center}\large
\textbf{Low Mach number limit for non-isentropic magnetohydrodynamic equations with ill-prepared data and zero magnetic diffusivity in bounded domains}
\end{center}

\begin{center}
{Yaobin Ou,  Lu Yang}
\vspace{1.6mm}\\
School of Mathematics, Renmin University of China, \\ Beijing 100872, P.R. China.\\
Institute of Applied Physics and Computational Mathematics, \\ Beijing 100088, P.R. China.
\end{center}
\begin{center}\small Emails:
ou@ruc.edu.cn,  yanglu96@ruc.edu.com.
\end{center}

\noindent{\bf Abstract}. In this article, we verify the low Mach number limit of strong solutions to the non-isentropic compressible magnetohydrodynamic equations with zero magnetic diffusivity and {\it ill-prepared} initial data in three-dimensional {\it bounded domains}, when the density and the temperature vary around constant states. Invoking a new weighted energy functional,  we establish the uniform estimates with respect to the Mach number, especially for the spatial derivatives of high order. Due to the vorticity-slip boundary condition of the velocity, we decompose the uniform estimates into the part for the fast variables and the other one for the slow variables.  In particular, the weighted estimates of  highest-order spatial derivatives of the fast variables are  crucial for the uniform bounds. Finally, the low Mach number limit is justified by the strong convergence of the density and the temperature, the divergence-free component of the velocity, and the weak convergence of other variables. The methods in this paper can be applied to singular limits of general hydrodynamic equations of hyperbolic-parabolic type,  including the full Navier-Stokes equations.

\vskip 2mm
\noindent {{{\textbf{Key Words:}} \
Low Mach number limit, non-isentropic magnetohydrodynamic equations, ill-prepared initial data, zero magnetic diffusivity,  strong solutions,  bounded domains
\vskip 2mm
\noindent {{{\textbf{AMS subject classification.}} \ 35B40, 35M33, 35Q35, 76W05

\section{Introduction}

This paper is concerned with the low Mach number limit of compressible fluid flows which are described by the following non-dimensional non-isentropic magnetohydrodynamic equations without magnetic diffusivity in a three-dimensional bounded domain $\Omega\subset \mathbb{R}^3$:
\begin{equation}\label{ns}\left\{\begin{split}
&\rho^\epsilon_t+{\rm div}(\rho^\epsilon\mathbf{u}^\epsilon)=0,\\
&(\rho^\epsilon\mathbf{u}^\epsilon)_t+{\rm div}(\rho^\epsilon\mathbf{u}^\epsilon\otimes\mathbf{u}^\epsilon)-{\rm div}(2\mu D(\mathbf{u}^\epsilon)+\xi{\rm div}\mathbf{u}^\epsilon I)+\frac{1}{\epsilon^2}\nabla p^\epsilon=({\rm curl}\mathbf{H}^\epsilon)\times\mathbf{H}^\epsilon,\\
&\partial_t\mathbf{H}^\epsilon-{\rm curl}(\mathbf{u}^\epsilon\times\mathbf{H}^\epsilon)=0,\quad {\rm div}\mathbf{H}^\epsilon=0,\\
&(\rho^\epsilon e^\epsilon)_t+{\rm div}(\rho^\epsilon\mathbf{u}^\epsilon e^\epsilon)+p^\epsilon{\rm div}\mathbf{u}^\epsilon-{\rm div}(\kappa\nabla\mathcal{T}^\epsilon)=\epsilon^2(2\mu|D(\mathbf{u}^\epsilon)|^2+\xi({\rm div}\mathbf{u}^\epsilon)^2),
\end{split}\right.\end{equation}
where the variables $\rho^\epsilon$, $\mathbf{u}^\epsilon=(u^\epsilon_1,u^\epsilon_2,u^\epsilon_3)^t$, $\mathbf{H}^\epsilon=(H^\epsilon_1,H^\epsilon_2,H^\epsilon_3)^t$,  $p^\epsilon$, $e^\epsilon$, $\mathcal{T}^\epsilon$ stand for the density, velocity field, magnetic field, pressure, internal energy and temperature, respectively. The Mach number $\epsilon\in (0,1]$ is a dimensionless parameter that measures the compressibility of the fluid. The constants $\mu$ and $\xi$ denote respectively the shear and bulk viscosity coefficients, which  satisfy the physical conditions $\mu>0$ and $\mu+\frac{3}{2}\xi\geq0$, and the constant  $\kappa>0$ is the heat conductivity coefficient. The notation $D(\mathbf{u}^\e)=(\nabla\mathbf{u}^\e+(\nabla\mathbf{u}^\e)^t)/2$ denotes the viscous stress tensor. Moreover, we assume that the fluid is ideal polytropic, which satisfies
\begin{equation}\label{nonisentropic}
e^\epsilon=c_v\mathcal{T}^\epsilon,\quad p^\epsilon=R\rho^\epsilon\mathcal{T}^\epsilon
\end{equation}
with $c_v>0$ and $R>0$ being the specific heat at constant volume and the generic gas constant, respectively.

From the physical standpoint, as the Mach number vanishes, compressible fluids behave similarly to the incompressible ones. By a formal derivation with asymptotic expansions (Lions \cite{Lions1}), as the Mach number $\epsilon$ tends to zero,
the solution $(\rho^\e, \mathbf{u}^\e, \mathbf{H}^\e,\frac{1}{\e^2}\nabla p^\e)$ of  \ef{ns} will converge to $(\bar{\rho}, \mathbf{v}, \mathbf{h}, \nabla\pi)$, where   $(\bar{\rho},\mathbf{v}, \mathbf{h}, \pi)$ satisfies the following incompressible MHD equations:
\begin{equation}\label{ins0}\left\{
  \begin{split}
& \bar{\rho}_t+{\rm div}(\bar{\rho}\mathbf{v})=0,\\
& (\bar{\rho}\mathbf{v})_t+{\rm div}(\bar{\rho}\mathbf{v}\otimes\mathbf{v})-{\rm div}(2\mu D(\mathbf{v})+\xi{\rm div}\mathbf{v}I)+\nabla\pi=({\rm curl}\mathbf{h})\times\mathbf{h},\\
&\p_t\textbf{h} -{\rm curl}(\textbf{v}\times\textbf{h})=0,
\quad\textrm{div}\textbf{h}=0,\\
&\gamma{\rm div}\mathbf{v}={\rm div}[\frac{\kappa}{R}\nabla(\frac{1}{\bar{\rho}})],
\end{split}\right.
\end{equation}
where $\gamma=1+R/c_v$ denotes the ratio of specific heats. The rigourous verification of this singular limit is named as the low Mach number  limit. In particular, when the right-hand side of $\eqref{ins0}_4$ equals to zero, this singular limit is also known as the incompressible limit.

The investigation on mathematical theory of low Mach number limit is very challenging, due to the singularity in the pressure term $\frac 1{\e^2}\nabla p^\e$ which prevents us from establishing the uniform estimates and the convergence.
Let us review some results on the low Mach number limit of the magnetohydrodynamic (MHD) equations as follows in several different contexts, depending on the regularity of solutions (strong or weak), the states of the system (isentropic or non-isentropic), the boundaries of the domain (whole space, torus or  domains with boundaries), as well as the generality of the initial data (well-prepared or ill-prepared).

For the hydrodynamic system of hyperbolic type in isentropic regime ($p^\e=p(\rho^\e)$),  Klainerman and Majda \cite{KM1,KM2} set up a general framework for studying the incompressible limit in $\mathbb{R}^n$ or $\mathbb{T}^n$ for the local smooth solutions with the  well-prepared  data  that the initial data of  ${\rm div}\mathbf{u}^\e$ is of $O(\epsilon)$. For the case of solid boundary,
Fan, Li and Nakamura \cite{FLN} obtained the convergence as Mach number tended to zero for
global strong solutions to the initial-boundary value problem of 3-D MHD equations with small initial data.
When the initial data are ill-prepared, which means that $\dv \textbf{u}^\e=O(1)$ initially, Hu and Wang \cite{HW2009} rigorously justified the incompressible limit of global weak solutions of the viscous compressible MHD equations to the weak solutions of the incompressible viscous systems with various boundary conditions. Moreover, the case of periodic boundary condition was verified by Jiang, Ju and Li \cite{JJL} via a modulated energy method. While
Li and Mu \cite{LM2018} showed the low Mach number limit for the strong solution to the compressible MHD equations in the life-time of the incompressible ones by filtering the acoustic waves. Recently, Gu, Ou and Yang \cite{GOY} justified the incompressible limit of local strong solutions to the isentropic compressible MHD equations with ill-prepared initial data in a bounded domain by establishing weighted energy estimates in both cases of viscous and inviscid magnetic field.  While Chen et al. \cite{CHPS} proved the global well-posedness of classical solution and the exponetial decay rate to the compressible MHD equations for the regular initial data with small energy but possibly large oscillations in 3-D bounded domains with slip boundary condition and vacuum.

In the non-isentropic regime, the problem is much more involved because the large pressure gradient of $O(1/\epsilon^2)$ in the momentum equation relies on both the density and the temperature. In the framework of hyperbolic-parabolic coupled systems, Volpert and Hudjaev \cite{VH} established the local existence of classical solutions to the Cauchy problem of \eqref{ns}. In a later work, Hu and Wang \cite{HW2010} investigated the global existence and large-time behavior of weak solutions to \eqref{ns} with large initial data in a three-dimensional bounded domain by using the Faedo-Galerkin method and vanishing viscosity method. For the case of well-prepared initial data, Wang and Xu \cite{WangXu2015} verified the incompressible limit of the non-isentropic MHD equations with zero thermal coefficient in bounded domains. While Cui, Ou and Ren \cite{COR2015} studied the low Mach number limit for the full compressible MHD equations with well-prepared initial data in a 3-D bounded domain. When the initial data is ill-prepared, Jiang, Ju and Li \cite{JJL2016} justified the low Mach number limit of full MHD equations without magnetic diffusivity in $\mathbb{R}^d$ $(d=2,3)$. For the case of large temperature variation, Jiang et al. \cite{JJLX} showed the incompressible limit for local classical solutions to the non-isentropic compressible ideal magnetohydrodynamic equations with general initial data in the whole space.  As the Mach number, viscosity and the heat conductivity go to zero, Kwon \cite{Kwon} studied the incompressible inviscid limit of the full MHD equations on expanding domains with general initial data and obtained the convergence rate.
Very recently, Sun \cite{Sun} showed the uniform estimates with respect to the Mach number for the full Navier-Stokes equations with large temperature variations in smooth domains with boundaries by the method of conormal Sobolev spaces which minimizes the use of normal derivatives near the boundary, and proved  the convergence in exterior domains as the Mach number tends to zero, due to the decay of local energy at far field. However, due to the oscillation of fast acoustic waves, the low Mach number limit for the strong solutions of the full MHD equations with ill-prepared data in bounded domains is not yet completely solved.

For other results on the singular limits of hydrodynamic equations with small Mach number, one may refer to \cite{CJS,Danchin,Feireisl2010,Novotny3,Hoff,HW2009,JJL,Ou} and the references therein for details. We remark that the mechanism for the incompressibility is different from the setting in \cite{CH,DM,HOWY} where the incompressibility is achieved by letting $\xi\rightarrow+\infty$ in \ef{ns}.

The purpose of this paper is to verify the low Mach number limit for the strong solutions to the non-isentropic MHD equations \ef{ns} with ill-prepared initial data and  vorticity-slip boundary conditions in three-dimensional bounded domains, in the case of zero magnetic diffusivity and small temperature variation. The main difficulties lie in establishing uniform high-order estimates for the strong solutions, in particular for  high-order spatial derivatives, and dealing with rapidly oscillating waves, which will prevent taking to the limit, especially for the nonlinear terms.

First, the uniform estimates for high-order spatial derivatives are difficult to obtain since integration by parts is not valid for the estimates in bounded domains. Unlike the case of whole space \cite{Alazard,JJLX}, the technique of pseudo-differential operators or the integration by parts for high-order spatial derivatives can not be applied to the case of bounded domain. To overcome this trouble, we develop the weighted uniform estimates for the strong solutions, which is motivated by \cite{OY}, from the low-order temporal and spatial derivatives to the high-order ones and use suitable multipliers to eliminate the singular terms. By making use of the vorticity-slip boundary conditions, we establish the estimates of fast variables {\bl$(\nabla\sigma^\e,\nabla\theta^\e,{\rm div} \mathbf{u}^\e)$} and slow variables {\bl$({\rm curl} \mathbf{u}^\e,\mathbf{H}^\e)$}, respectively, which is different from the Dirichlet case {\bl for weak solutions} in \cite{DGLM}. {\bl Since integrating by parts is usually invalid for estimating the high-order spatial derivatives in the boundary case, we have to control the vorticity and the divergence of the velocity by different ways. When estimating each derivative, we carefully balance the singular differential operator of order $\epsilon^{-1}$, which leads to the uniform estimates for the full norm. Moreover, weighted estimates for the time derivatives and some high-order spatial derivatives are introduced to close the estimates.}  We remark that  the weighted energy space in this paper is essentially different from the one in the isentropic case \cite{GOY}. Since the  density and the temperature are strongly coupled in the singular pressure term {\bl$\frac{1}{\epsilon^2}\nabla p^\epsilon$,} the norms of highest-order spatial derivatives of their variations are not uniformly bounded. Compared to \cite{OY} for the Navier-Stokes equations, the nonlinear terms $({\rm curl}\mathbf{H}^\epsilon)\times\mathbf{H}^\epsilon$ in $\eqref{ns}_2$  and  ${\rm curl}(\mathbf{u}^\epsilon\times\mathbf{H}^\epsilon)$ in  $\eqref{ns}_3$ change the structure of the system and create much more difficulty in establishing the estimates of highest-order spatial derivatives, and the zero magnetic diffusivity decreases the regularity of the magnetic field. Therefore, we establish the uniform estimates by using the special structure of the MHD equations.

{\bl Second,} the highly oscillating waves make trouble in the strong convergence of velocity field and the boundary effects of the acoustic waves also occur in the bounded domains. To solve this difficulty, we use the methods by  Secchi \cite{Secchi}  to get the strong convergence of the density and the temperature, the ``incompressible part'' $P\mathbf{u}^\e$ of the velocity $\mathbf{u}^\e$ and the weak convergence of $(I-P)\mathbf{u}^\e$, which gives the convergence of the nonlinear term  ${\rm div}(\rho^\e\mathbf{u}^\e\otimes \mathbf{u}^\e)$ and then yields the limiting equations. Note that the previous convergence results of non-isentropic Navier-Stokes system in bounded domains \cite{Dou,JO,JuOu} were established  for {\bl the situations of well-prepared initial data only, where $\dv \textbf{u}^\e_0 =O(\epsilon)$.} In those cases, the fast acoustic waves were not encountered, thus the strong convergence of the velocity followed, which is the main difference from this paper.

In what follows, for the simplicity, we will set $R=c_v=1$. Since the low Mach number fluid can be regarded as a perturbation near the background incompressible fluid \cite{Lions1}, where the density and temperature are usually set to be constants, we denote the density and temperature variations by $\sigma^\epsilon$ and $\theta^\epsilon$:
$$
\rho^\epsilon=1+\epsilon\sigma^\epsilon,\quad \mathcal{T}^\epsilon=1+\epsilon \theta^\epsilon.
$$
Then by virtue of the usual vorticity identities
$$
{\rm curl}(\mathbf{u}^\e\times\mathbf{H}^\e)=\mathbf{u}^\e{\rm div}\mathbf{H}^\e-\mathbf{H}^\e{\rm div}\mathbf{u}^\e+(\mathbf{H}^\e\cdot\nabla)\mathbf{u}^\e-(\mathbf{u}^\e\cdot\nabla)\mathbf{H}^\e,
$$
$$({\rm curl}\mathbf{H}^\e)\times\mathbf{H}^\e=(\mathbf{H}^\e\cdot\nabla)\mathbf{H}^\e
-\frac{1}{2}\nabla(|\mathbf{H}^\e|^2),$$
and {\bl${\rm div}\mathbf{H}^\e = 0$,} we can equivalently write the non-dimensional system \ef{ns} as follows
\begin{equation}\label{ns0}\left\{\begin{split}
&\sigma^\epsilon_t+{\rm div}(\sigma^\epsilon\mathbf{u}^\epsilon)+\frac{1}{\epsilon}{\rm div}\mathbf{u}^\epsilon=0,\quad & (t,x)\in(0,T]\times\Omega,\\
&\rho^\epsilon(\mathbf{u}^\epsilon_t+\mathbf{u}^\epsilon\cdot\nabla\mathbf{u}^\epsilon)+\frac{1}{\epsilon}
(\nabla\sigma^\epsilon+\nabla\theta^\epsilon)+\nabla(\sigma^\epsilon\theta^\epsilon)\\
&\quad=\mu\Delta\mathbf{u}^\epsilon
+\nu\nabla{\rm div}\mathbf{u}^\epsilon+(\mathbf{H}^\e\cdot\nabla)\mathbf{H}^\e-\frac{1}{2}\nabla(|\mathbf{H}^\e|^2),\quad & (t,x)\in(0,T]\times\Omega,\\
&\mathbf{H}^\epsilon_t+({\rm div}\mathbf{u}^\epsilon)\mathbf{H}^\epsilon+(\mathbf{u}^\epsilon\cdot\nabla)\mathbf{H}^\epsilon
-(\mathbf{H}^\epsilon\cdot\nabla)\mathbf{u}^\epsilon=0,\quad{\rm div}\mathbf{H}^\epsilon=0,\quad & (t,x)\in(0,T]\times\Omega,\\
&\rho^\epsilon(\theta^\epsilon_t+\mathbf{u}^\epsilon\cdot\nabla\theta^\epsilon) +(\rho^\epsilon\theta^\epsilon+\sigma^\epsilon){\rm div}\mathbf{u}^\epsilon+\frac{1}{\epsilon}{\rm div}\mathbf{u}^\epsilon=\kappa\Delta\theta^\epsilon+\epsilon(2\mu|D(\mathbf{u}^\epsilon)|^2+\xi({\rm div}\mathbf{u}^\epsilon)^2),\quad & (t,x)\in(0,T]\times\Omega,
\end{split}\right.\end{equation}
where $\nu=\mu+\xi$. We impose the following initial and boundary conditions:
\begin{equation}\label{initial}
(\sigma^\epsilon,\mathbf{u}^\epsilon,\mathbf{H}^\epsilon,\theta^\epsilon)|_{t=0}
=(\sigma^\epsilon_0,\mathbf{u}^\epsilon_0,\mathbf{H}^\epsilon_0,\theta^\epsilon_0)\quad {\rm in}\ \Omega,
\end{equation}
\begin{equation}\label{boundary}
\mathbf{u}^\epsilon\cdot\mathbf{n}=0,\quad \mathbf{n}\times{\rm curl}\mathbf{u}^\epsilon=0,\quad \frac{\partial\theta^\epsilon}{\partial\mathbf{n}}=0\quad {\rm on}\ \partial\Omega,
\end{equation}
where $\mathbf{n}$ is the unit outer normal vector to the boundary.

The local existence and uniqueness of the strong solution to \ef{ns0} for fixed Mach number can be obtained by the energy estimates (e.g., established in Theorem \ref{unithm}), the method of characteristics, the Galerkin approximation and the classical fixed point arguments as the works by Hoff \cite{Hoff2012} and Gu, Ou and Yang \cite{GOY}, thus we omit the details here.
{\thm\label{main1} {\rm (Local wellposedness for fixed $\epsilon$)} Let $\e \in (0,1]$ be a fixed number and $\Omega\subset \mathbb{R}^3$
be a simply connected, bounded domain with smooth boundary  $\partial\Omega$.
Suppose that the initial datum $(\sigma_0^\epsilon,\mathbf{u}_0^\epsilon,\theta_0^\epsilon)$ of \eqref{ns0} satisfies
\begin{equation*}
(\sigma_0^\epsilon,\mathbf{u}_0^\epsilon,\theta_0^\epsilon)\in H^3(\Omega),\quad \mathbf{H}_0^\epsilon\in H^2(\Omega),
\quad (\partial^k_t\sigma^\epsilon(0),\partial^k_t\mathbf{u}^\epsilon(0),
\partial^k_t\mathbf{H}^\epsilon(0),\partial^k_t\theta^\epsilon(0))\in H^{2-2k}(\Omega),\quad k=0,1,
\end{equation*}
with $1+\epsilon \sigma^\epsilon_0\geq m^\epsilon>0$ and $1+\epsilon \theta^\epsilon_0\geq m^\epsilon>0$ for some positive constant $m^\e$. Moreover, we assume the following compatibility conditions are satisfied:
\begin{equation}\label{compat}
\begin{split}
&\mathbf{u}^\epsilon_0\cdot\mathbf{n}=\mathbf{u}^\epsilon_t(0)\cdot\mathbf{n}=0,\quad {\rm curl}\mathbf{u}^\epsilon_0\times \mathbf{n}= {\rm curl}\mathbf{u}^\epsilon_t(0)\times\mathbf{n}=0,\quad
{\frac{\partial\theta^\epsilon_0}{\partial\mathbf{n}}
=\frac{\partial\theta^\epsilon_t(0)}{\partial\mathbf{n}}=0,}
\quad {\rm on}\ \partial\Omega.
\end{split}
\end{equation}
Then there exists a positive constant $T^\epsilon=T^\epsilon(\sigma_0^\epsilon,\mathbf{u}_0^\epsilon,\mathbf{H}_0^\epsilon,\theta_0^\epsilon, m^\epsilon,\epsilon)$
such that the initial-boundary value problem \eqref{ns0}
admits a unique {\bl strong} solution $(\sigma^\epsilon,\mathbf{u}^\epsilon,\mathbf{H}^\epsilon,\theta^\epsilon)$,
satisfying that $1+\epsilon \sigma^\epsilon>0$ and  $1+\epsilon \theta^\epsilon>0$ in $\Omega\times(0,T^\epsilon)$, and
\begin{equation}
\begin{split}\label{regularity}
& (\sigma^\epsilon,\mathbf{H}^\epsilon)\in C([0,T^\epsilon],H^2(\Omega))\cap L^\infty (0,T^\epsilon;H^3(\Omega)),\quad (\mathbf{u}^\epsilon,\theta^\epsilon)\in C([0,T^\epsilon],H^3(\Omega))\cap L^2(0,T^\epsilon;H^4(\Omega)),\\
& \sigma^\epsilon_t\in C([0,T^\epsilon],H^1(\Omega))\cap L^2(0,T^\epsilon;H^2(\Omega)),\quad  \mathbf{H}^\epsilon_t\in C([0,T^\epsilon],H^1(\Omega))\cap L^\infty(0,T^\epsilon;H^2(\Omega)),\\
& (\sigma^\epsilon_{tt},\mathbf{H}^\epsilon_{tt})\in L^\infty( 0,T^\epsilon;L^2(\Omega)),\quad (\mathbf{u}_t^\epsilon,\theta_t^\epsilon)\in C([0,T^\epsilon],H^1(\Omega))\cap L^2(0,T^\epsilon;H^2(\Omega)),\\ &(\mathbf{u}_{tt}^\epsilon,\theta_{tt}^\epsilon)\in {L^\infty( 0,T^\epsilon;L^2(\Omega))}\cap L^2(0,T^\epsilon;H^1(\Omega)).
\end{split}
\end{equation}
}

\begin{rmk} The notation  $\mathbf{u}^\epsilon_t(0)$ is indeed defined by taking $t=0$ in $\ef{ns0}_2$ as follows:
\begin{equation*}
\begin{split}
\mathbf{u}^\epsilon
_t(0)=[-\frac{1}{\epsilon}
(\nabla\sigma^\epsilon_0+\nabla\theta^\epsilon_0)-\nabla(\sigma^\epsilon_0\theta^\epsilon_0)
+\mu\Delta\mathbf{u}^\epsilon_0
+\nu\nabla{\rm div}\mathbf{u}^\epsilon_0]/\rho^\epsilon_0-\mathbf{u}^\epsilon_0\cdot\nabla\mathbf{u}^\epsilon_0.
\end{split}
\end{equation*}
The notation $\sigma^\epsilon_t(0)$, $\mathbf{H}^\epsilon_t(0)$ and $\theta^\epsilon_t(0)$ are defined by letting $t=0$ in $\ef{ns0}_{1,3,4}$, respectively. Moreover, $\mathbf{u}^\epsilon_{tt}(0)$, $\sigma^\epsilon_{tt}(0)$, $\mathbf{H}^\epsilon_{tt}(0)$ and $\theta^\epsilon_{tt}(0)$ are defined inductively by applying $\partial_t$ to $\ef{ns0}_{1,2,3,4}$ and taking $t=0$.
\end{rmk}

The purpose of this paper is to prove the uniform estimates of local strong solutions  obtained in Theorem \ref{main1}  and the corresponding incompressible limit. In order to state the theorem precisely, we introduce the following notation.

\begin{defn}
We define the weighted energy as
\begin{equation*}\begin{split}
M^\epsilon(t):=& \esssup_{0\leq s\leq t} \{\|(\mathbf{u}^\epsilon,\rho^\epsilon,\mathbf{H}^\epsilon,\epsilon\sigma^\epsilon,\epsilon\theta^\epsilon)||_{H^2(\Omega)}^2(s)
+\|(\sigma^\epsilon,\theta^\epsilon)\|_{H^1(\Omega)}^2(s)+\|\epsilon( \sigma^\epsilon_t,  \mathbf{u}^\epsilon_t, \theta^\epsilon_t)\|_{H^1(\Omega)}^2(s)\\
&\qquad+\|\mathbf{H}^\epsilon_t\|_{H^1(\Omega)}^2(s)
+\|\epsilon^2\mathbf{H}^\epsilon_t\|_{H^2(\Omega)}^2(s)
+\|\epsilon^2 (\sigma^\epsilon_{tt}, \mathbf{u}^\epsilon_{tt}, \mathbf{H}^\epsilon_{tt}, \theta^\epsilon_{tt})\|_{L^2(\Omega)}^2(s)\\
&\qquad
 +\| (\epsilon \mathbf{u}^\epsilon,\epsilon\mathbf{H}^\epsilon,\epsilon^2 \sigma^\epsilon,\epsilon^2 \theta^\epsilon) \|_{H^3(\Omega)}^2(s)+\|(\rho^\epsilon,(\rho^{\epsilon})^{-1})\|_{L^\infty(\Omega)}^2(s)\}\\
&+\int_0^t \{\|(\sigma^\epsilon,\theta^\epsilon)\|_{H^2(\Omega)}^2
+\|(\epsilon \sigma^\epsilon,\mathbf{u}^\epsilon,\epsilon \theta^\epsilon)\|_{H^3(\Omega)}^2
+\|\epsilon( \sigma^\epsilon_t, \mathbf{u}^\epsilon_t,\theta^\epsilon_t)\|_{H^2(\Omega)}^2\\
&\qquad
+\|\epsilon^2(\mathbf{u}^\epsilon_{tt},\theta^\epsilon_{tt})\|_{H^1(\Omega)}^2 +\|(\epsilon^2\mathbf{u}^\epsilon,\epsilon^2\theta^\epsilon)\|_{H^4(\Omega)}^2\}ds.
\end{split}\end{equation*}
 and the initial weighted energy as
$$
M^\epsilon_0:=M^\epsilon(t=0).
$$
\end{defn}

{\begin{rmk}
We can see the initial data is ``ill-prepared" with the assumption that $M^\epsilon_0\leq C$, in particular,
$$
\|{\rm div}\mathbf{u}^\epsilon_0\|_{H^2}+\epsilon\|(\sigma_t^\epsilon(0),
\mathbf{u}_t^\epsilon(0),\theta_t^\epsilon(0))\|_{H^1}\leq C.
$$
In this situation, the fast oscillating waves are encountered in the low Mach number limit.
\end{rmk} }

In this paper, we aim to derive the following energy inequality by low-order and high-order temporal and spatial derivative estimates of $(\sigma^\epsilon,\mathbf{u}^\epsilon,\theta^\epsilon)$.

\begin{prop}\label{prop}
Let $\epsilon\in (0,1]$ be an arbitrary constant and $T^\epsilon>0$ be the maximal time of existence for the solution to the initial-boundary value problem $\ef{ns0}$ established in Theorem \ref{main1}.  Then for any  $\epsilon \in (0,1]$ and $t\in[0,T^\epsilon)$, we have
\begin{equation*}
M^\epsilon(t)\leq C_0(M^\epsilon_0){\rm exp}\{(t^{\frac{1}{8}}+ \epsilon^\frac 1 2 )C_1(M^\epsilon(t))\}.
\end{equation*}
 Here $C_0(\cdot)$ and $C_1(\cdot)$ are  given positive nondecreasing continuous functions.
\end{prop}

The following theorem for the uniform estimates of the strong solutions with ill-prepared initial data is based on Theorem \ref{main1} and Proposition \ref{prop}.

\begin{thm}\label{unithm} {\rm (Uniform estimates for ill-prepared initial data)}
 Let $\Omega\subset \mathbb{R}^3$
be a simply connected, bounded domain with smooth boundary  $\partial\Omega$ and $\epsilon\in (0,1]$ is an arbitrary constant. Suppose that the initial datum $(\sigma_0^\epsilon,\mathbf{u}_0^\epsilon,\mathbf{H}_0^\epsilon,\theta_0^\epsilon)\in H^3(\Omega)$ satisfies the assumptions in Theorem \ref{main1} and that $M^\epsilon_0\leq D_0$, where $D_0>0$ is a constant independent of $\e$. Assume in addition that there exists a sequence $(\sigma_0^\epsilon,\mathbf{u}_0^\epsilon,\mathbf{H}_0^\epsilon,\theta_0^\epsilon)\in H^5(\Omega)$ such that $(\sigma_0^{\epsilon,n},\mathbf{u}_0^{\epsilon,n},\mathbf{H}_0^{\epsilon,n},\theta_0^{\epsilon,n})\to (\sigma_0^\epsilon,\mathbf{u}_0^\epsilon,\mathbf{H}_0^\epsilon,\theta_0^\epsilon)$ in $H^3(\Omega)$ as $n\to\infty$ and $(\sigma_0^{\epsilon,n},\mathbf{u}_0^{\epsilon,n},\mathbf{H}_0^{\epsilon,n},\theta_0^{\epsilon,n})$ satisfies the compatibility conditions through order two. If $(\sigma^\e,\mathbf{u}^\e,\mathbf{H}^\e,\theta^\e)$ is the unique strong solution  to  ${\eqref{ns0}}$ for each $\e$, which is obtained in Theorem \ref{main1},   then there exist  positive constants $\epsilon_0$, $T$ and $C$, independent of $\e$, such that the following uniform estimate  holds:
\begin{equation}\label{uniform}
M^\epsilon(t)\leq C,\quad \epsilon\in(0,\epsilon_0],\; t\in[0,T].
\end{equation}
 \end{thm}

{\begin{rmk} Provided that $\mathbf{u}^\e\in H^2(\Omega)$ with $\mathbf{u}^\e\cdot\mathbf{n}=0$ on $\partial\Omega$, Xiao and Xin \cite{XX2013} showed that the vorticity-slip boundary condition $\mathbf{n}\times{\rm curl}\mathbf{u}^\e=0 $ and the Navier-slip boundary condition
$\tau\cdot(2S(\mathbf{u}^\e)\mathbf{n})+\alpha(x)\mathbf{u}^\e\cdot\tau=0$
 only differ from a low-order term in the following sense:
\begin{equation}\label{relation}
\tau\cdot(2S(\mathbf{u}^\e)\mathbf{n}-\mathbf{n}\times{\rm curl}\mathbf{u}^\e)=-2\tau\cdot(S(\mathbf{n})\mathbf{u}^\e)\q {\rm on}\q \partial\Omega.
\end{equation}
Here $S(\mathbf{u}^\e)=\frac 1 2(\nabla\mathbf{u}^\e+(\nabla\mathbf{u}^\e)^\top)$ is the stress tensor, $\tau$ is the unit tangential vector to $\partial\Omega$, and $\alpha(x)$ is a given function. Therefore, the results obtained in this paper can be extended to the case of Navier-slip boundary condition by small adaptations. In addition, our results can be also generalized to other hydrodynamic equations {\bl of hyperbolic-parablic type} with some extra efforts.
\end{rmk} }

We denote by $P$ the orthogonal projection onto incompressible vector fields, i.e.
$$v=Pv+Qv,\quad{\rm with}~~ \ {\rm div}(Pv)=0,\ {\rm curl}(Qv)=0 $$
for all $v\in L^2(\Omega).$ From the results by Galdi
\cite{Galdi}, the operators $P$ and $Q$ are linear and bounded in $W^{s,p}(\Omega)$ for all $s\ge 0$ and $1<p<\infty$ in a bounded domain $\Omega$ with smooth boundary $\partial \Omega$.

\begin{thm}\label{convergence} {\rm (Singular limit)}
Let the assumptions in Theorem \ref{unithm} be satisfied. Suppose that
$(\mathbf{u}_0^\epsilon,\mathbf{H}_0^\epsilon) \rightharpoonup  (\mathbf{w}_0,\mathbf{h}_0)$ in $H^2(\Omega)$ and $(\mathbf{u}_0^\epsilon,\mathbf{H}_0^\epsilon) \rightarrow  (\mathbf{w}_0,\mathbf{h}_0)$ in $H^1(\Omega)$ as $\epsilon \to 0$.
Then there exist a subsequence of $\{\mathbf{u}^\epsilon,\mathbf{H}^\epsilon\}$, still denoted as $\{\mathbf{u}^\epsilon,\mathbf{H}^\epsilon\}$, and
functions $ \mathbf{v} \in L^\infty(0,T;H^2(\Omega))\cap L^2(0,T;H^3(\Omega))$ and $ \mathbf{h} \in   L^\infty(0,T;H^2(\Omega))$, such that
$$
\mathbf{u}^\epsilon\rightharpoonup\mathbf{v}\quad  weakly\ in\ L^2(0,T;H^3(\Omega))\ and \ weakly-* \ in \ L^\infty(0,T;H^2 (\Omega)),
$$
$$
P\mathbf{u}^\epsilon\rightarrow\mathbf{v}\quad strongly\ in\ L^2(0,T;H^{3-\delta}(\Omega))\cap C([0,T],H^{2-\delta}(\Omega)),{\bl\ for\ any\ \delta\in(0,1),}
$$
$$
Q\mathbf{u}^\epsilon\rightharpoonup 0\quad weakly\ in\ L^2(0,T;H^3(\Omega)) \ and \ weakly-* \ in \ L^\infty(0,T;H^2 (\Omega)),
$$
$$
\mathbf{H}^\epsilon\rightharpoonup \mathbf{h}\quad strongly\ in\ C([0,T],H^{2-\delta}(\Omega)),{\bl\ for\ any\  \delta\in(0,1),}
$$
as $\epsilon \to 0$. And there exists a function $P(x,t)$, such that {\bl $(\mathbf{v},P,\mathbf{h})$} is the unique solution of the following initial-boundary value problem of the incompressible {\bl MHD} equations:
\begin{equation}\label{ins}\left\{
  \begin{split}
& \mathbf{v}_t+\mathbf{v}\cdot\nabla\mathbf{v}+\nabla P=\mu\Delta\mathbf{v},\quad {\rm div}\mathbf{v}=0\quad {\rm in}\ \Omega\times(0,T],\\
&\mathbf{h}_t-\nabla\times(\mathbf{v}\times\mathbf{h})=0,
\quad{\rm div}\mathbf{h}=0,\quad {\rm in}\ \Omega\times(0,T],\\
&\mathbf{v}\cdot\mathbf{n}=0,\quad \mathbf{n}\times{\rm curl}\mathbf{v}=0\ \ {\rm on}\ \partial\Omega\times(0,T],\\
&\mathbf{v}(x,0)=P \mathbf{w}_0,\quad x\in\Omega.
\end{split}\right.
\end{equation}
\end{thm}

{\rmk Unlike the Dirichlet case \cite{JN2015} {\bl for non-isentropic Navier-Stokes equations}} that $\mathbf{u}^\e=0$ and $\theta^\e=0$ on the boundary, the strong convergence of $\mathbf{u}^\e$ is not available in this paper since there is no dissipation in the boundary layers of the acoustic waves for the boundary conditions \eqref{boundary}. In addition, the entire family $\mathbf{u}^\e$ converges to $\mathbf{v}$ since the local strong solution $\mathbf{v}$ of \eqref{ins} is unique \cite{OY}.
}

The rest of this paper is organized as follows. In Section 2, we state some elementary lemmas to be used in estimating the Sobolev norms and present some inequalities to deal with the {\bl vorticity-slip} boundary condition. In Section 3, the non-standard uniform estimates with respect to the Mach number are shown by employing a careful weighted analysis on both temporal and spatial derivatives. We first derive the estimates for the density $\rho^\epsilon$ and then the basic $L^2$ estimates for $(\sigma^\epsilon,\mathbf{u}^\epsilon,\mathbf{H}^\epsilon,\theta^\epsilon)$. Next, the weighted estimates for the first-order time derivatives are shown, in particular for $\e (\nabla \sigma_t^\e,\nabla \theta_t^\e,{\rm div} \mathbf{u}_t^\e)$. Then we carry out the estimates for the slow component ${\rm curl}\mathbf{u}^\e$, including the weighted estimates for its time derivatives and the highest-order spatial derivatives. Finally, we establish the weighted estimates for the second-order time derivatives $\e^2 (\sigma_{tt}^\e,\mathbf{u}_{tt}^\e,\mathbf{H}_{tt}^\e,\theta_{tt}^\e)$ and some subtle weighted estimates for the highest-order derivatives of the fast components  $ (\nabla \sigma^\e,\nabla \theta^\e)$ and ${\rm div} \mathbf{u}^\e$ (the weight could be $\e$ or $\e^2$), which are critical to close the entire uniform estimates.  By collecting all the spatial-temporal estimates, we obtain the uniform estimates with respect to the Mach number. In Section 4, we give the proofs of Theorem \ref{unithm} and Theorem \ref{convergence}.

\indent Before ending this section, we give the notations used throughout this paper. We denote the usual Sobolev spaces $W^{k,p}(\Omega)$ and {\bl $H^k(\Omega)\equiv W^{k,2}(\Omega)$} endowed with the norms $\|\cdot\|_{W^{k,p}}$ and $\|\cdot\|_{H^k}$, respectively. We shall also use the following abbreviations for the Sobolev spaces involving time:
$$
L^p_t(H^k)\equiv L^p(0,t;H^k(\Omega)),\quad C_t(H^k)\equiv C([0,t],H^k(\Omega)),
$$
with
$$
\|\cdot\|_{L_t^p(H^k)}\equiv\|\cdot\|_{L^p(0,t;H^k)},
\quad\|\cdot\|_{C_t(H^k)}\equiv\|\cdot\|_{C([0,t],H^k)},
$$
where $H^0(\Omega)\equiv L^2(\Omega)$.
 We will utilize the notations $C$, $\delta$ and $C_\delta$ to signify a generic positive constant independent of $\epsilon$, a small positive constant and a positive constant depending on $\delta$, respectively. In addition, we use $f\lesssim g$ to represent the relation that $f\le C g$.

\section{Preliminaries}

In this section, we give some lemmas and inequalities that are often used in this paper.

\begin{lem}\label{hs1} {\rm (see Bourguignon and Brezis \cite{Brezis})} Let $\Omega$ be a bounded domain in $\mathbb{R}^N$ with smooth boundary $\partial\Omega$ and its unit outward normal $\mathbf{n}$. Then there exists a constant $C>0$ independent of $\mathbf{u}$, such that
$$
\|\mathbf{u}\|_{H^s(\Omega)}\leq C(\|{\rm div}\mathbf{u}\|_{H^{s-1}(\Omega)}+\|{\rm curl}\mathbf{u}\|_{H^{s-1}(\Omega)}
+\|\mathbf{u}\cdot\mathbf{n}\|_{H^{s-\frac{1}{2}}(\partial\Omega)}
+\|\mathbf{u}\|_{H^{s-1}(\Omega)}),
$$
for any $\mathbf{u}\in [H^s(\Omega)]^N$, $s\geq1$.
\end{lem}

\begin{lem}\label{interpolation}{\rm (see Friedman \cite{Frie}, Part 1, Theorem 10.1)} Let $\Omega\subset {\mathbb R}^N$
be a bounded domain with a $C^k$-boundary and  $u$ be any function
in $W^{k,r}(\Omega)\cap L^q(\Omega)$ with $1\le r,q\le\infty$. For
any integer $j$ with $0\le j<k$, and for any number $a$ in the
interval $[j/k, 1]$, set
$$\frac{1}{p}=\frac{j}{N}+a(\frac{1}{r}-\frac{k}{N})+(1-a)\frac{1}{q}.$$
 If $k-j-N/r$ is not a nonnegative integer, then
\begin{equation}\label{inter1}
\|D^j u\|_{L^p(\Omega)}\le C
\|u\|_{W^{k,r}(\Omega)}^a\|u\|_{L^q(\Omega)}^{1-a}.
\end{equation}
 If $k-j-N/r$ is a nonnegative integer, then \eqref{inter1} only holds for $a=j/k$.
 The constant $C$ depends only on $\Omega,r,q,k,j,a$.
\end{lem}

\begin{lem} \label{hs2}{\rm (see Xiao and Xin \cite{Xiao})} Let $\Omega$ be a bounded domain in $\mathbb{R}^N$ with smooth boundary $\partial\Omega$ and its unit outward normal $\mathbf{n}$. Then there exists a constant $C>0$ independent of $\mathbf{u}$, such that
$$
\|\mathbf{u}\|_{H^s(\Omega)}\leq C(\|{\rm div}\mathbf{u}\|_{H^{s-1}(\Omega)}+\|{\rm curl}\mathbf{u}\|_{H^{s-1}(\Omega)}
+\|\mathbf{u}\times\mathbf{n}\|_{H^{s-\frac{1}{2}}(\partial\Omega)}
+\|\mathbf{u}\|_{H^{s-1}(\Omega)}),
$$
for any $\mathbf{u}\in [H^s(\Omega)]^N$, $s\geq1$.
\end{lem}

For convenience, we introduce some special cases of interpolation inequality in $\Omega\subset\mathbb{R}^3$ which will be applied frequently in the following sections.

\begin{coro}\label{inter} Let $\Omega$ be a bounded, open domain of $\mathbb{R}^3$  with a $C^3$ boundary. Then the following inequalities are satisfied:
$$
\|u\|_{L^3}\leq C\|u\|_{H^1}^{\frac{1}{2}}\|u\|_{L^2}^{\frac{1}{2}},\quad\|\nabla^2u\|_{L^2}\leq C\|u\|_{H^3}^{\frac{1}{2}}\|\nabla u\|_{L^2}^{\frac{1}{2}},
$$
$$
\|\nabla^2u\|_{L^3}\leq C\|u\|_{H^3}^{\frac{1}{2}}\|u\|_{H^2}^{\frac{1}{2}},\quad\|\nabla^2u\|_{L^4}\leq C\|u\|_{H^3}^{\frac{3}{4}}\|u\|_{H^2}^{\frac{1}{4}},
$$
$$
\|\nabla u\|_{L^\infty}\leq C\|u\|_{W^{2,4}}\leq C\|u\|_{H^3}^{\frac{3}{4}}\|u\|_{H^2}^{\frac{1}{4}}, {\quad \|u\|_{L^\infty}\leq C\|u\|_{H^2}^{\frac{3}{4}}\|u\|_{H^1}^{\frac{1}{4}},}
$$
$$
\|u\|_{L^2_t(L^3)}\leq Ct^\frac{1}{4}\|u\|_{L^2_t(H^1)}^{\frac{1}{2}}\|u\|_{L^\infty_t(L^2)}^{\frac{1}{2}},
\quad
\|u\|_{L^2_t(L^\infty)}\leq Ct^\frac{1}{8}\|u\|_{L^2_t(H^2)}^{\frac{3}{4}}\|u\|_{L^\infty_t(H^1)}^{\frac{1}{4}}.
$$
\end{coro}

\begin{lem} \label{lembdry} {{\rm (see Ou and Yang \cite{OY})}}   Let $\Omega$ be a bounded domain in $\mathbb{R}^3$ with smooth boundary $\partial\Omega$ and its unit outward normal $n$. Suppose that $n\times {\rm curl}u=0$ on $\partial\Omega$ for any vector ${\bl u}\in { H^{k+2}}(\Omega)$, $k=0,1,2,\cdots$. Then there exists a constant $C=C(k,\Omega)>0$ independent of $u$, such that
$$
 \|{\rm curl}{\rm curl}u\cdot n\|_{H^{k+\frac{1}{2}}(\partial\Omega)}
\le C \|u\|_{H^{k+2}(\Omega)}.
$$
 \end{lem}

{
\begin{lem}\label{als}{\rm (see Simon \cite{Simon})} Let $X_0$, $X$ and $X_1$ be three Banach spaces with $X_0\subseteq X\subseteq X_1$. Suppose that $X_0$ is compactly embedded in $X$ and that $X$ is continuously embedded in $X_1$. For $1\leq p,q\leq+\infty$, let
$$
W=\{u\in L^p(0,T;X_0)|u_t\in L^q(0,T;X_1)\}.
$$
(i) If $p<+\infty$, then the embedding of $W$ into $L^p(0,T;X)$ is compact.\\
(ii) If $p=+\infty$ and $q>1$, then the embedding of $W$ into $C([0,T],X)$ is compact.
\end{lem}
}

\section{Uniform estimates}

\indent In what follows, for simplicity, we omit the superscript $\epsilon$ for the solutions.
We will show that,  for  $t\in [0,T]$ with $T=\min\{T^\epsilon,1\}$ and $\epsilon \in (0,1]$,
\begin{equation*}
M(t)\leq C_0(M_0)\exp\{C(t^\frac{1}{8}+{\epsilon^\frac{1}{2}})C_1(M(t))\},
\end{equation*}
which gives the uniform bound of $M(t)$ with respect to $\epsilon$, where $C_0(\cdot)$ and $C_1(\cdot)$ are generic nondecreasing positive and continuous functions.

\subsection{Upper and lower bounds of $\rho$}

It is convenient to derive the estimate of $\rho$, although $\rho$ only appears as the coefficient of $\mathbf{u}_t$. The following estimates can be obtained by the standard energy method since the density equation doesn't contain the large parameters of reciprocal of $\epsilon$.
\begin{lem}\label{lemrho}  Suppose that the assumptions in Theorem \ref{main1} are satisfied and $(\sigma,\mathbf{u},\mathbf{H},\theta)$ obtained in Theorem \ref{main1} is sufficiently smooth. Then for any $0<\epsilon\leq1$ and $0\leq t\leq T$, we have
\begin{equation}\label{lemma1}
\|(\rho,\rho^{-1})\|^2_{L^\infty}(t)+\|\rho\|^2_{H^2}(t)\leq C_0(M_0){\rm exp}\{t^\frac 1 2 C_1(M(t)) \}.
\end{equation}
\end{lem}
\textbf{Proof.} The proof is standard {(see also \cite{OY})} and we sketch it here for completeness. We multiply  ${\eqref{ns}}_1$ by $\rho^{-k}$ and integrate to derive that
\begin{equation}\label{rhola}
\frac{d}{dt}\|\rho^{-1}\|_{L^{k-1}}\leq \frac{k}{k-1}\|{\rm div}\mathbf{u}\|_{L^\infty}\|\rho^{-1}\|_{L^{k-1}}.
\end{equation}
By applying Gr\"{o}nwall's inequality, and letting $k\rightarrow+\infty$, we get
\begin{equation}\label{rholb}\begin{split}
\|\rho^{-1}\|_{L^\infty}(t)&\leq\|\rho_0^{-1}\|_{L^\infty}{\rm exp}\left(C\int_0^t
\|\mathbf{u}\|_{H^3}ds\right)\\
&\leq  C_0(M_0){\rm exp}\{t^\frac 1 2 C_1(M(t)) \}.
\end{split}\end{equation}
Similarly, one can derive  the upper bound of $\rho$.  Applying $\partial^\alpha$ to ${\eqref{ns}}_1$, $0\leq|\alpha|\leq2$, then multiplying the resulting equality by $\partial^\alpha\rho$ in  $L^2(\Omega)$ and summing over $\alpha$, we get that
\begin{equation*}\begin{split}
\frac{d}{dt}\|\rho\|^2_{H^2}\leq C\|\mathbf{u}\|_{H^3}\|\rho\|^2_{H^2},
\end{split}\end{equation*}
which implies that, for any $0<\epsilon\leq1$ and  $0\leq t\leq T$,
\begin{equation}\label{rhoh2}\begin{split}
\|\rho\|^2_{H^2}(t)\leq\|\rho_0\|^2_{H^2}{\rm exp}\left(C\int_0^t
\|\mathbf{u}\|_{H^3}ds\right)&\le C_0(M_0){\rm exp}\{t^\frac 1 2 C_1(M(t)) \}.
\end{split}\end{equation}
From \ef{rholb} and {\bl\ef{rhoh2}}, we can conclude \ef{lemma1} immediately.
\hfill$\Box$

\subsection{$L^2-$ estimates for the solution $(\sigma,\mathbf{u},\theta)$}

\begin{lem} Suppose that the assumptions in Lemma \ref{lemrho} are satisfied. For $\epsilon \in (0,1]$ and  $t\in[0,T]$, the solution to the system ${\eqref{ns0}}$ satisfies the following estimate
\begin{equation}\begin{split}\label{basicestimate}
&\|(\sigma,\mathbf{u},\mathbf{H},\theta)\|^2_{L^2}(t)
+\|\mathbf{u}\|^2_{L^2_t(H^1)}
+\|\nabla\theta\|^2_{L^2_t(L^2)}\leq  C_0(M_0){\rm exp}\{t^\frac 1 2 C_1(M(t)) \}.
\end{split}\end{equation}
\end{lem}
\textbf{Proof}. We multiply ${\eqref{ns0}}_1$ with $\sigma$ and integrate the resulting equation in $L^2(\Omega)$ to get
\begin{equation*}
\frac{1}{2}\frac{d}{dt}\|\sigma\|^2_{L^2}-\frac{1}{\epsilon}\int_\Omega\mathbf{u}\cdot\nabla\sigma dx=-\int_\Omega\sigma{\rm div}(\sigma\mathbf{u})dx\lesssim\|\mathbf{u}\|_{H^1}\|\sigma\|^2_{H^1}.
\end{equation*}
Next taking the inner product of $\langle{\eqref{ns0}}_2,\mathbf{u}\rangle$ and $\langle{\eqref{ns0}}_4,\theta\rangle$ in $L^2(\Omega)$, respectively, we derive that
\begin{equation*}\begin{split}
&\frac{1}{2}\frac{d}{dt} \|\sqrt{\rho}\mathbf{u}\|^2_{L^2} +\frac{1}{\epsilon}\int_\Omega\mathbf{u}\cdot\nabla(\sigma+\theta) dx+(\mu+\nu)\|{\rm div}\mathbf{u}\|^2_{L^2}+\mu\|{\rm curl}\mathbf{u}\|^2_{L^2}\\
&\quad=-\int_\Omega\mathbf{u}\cdot\nabla(\sigma\theta)dx
-\int_\Omega({\rm curl}\mathbf{H}\times\mathbf{H})\cdot\mathbf{u}dx\\
&\lesssim\|\mathbf{u}\|_{H^1}\|\sigma\|_{H^1}\|\theta\|_{H^1}
-\int_\Omega({\rm curl}\mathbf{H}\times\mathbf{H})\cdot\mathbf{u}dx,
\end{split}\end{equation*}
and
\begin{equation*}\begin{split}
&\frac{1}{2}\frac{d}{dt}\|\sqrt{\rho}\theta\|^2_{L^2}
+\kappa\|\nabla\theta\|^2_{L^2}-
\frac{1}{\epsilon}\int_\Omega\mathbf{u}\cdot\nabla\theta dx\\
&\quad=-\int_\Omega[(\rho\theta+\sigma){\rm div}\mathbf{u}-\epsilon(2\mu|D(\mathbf{u})|^2+\xi({\rm div}\mathbf{u})^2)]\cdot\theta dx\\
&\quad\lesssim\|\theta\|_{H^1}\|\mathbf{u}\|_{H^1}(\|\sigma\|_{H^1}+
\|\rho\|_{L^\infty}\|\theta\|_{H^1}+\|\mathbf{u}\|_{H^3}).
\end{split}\end{equation*}
To deal with the magnetic field equations, we take the inner product of $\langle{\eqref{ns0}}_3,\mathbf{H}\rangle$ to obtain that
\begin{equation*}\begin{split}
\frac{1}{2}\frac{d}{dt}\|\mathbf{H}\|^2_{L^2}=\int_\Omega(\mathbf{u}\times\mathbf{H})\cdot{\rm curl}\mathbf{H}dx.
\end{split}\end{equation*}
Finally, by use of $({\rm curl}\mathbf{H}\times \mathbf{H})\cdot \mathbf{u}
=-(\mathbf{u}\times\mathbf{H})\cdot{\rm curl}\mathbf{H}$, we put the above estimates together and integrate the resulting inequality on $[0,t]$ to obtain that
\begin{equation}\label{basic}\begin{split}
&\|\sigma\|^2_{L^2}(t)+\|\sqrt{\rho}\mathbf{u}\|^2_{L^2}(t)
+\|\sqrt{\rho}\theta\|^2_{L^2}(t)+\|\mathbf{H}\|^2_{L^2}(t)
+\|\mathbf{u}\|^2_{L^2_t(H^1)}
+\kappa\|\nabla\theta\|^2_{L^2_t(L^2)}\\
&\leq C_0(M_0)+
C\int_0^t\|\mathbf{u}\|_{H^1}(\|\sigma\|^2_{H^1}
+\|\sigma\|_{H^1}\|\theta\|_{H^1}+\|\rho\|_{L^\infty}\|\theta\|^2_{H^1}
+\|\mathbf{u}\|_{H^3}\|\theta\|_{H^1})ds\\
&\leq C_0(M_0)+Ct\|\mathbf{u}\|_{L^\infty_t(H^1)}(\|\sigma\|^2_{L^\infty_t(H^1)}
+\|\sigma\|_{L^\infty_t(H^1)}\|\theta\|_{L^\infty_t(H^1)}
+\|\rho\|_{L^\infty_{x,t}}\|\theta\|^2_{L^\infty_t(H^1)})\\
&\quad+Ct^{\frac{1}{2}}\|\mathbf{u}\|_{L^2_t(H^3)}\|\mathbf{u}\|_{L^\infty_t(H^1)}
\|\theta\|_{L^\infty_t(H^1)}\\
&\leq  C_0(M_0){\rm exp}\{t^\frac 1 2 C_1(M(t)) \}.
\end{split}\end{equation}
From \ef{lemma1} and \ef{basic} and Lemma \ref{hs1}, we conclude \ef{basicestimate} immediately.
\hfill$\Box$

\begin{lem} Assume that the assumptions in Lemma \ref{lemrho} are satisfied. Then the solution to the system ${\eqref{ns0}}$ satisfies the following estimate
\begin{equation}\begin{split}\label{basicestimate1}
\|(\nabla\sigma,{\rm div}\mathbf{u},\nabla\theta,\nabla^2\mathbf{H})\|^2_{L^2}(t)
+\|\nabla{\rm div}\mathbf{u}\|^2_{L^2_t(L^2)}+\|\Delta\theta\|^2_{L^2_t(L^2)}\leq  C_0(M_0){\rm exp}\{t^\frac 1 2 C_1(M(t)) \},
\end{split}\end{equation}
for $\epsilon \in (0,1]$ and  $t\in[0,T]$.
\end{lem}
\textbf{Proof}. By virtue of the boundary condition $\mathbf{u}\cdot\mathbf{n}|_{\partial\Omega}=0$, multiplying $\nabla{\eqref{ns0}}_1$ with $\nabla\sigma$ and integrating the resulting equality in $L^2(\Omega)$, we get
\begin{equation*}\begin{split}
&\frac{1}{2}\frac{d}{dt}\|\nabla\sigma\|^2_{L^2}+\frac{1}{\epsilon}\int_\Omega\nabla{\rm div}\mathbf{u}\cdot\nabla\sigma dx\\
&\quad=-\int_\Omega({\bl{\rm div}\mathbf{u}\nabla\sigma}+\sigma\nabla{\rm div}\mathbf{u}+\nabla\mathbf{u}\cdot\nabla\sigma+\mathbf{u}\cdot\nabla^2\sigma)\cdot\nabla\sigma dx\\
&\quad{=-\int_\Omega(\frac{1}{2}{\bl{\rm div}\mathbf{u}\nabla\sigma}+\sigma\nabla{\rm div}\mathbf{u}+\nabla\mathbf{u}\cdot\nabla\sigma)\cdot\nabla\sigma dx}\\
&\quad{\lesssim\|\mathbf{u}\|_{H^3}\|\sigma\|^2_{H^1}}.
\end{split}\end{equation*}
Then we take the inner product of $\langle{\eqref{ns0}}_2,-\nabla{\rm div}\mathbf{u}\rangle$ and $\langle{\nabla\eqref{ns0}}_4,\nabla\theta\rangle$ in $L^2(\Omega)$, respectively, to derive that
\begin{equation*}\begin{split}
&\frac{1}{2}\frac{d}{dt} \|\sqrt{\rho}{\rm div}\mathbf{u}\|^2_{L^2} +(\mu+\nu)\|\nabla{\rm div}\mathbf{u}\|^2_{L^2} -\frac{1}{\epsilon}\int_\Omega\nabla(\sigma+\theta)\cdot\nabla{\rm div}\mathbf{u} dx\\
&\quad=\int_\Omega[\nabla{\rm div}\mathbf{u}\cdot\nabla(\sigma\theta)+\frac{1}{2}\rho_t({\rm div}\mathbf{u})^2]dx
-\int_\Omega{\bl\nabla\rho\cdot\mathbf{u}_t{\rm div}\mathbf{u}}dx\\
&\quad\quad+\int_\Omega[\rho\mathbf{u}\cdot\nabla\mathbf{u}
+(\mathbf{H}\cdot\nabla)\mathbf{H}-\frac{1}{2}\nabla(|\mathbf{H}|^2)]\cdot
\nabla{\rm div}\mathbf{u}dx\\
&\quad\lesssim\|\mathbf{u}\|_{H^3}\|\sigma\|_{H^1}\|\theta\|_{H^1}+\|\mathbf{u}\|_{H^2}
(\|\epsilon\sigma_t\|_{H^1}\|\mathbf{u}\|_{H^2}
+\|\sigma\|_{H^1}\|\epsilon\mathbf{u}_t\|_{H^1}\\
&\quad\quad+\|\rho\|_{L^\infty}\|\mathbf{u}\|_{H^1}\|\mathbf{u}\|_{H^2}
+\|\mathbf{H}\|_{H^1}\|\mathbf{H}\|_{H^2}),
\end{split}\end{equation*}
and
\begin{equation*}\begin{split}
&\frac{1}{2}\frac{d}{dt}\|\sqrt{\rho}\nabla\theta\|^2_{L^2}
+\kappa\|\Delta\theta\|^2_{L^2}+
\frac{1}{\epsilon}\int_\Omega\nabla{\rm div}\mathbf{u}\cdot\nabla\theta dx\\
&\quad=-\int_\Omega[{\bl(\theta_t+\mathbf{u}\cdot\nabla\theta)\nabla\rho}
+\rho\nabla\mathbf{u}\cdot\nabla\theta+
{\bl{\rm div}\mathbf{u}\nabla(\rho\theta+\sigma)}\\
&\quad\quad+(\rho\theta+\sigma)\nabla{\rm div}\mathbf{u}-\epsilon\nabla(2\mu|D(\mathbf{u})|^2+\xi({\rm div}\mathbf{u})^2)]\cdot\nabla\theta dx\\
&\quad\lesssim\|\theta\|_{H^2}(\|\sigma\|_{H^1}\|\epsilon\theta_t\|_{H^1}+
\|\rho\|_{H^2}\|\mathbf{u}\|_{H^2}\|\theta\|_{H^1}+\|\mathbf{u}\|_{H^2}\|\sigma\|_{H^1}+\epsilon\|\mathbf{u}\|^2_{H^2}).
\end{split}\end{equation*}
Next, we take $\langle\partial_{ij}{\eqref{ns0}}_3,\partial_{ij}\mathbf{H}\rangle$ to get
\begin{equation*}\begin{split}
\frac{1}{2}\frac{d}{dt}\|\partial_{ij}\mathbf{H}\|^2_{L^2}=&\frac 1 2 \int_\Omega {\rm div}\mathbf{u}
 |\partial_{ij}\mathbf{H}|^2 dx-\int_\Omega(\partial_{ij}\mathbf{u}
\cdot\nabla\mathbf{H}-\partial_i\mathbf{u}\cdot\nabla\partial_j\mathbf{H}
-\partial_j\mathbf{u}\cdot\nabla\partial_i\mathbf{H}){\bl\cdot\partial_{ij}\mathbf{H}}dx\\
&\quad-\int_\Omega\partial_{ij}(\mathbf{H}
{\rm div}\mathbf{u}-\mathbf{H}\cdot\nabla\mathbf{u}){\bl\cdot\partial_{ij}\mathbf{H}}dx\\
&\lesssim\|\mathbf{H}\|^2_{H^2}\|\mathbf{u}\|_{H^3}.
\end{split}\end{equation*}
Putting the above estimates together and integrate {the resulting} inequality on $[0,t]$, we obtain
\begin{equation}\label{basic2}\begin{split}
&\|\nabla\sigma\|^2_{L^2}(t)+\|\sqrt{\rho}{\rm div}\mathbf{u}\|^2_{L^2}(t)
+\|\sqrt{\rho}\nabla\theta\|^2_{L^2}(t)+\|\nabla^2\mathbf{H}\|^2_{L^2}(t)\\
&\quad+(\mu+\nu)\|\nabla{\rm div}\mathbf{u}\|^2_{L^2_t(L^2)}
+\|\sqrt{\kappa}\Delta\theta\|^2_{L^2_t(L^2)}\\
&\leq C_0(M_0)+
C\int_0^t[\|\mathbf{u}\|_{H^3}\|\sigma\|_{H^1}(\|\sigma\|_{H^1}+\|\theta\|_{H^1})
+\|\mathbf{u}\|_{H^2}(\|\epsilon\sigma_t\|_{H^1}\|\mathbf{u}\|_{H^2}\\
&\quad
+\|\sigma\|_{H^1}\|\epsilon\mathbf{u}_t\|_{H^1}
+\|\rho\|_{L^\infty}\|\mathbf{u}\|_{H^1}\|\mathbf{u}\|_{H^2}+\|\mathbf{H}\|_{H^1}\|\mathbf{H}\|_{H^2})+\|\theta\|_{H^2}(\|\rho\|_{H^2}\|\mathbf{u}\|_{H^2}\|\theta\|_{H^1}\\
&\quad+\|\sigma\|_{H^1}\|\epsilon\theta_t\|_{H^1}+\|\mathbf{u}\|_{H^2}\|\sigma\|_{H^1}+\epsilon\|\mathbf{u}\|^2_{H^2})
+\|\mathbf{H}\|^2_{H^2}\|\mathbf{u}\|_{H^3}]ds\\
&\leq C_0(M_0)+Ct^{\frac{1}{2}}[\|\mathbf{u}\|_{L^2_t(H^3)}\|\sigma\|_{L^\infty_t(H^1)}
(\|\sigma\|_{L^\infty_t(H^1)}+\|\theta\|_{L^\infty_t(H^1)})\\
&\quad
+\|\mathbf{u}\|_{L^\infty_t(H^2)}(\|\epsilon\sigma_t\|_{L^\infty_t(H^1)}\|\mathbf{u}\|_{L^\infty_t(H^2)}
+\|\sigma\|_{L^\infty_t(H^1)}\|\epsilon\mathbf{u}_t\|_{L^\infty_t(H^1)}\\
&\quad+\|\rho\|_{L^\infty_{x,t}}\|\mathbf{u}\|_{L^\infty_t(H^1)}\|\mathbf{u}\|_{L^\infty_t(H^2)}
+\|\mathbf{H}\|_{L^\infty_t(H^1)}\|\mathbf{H}\|_{L^\infty_t(H^2)})\\
&\quad
+\|\theta\|_{L^2_t(H^2)}(\|\rho\|_{L^\infty_t(H^2)}\|\mathbf{u}\|_{L^\infty_t(H^2)}\|\theta\|_{L^\infty_t(H^1)}
+\|\sigma\|_{L^\infty_t(H^1)}\|\epsilon\theta_t\|_{L^\infty_t(H^1)}\\
&\quad+\|\mathbf{u}\|_{L^\infty_t(H^2)}\|\sigma\|_{L^\infty_t(H^1)}+\epsilon\|\mathbf{u}\|^2_{L^\infty_t(H^2)})
+\|\mathbf{H}\|^2_{L^\infty_t(H^2)}\|\mathbf{u}\|_{L^2_t(H^3)}]\\
&\leq  C_0(M_0){\rm exp}\{t^\frac 1 2 C_1(M(t)) \}.
\end{split}\end{equation}
From \ef{lemma1} and \ef{basic2}, we conclude \ef{basicestimate1} immediately.
\hfill$\Box$

\subsection{Estimates for $\mathbf{u}$, $\e \sigma_t$, $\e \theta_t$, and $\e \mathbf{u}_t$}

\begin{lem}Let the assumptions in Lemma \ref{lemrho} be satisfied. Then for $\epsilon \in (0,1]$ and  $t\in[0,T]$ we have
\begin{equation}\label{firstorder1}\begin{split}
&\|\epsilon(\sigma_t,\mathbf{u}_t,\mathbf{H}_t,\theta_t)\|^2_{L^2}(t)
+\|\epsilon\mathbf{u}_t\|^2_{L^2_t(H^1)} +\|\epsilon\nabla\theta_t\|^2_{L^2_t(L^2)}\leq C_0(M_0) {\rm exp} \{t C_1(M(t))\},
\end{split}\end{equation}
and
\begin{equation}\label{firstorder2}\begin{split}
&\|\epsilon(\nabla\sigma_t,{\rm div}\mathbf{u}_t,\nabla\theta_t)\|^2_{L^2}(t)+\|\epsilon\nabla{\rm div}\mathbf{u}_t\|^2_{L^2_t(L^2)}+\|\epsilon\Delta\theta_t\|^2_{L^2_t(L^2)}\leq C_0(M_0) {\rm exp} \{t^{\frac{1}{4}} C_1(M(t))\}.
\end{split}\end{equation}
\end{lem}
\textbf{Proof}. Differentiating ${\eqref{ns0}}$ with respect to $t$, we deduce that
\begin{equation}\label{nstime}\left\{\begin{split}
&\sigma_{tt}+\frac{1}{\epsilon}{\rm div}\mathbf{u}_t=-(\mathbf{u}_t\cdot\nabla\sigma+\mathbf{u}\cdot\nabla\sigma_t+\sigma_t\cdot{\rm div}\mathbf{u}+\sigma{\rm div}\mathbf{u}_t),\quad &(t,x)\in(0,T]\times\Omega,\\
&\rho(\mathbf{u}_{tt}+\mathbf{u}\cdot\nabla\mathbf{u}_t)+\frac{1}{\epsilon}
(\nabla\sigma_t+\nabla\theta_t)-\mu\Delta\mathbf{u}_t
-\nu\nabla{\rm div}\mathbf{u}_t\\
&\quad=-\rho_t(\mathbf{u}_t+\mathbf{u}\cdot\nabla\mathbf{u})-\rho\mathbf{u}_t\cdot\nabla\mathbf{u}
-\nabla(\sigma\theta)_t+{\rm curl}\mathbf{H}_t\times\mathbf{H}+{\rm curl}\mathbf{H}\times\mathbf{H}_t,\quad &(t,x)\in(0,T]\times\Omega,\\
&\mathbf{H}_{tt}={\rm div}\mathbf{u}_t\mathbf{H}-{\rm div}\mathbf{u}\mathbf{H}_t-\mathbf{u}_t\cdot\nabla\mathbf{H}
-\mathbf{u}\cdot\nabla\mathbf{H}_t+\mathbf{H}_t\cdot\nabla\mathbf{u}+\mathbf{H}\cdot\nabla\mathbf{u}_t,
\quad &(t,x)\in(0,T]\times\Omega,\\
&\rho(\theta_{tt}+\mathbf{u}\cdot\nabla\theta_t)+\frac{1}{\epsilon}{\rm div}\mathbf{u}_t-\kappa\Delta\theta_t\\
&\quad=\epsilon(2\mu|D(\mathbf{u})|^2+\xi({\rm div}\mathbf{u})^2)_t-((\rho\theta+\sigma){\rm div}\mathbf{u})_t-\rho_t\theta_t-(\rho\mathbf{u})_t\cdot\nabla\theta,\quad &(t,x)\in(0,T]\times\Omega,\\
&\mathbf{u}_t\cdot\mathbf{n}=0,\ \mathbf{n}\times{\rm curl}\mathbf{u}_t=0,\ \frac{\partial\theta_t}{\partial\mathbf{n}}=0, \quad&{\rm on} \quad(0,T]\times\partial\Omega,\\
&(\sigma_t,\mathbf{u}_t,\mathbf{H}_t,\theta_t)(0,x)=(\sigma_t(0),\mathbf{u}_t(0),\mathbf{H}_t(0),\theta_t(0)),\quad &x\in\Omega,
\end{split}\right.\end{equation}
where
$$
\sigma_t(0):=-\frac{1}{\epsilon}{\rm div}\mathbf{u}_0-\sigma_0{\rm div}\mathbf{u}_0-\mathbf{u}_0\cdot\nabla\sigma_0,
$$
$$
\mathbf{u}_t(0):=[\mu\Delta\mathbf{u}_0+\nu\nabla{\rm div}\mathbf{u}_0+{\rm curl}\mathbf{H}_0\cdot\mathbf{H}_0 -\nabla(\sigma_0\theta_0)-\frac{1}{\epsilon}(\nabla\sigma_0+\nabla\theta_0)]/\rho_0
-\mathbf{u}_0\cdot\nabla\mathbf{u}_0,
$$
$$
\theta_t(0):=[\kappa\Delta\theta_0+\epsilon(2\mu|D(\mathbf{u}_0)|^2+\xi({\rm div}\mathbf{u}_0)^2)-\frac{1}{\epsilon}{\rm div}\mathbf{u}_0-(\rho_0\theta_0+\sigma_0){\rm div}\mathbf{u}_0]/\rho_0
-\mathbf{u}_0\cdot\nabla\theta_0,
$$
$$
\mathbf{H}_t(0):=(\mathbf{H}_0\cdot\nabla)\mathbf{u}_0-(\mathbf{u}_0\cdot\nabla)\mathbf{H}_0-\mathbf{H}_0{\rm div}\mathbf{u}_0.
$$
Taking the inner product of ${\eqref{nstime}}_1$ and $\epsilon^2\sigma_t$ in $L^2(\Omega)$ and integrating by parts, we get
\begin{equation*}\begin{split}
\frac{1}{2}\frac{d}{dt}\|\epsilon\sigma_t\|^2_{L^2}+\epsilon\int_\Omega{\rm div}\mathbf{u}_t\cdot\sigma_tdx&=-\int_\Omega\epsilon^2(\mathbf{u}_t\cdot\nabla\sigma+\mathbf{u}\cdot\nabla\sigma_t+\sigma_t\cdot{\rm div}\mathbf{u}+\sigma{\rm div}\mathbf{u}_t)\sigma_tdx\\
&\lesssim\|\epsilon\sigma_t\|_{H^1}(\|\epsilon\mathbf{u}_t\|_{H^1}\|\sigma\|_{H^1}
+\|\mathbf{u}\|_{H^2}\|\epsilon\sigma_t\|_{H^1}).
\end{split}\end{equation*}
Then we multiply ${\eqref{nstime}}_2$ by $\epsilon^2\mathbf{u}_t$ in $L^2(\Omega)$ to find
\begin{equation*}\begin{split}
\frac{1}{2}&\frac{d}{dt}\|\epsilon\sqrt{\rho}\mathbf{u}_t\|^2_{L^2}+(\mu+\nu)\|\epsilon{\rm div}\mathbf{u}_t\|^2_{L^2}+\mu\|\epsilon{\rm curl}\mathbf{u}_t\|^2_{L^2}
+\epsilon\int_\Omega\nabla(\sigma_t+\theta_t)\cdot\mathbf{u}_tdx\\
&=-\int_\Omega\epsilon^2[\rho_t\mathbf{u}_t+(\rho\mathbf{u})_t\cdot\nabla\mathbf{u}
+\nabla(\sigma\theta)_t+\frac{1}{2}\nabla(|\mathbf{H}|^2)_t-\mathbf{H}_t\cdot\nabla\mathbf{H}
-\mathbf{H}\cdot\nabla\mathbf{H}_t]\cdot\mathbf{u}_tdx\\
&\lesssim\|\epsilon\mathbf{u}_t\|_{H^1}[\|\epsilon\sigma_t\|_{H^1}(\|\epsilon\mathbf{u}_t\|_{H^1}+\epsilon\|\mathbf{u}\|^2_{H^2})
+\|\epsilon\mathbf{H}_t\|_{L^2}\|\mathbf{H}\|_{H^2}\\
&\quad+\|\rho\|_{L^\infty}\|\epsilon\mathbf{u}_t\|_{H^1}\|\mathbf{u}\|_{H^2}
+\|\epsilon\sigma_t\|_{H^1}\|\theta\|_{H^1}+\|\epsilon\theta_t\|_{H^1}\|\sigma\|_{H^1}],
\end{split}\end{equation*}
where the integration by parts is employed.
We perform $\langle{\eqref{nstime}}_4,\epsilon^2\theta_t\rangle$ in $L^2(\Omega)$ to derive
\begin{equation*}\begin{split}
&\frac{1}{2}\frac{d}{dt}\|\epsilon\sqrt{\rho}\theta_t\|^2_{L^2}
+\kappa\|\epsilon\nabla\theta_t\|^2_{L^2}
+\epsilon\int_\Omega{\rm div}\mathbf{u}_t\cdot\theta_tdx\\
&=\int_\Omega\epsilon^2[\epsilon(2\mu|D(\mathbf{u})|^2+\xi({\rm div}\mathbf{u})^2)_t-(\rho_t\theta_t+(\rho\mathbf{u})_t\cdot\nabla\theta)
-((\rho\theta+\sigma){\rm div}\mathbf{u})_t]\cdot\theta_tdx\\
&\lesssim\|\epsilon\theta_t\|_{H^1}[\epsilon\|\epsilon\mathbf{u}_t\|_{H^1}\|\mathbf{u}\|_{H^2}+
(\|\epsilon\sigma_t\|_{H^1}\|\theta\|_{H^1}+\|\rho\|_{L^\infty}\|\epsilon\theta_t\|_{H^1}+\|\epsilon\sigma_t\|_{H^1})
\|\mathbf{u}\|_{H^2}\\
&\quad+(\|\rho\|_{L^\infty}\|\theta\|_{H^1}+\|\sigma\|_{H^1})\|\epsilon\mathbf{u}_t\|_{H^1}
+\|\epsilon\sigma_t\|_{H^1}\|\epsilon\theta_t\|_{H^1}].
\end{split}\end{equation*}
Taking the inner product of $\langle{\eqref{nstime}}_3,\epsilon^2\mathbf{H}_t\rangle$ in $L^2(\Omega)$ to deduce
\begin{equation*}\begin{split}
\frac{1}{2}\frac{d}{dt}\|\epsilon\mathbf{H}_t\|^2_{L^2}=&=-\epsilon^2\int_\Omega({\rm div}\mathbf{u}_t\mathbf{H}+ \frac 1 2 {\rm div}\mathbf{u}\mathbf{H}_t+\mathbf{u}_t\cdot\nabla\mathbf{H}-\mathbf{H}_t\cdot\nabla\mathbf{u}-\mathbf{H}\cdot\nabla\mathbf{u}_t)\cdot\mathbf{H}_tdx\\
&\lesssim\|\epsilon\mathbf{H}_t\|_{L^2}(\|\epsilon\mathbf{u}_t\|_{H^1}\|\mathbf{H}\|_{H^2}+
\|\epsilon\mathbf{H}_t\|_{L^2}\|\mathbf{u}\|_{H^3}).
\end{split}\end{equation*}
Finally, we put all the estimates together and integrate the resulting inequality to obtain
\begin{equation*}\begin{split}
&\|\epsilon\sigma_t\|^2_{L^2}(t)+\|\epsilon\sqrt{\rho}\mathbf{u}_t\|^2_{L^2}(t)
+\|\epsilon\sqrt{\rho}\theta_t\|^2_{L^2}(t)+\|\epsilon\mathbf{H}_t\|^2_{L^2} (t)\\
&\quad\quad+(\mu+\nu)\|\epsilon{\rm div}\mathbf{u}_t\|^2_{L^2_t(L^2)}+\mu\|\epsilon{\rm curl}\mathbf{u}_t\|^2_{L^2_t(L^2)}+\kappa\|\epsilon\nabla\theta_t\|^2_{L^2_t(L^2)}\\
&\leq C_0(M_0)+ Ct\{\|\epsilon\sigma_t\|_{L^\infty_t(H^1)}(\|\mathbf{u}\|_{L^\infty_t(H^2)}\|\epsilon\sigma_t\|_{L^\infty_t(H^1)}+
\|\epsilon\mathbf{u}_t\|_{L^\infty_t(H^1)}\|\sigma\|_{L^\infty_t(H^1)})\\
&\quad+\|\epsilon\mathbf{u}_t\|_{L^\infty_t(H^1)}
[\|\epsilon\sigma_t\|_{L^\infty_t(H^1)}(\|\epsilon\mathbf{u}_t\|_{L^\infty_t(H^1)}+\epsilon\|\mathbf{u}\|^2_{L^\infty_t(H^2)})
+\|\epsilon\mathbf{H}_t\|_{L^\infty_t(L^2)}\|\mathbf{H}\|_{L^\infty_t(H^2)}\\
&\quad+\|\rho\|_{L^\infty_{x,t}}\|\epsilon\mathbf{u}_t\|_{L^\infty_t(H^1)}\|\mathbf{u}\|_{L^\infty_t(H^2)}
+\|\epsilon\sigma_t\|_{L^\infty_t(H^1)}\|\theta\|_{L^\infty_t(H^1)}
+\|\epsilon\theta_t\|_{L^\infty_t(H^1)}\|\sigma\|_{L^\infty_t(H^1)}]\\
&\quad+\|\epsilon\theta_t\|_{L^\infty_t(H^1)}[\epsilon\|\epsilon\mathbf{u}_t\|_{L^\infty_t(H^1)}
\|\mathbf{u}\|_{L^\infty_t(H^2)}+\|\epsilon\sigma_t\|_{L^\infty_t(H^1)}\|\epsilon\theta_t\|_{L^\infty_t(H^1)}\\
&\quad
+(\|\epsilon\sigma_t\|_{L^\infty_t(H^1)}\|\theta\|_{L^\infty_t(H^1)}+\|\rho\|_{L^\infty_{x,t}}\|\epsilon\theta_t\|_{L^\infty_t(H^1)}
+\|\epsilon\sigma_t\|_{L^\infty_t(H^1)})\|\mathbf{u}\|_{L^\infty_t(H^2)}\\
&\quad+(\|\rho\|_{L^\infty_{x,t}}\|\theta\|_{L^\infty_t(H^1)}+\|\sigma\|_{L^\infty_t(H^1)})
\|\epsilon\mathbf{u}_t\|_{L^\infty_t(H^1)}
]\}\\
&\quad+Ct^\frac{1}{2}\|\epsilon\mathbf{H}_t\|_{L^\infty_t(L^2)}
(\|\epsilon\mathbf{u}_t\|_{L^\infty_t(H^1)}\|\mathbf{H}\|_{L^\infty_t(H^2)}+
\|\epsilon\mathbf{H}_t\|_{L^\infty_t(L^2)}\|\mathbf{u}\|_{L^2_t(H^3)})\\
&\leq C_0(M_0) {\rm exp}\{t^\frac{1}{2} C_1(M(t))\},\ t\in[0,T],\ \epsilon\in(0,1].
\end{split}\end{equation*}
From the above inequality and Lemma \ref{hs1}, we obtain \ef{firstorder1}.

Next, we take the inner product of $\nabla{\eqref{nstime}}_1$ and $\epsilon^2\nabla\sigma_t$ in $L^2(\Omega)$ to get
\begin{equation*}\begin{split}
\frac{1}{2}\frac{d}{dt}&\|\epsilon\nabla\sigma_t\|^2_{L^2}+\epsilon\int_\Omega\nabla{\rm div}\mathbf{u}_t\cdot\nabla\sigma_tdx\\
&=-\epsilon^2\int_\Omega(\nabla\sigma_t{\rm div}\mathbf{u}+\sigma_t\nabla{\rm div}\mathbf{u}
+\sigma\nabla{\rm div}\mathbf{u}_t+\nabla\sigma{\rm div}\mathbf{u}_t\\
&\quad+\nabla\mathbf{u}_t\cdot\nabla\sigma
+\mathbf{u}_t\cdot\nabla^2\sigma+\nabla\mathbf{u}\cdot\nabla\sigma_t+\mathbf{u}\cdot\nabla^2\sigma_t)\cdot\nabla\sigma_tdx\\
&\lesssim{\|\mathbf{u}\|_{H^3}\|\epsilon\sigma_t\|^2_{H^1}}
+\|\epsilon\nabla\sigma_t\|_{L^3}(\|\sigma\|_{H^1}\|\epsilon\mathbf{u}_t\|_{H^2}
+\|\sigma\|_{H^2}\|\epsilon\mathbf{u}_t\|_{H^1}),
\end{split}\end{equation*}
 where we have used that $$\epsilon^2\int_\Omega\mathbf{u}_t\cdot\nabla^2\sigma\cdot\nabla\sigma_tdx
=-\epsilon^2\int_\Omega\nabla\sigma\cdot(\nabla\sigma_t{\rm div}\mathbf{u}_t+\mathbf{u}_t\cdot\nabla^2\sigma_t)dx.$$  We also perform $\langle{\eqref{nstime}}_2,-\epsilon^2\nabla{\rm div}\mathbf{u}_t\rangle$ in $L^2(\Omega)$ to derive that
\begin{equation*}\begin{split}
&\frac{1}{2}\frac{d}{dt}\|\epsilon\sqrt{\rho}{\rm div}\mathbf{u}_t\|^2_{L^2}+(\mu+\nu)\|\epsilon\nabla{\rm div}\mathbf{u}_t\|^2_{L^2}-\epsilon\int_\Omega\nabla(\sigma_t+\theta_t)\cdot\nabla{\rm div}\mathbf{u}_tdx\\
&=\epsilon^2\int_\Omega(\rho_t\mathbf{u}_t+\rho_t\mathbf{u}\cdot\nabla\mathbf{u}+\rho\mathbf{u}_t\cdot\nabla\mathbf{u}
+\rho\mathbf{u}\cdot\nabla\mathbf{u}_t+\sigma_t\nabla\theta+\sigma\nabla\theta_t
+\theta_t\nabla\sigma+\theta\nabla\sigma_t)\cdot\nabla{\rm div}\mathbf{u}_tdx\\
&\quad+\epsilon^2\int_\Omega\frac{1}{2}\rho_t({\rm div}\mathbf{u}_t)^2dx-\epsilon^2\int_\Omega\nabla\rho\mathbf{u}_{tt}\cdot{\rm div}\mathbf{u}_tdx-\epsilon^2\int_\Omega({\rm curl}\mathbf{H}_t\times\mathbf{H}+{\rm curl}\mathbf{H}\times\mathbf{H}_t)\cdot\nabla{\rm div}\mathbf{u}_tdx\\
&\lesssim\|\epsilon\mathbf{u}_t\|_{H^2}[\|\epsilon\sigma_t\|_{H^1}(\|\epsilon\mathbf{u}_t\|_{H^1}+\epsilon\|\mathbf{u}\|^2_{H^2})
+\|\rho\|_{L^\infty}\|\epsilon\mathbf{u}_t\|_{H^1}\|\mathbf{u}\|_{H^2}+\|\sigma\|_{H^1}\|\epsilon^2\mathbf{u}_{tt}\|_{L^3}\\
&\quad
+\|\epsilon\sigma_t\|_{H^1}\|\nabla\theta\|_{L^3}
+\|\theta\|_{H^1}\|\epsilon\nabla\sigma_t\|_{L^3}
+\|\epsilon\theta_t\|_{H^1}\|\nabla\sigma\|_{L^3}
+\|\sigma\|_{H^1}\|\epsilon\nabla\theta_t\|_{L^3}
+\|\mathbf{H}\|_{H^2}\|\epsilon\mathbf{H}_t\|_{H^1}],
\end{split}\end{equation*}
Then we multiply $\nabla{\eqref{nstime}}_3$ and $\epsilon^2\nabla\theta_t$ in $L^2(\Omega)$ to find
\begin{equation*}\begin{split}
&\frac{1}{2}\frac{d}{dt}\|\epsilon\sqrt{\rho}\nabla\theta_t\|^2_{L^2}
+\kappa\|\epsilon\Delta\theta_t\|^2_{L^2}
+\epsilon\int_\Omega\nabla{\rm div}\mathbf{u}_t\cdot\nabla\theta_tdx\\
&={ -\epsilon^2\int_\Omega[\nabla\rho(\theta_{tt}+\mathbf{u}\cdot\nabla\theta_t)+\rho\nabla\mathbf{u}\cdot\nabla\theta_t]\cdot\nabla\theta_tdx}\\
&\quad+\epsilon^2\int_\Omega\nabla[\epsilon(2\mu|D(\mathbf{u})|^2+\xi({\rm div}\mathbf{u})^2)_t-((\rho\theta+\sigma){\rm div}\mathbf{u})_t-\rho_t\theta_t-(\rho\mathbf{u})_t\cdot\nabla\theta]\cdot\nabla\theta_tdx\\
&\lesssim\|\epsilon\theta_t\|_{H^1}[\|\rho\|_{H^2}
\|\mathbf{u}\|_{H^2}\|\epsilon\theta_t\|_{H^2}
+\|\epsilon\mathbf{u}\|_{H^3}\|\epsilon\mathbf{u}_t\|_{H^2}+\|\epsilon\sigma_t\|_{H^1}\|\epsilon\theta_t\|_{H^2}\\
&\quad+\|\mathbf{u}\|_{H^3}(\|\epsilon\sigma_t\|_{H^1}\|\theta\|_{H^1}+\|\rho\|_{H^2}\|\epsilon\theta_t\|_{H^1}
+\|\epsilon\sigma_t\|_{H^1})]\\
&\quad+\|\epsilon\nabla\theta_t\|_{L^3}
[\|\sigma\|_{H^1}\|\epsilon^2\theta_{tt}\|_{H^1}
+\epsilon\|\mathbf{u}\|_{H^3}
\|\epsilon\sigma_t\|_{H^1}\|\theta\|_{H^1}\\
&\quad+\|\epsilon\mathbf{u}_t\|_{H^2}(\|\rho\|_{L^\infty}\|\theta\|_{H^1}
+\|\sigma\|_{H^1})+\|\theta\|_{H^2}
(\|\epsilon\sigma_t\|_{H^1}\|\mathbf{u}\|_{H^2}+\|\rho\|_{H^2}\|\epsilon\mathbf{u}_t\|_{H^1})].
\end{split}\end{equation*}
Finally we put all the estimates together and integrate the resulting inequality on $[0,t]$ to obtain
\begin{equation*}\begin{split}
&\|\epsilon\nabla\sigma_t\|^2_{L^2}(t)+\|\epsilon\sqrt{\rho}{\rm div}\mathbf{u}_t\|^2_{L^2}(t)+\|\epsilon\sqrt{\rho}\nabla\theta_t\|^2_{L^2}(t)+(\mu+\nu)\|\epsilon\nabla{\rm div}\mathbf{u}_t\|^2_{L^2_t(L^2)}+\kappa\|\epsilon\Delta\theta_t\|^2_{L^2_t(L^2)}\\
&\leq
C_0(M_0)+C{\bl t^{\frac{1}{2}}}\{\|\mathbf{u}\|_{L^2_t(H^3)}\|\epsilon\sigma_t\|^2_{L^\infty_t(H^1)}
+\|\epsilon\mathbf{u}_t\|_{L^2_t(H^2)}[\|\epsilon\sigma_t\|_{L^\infty_t(H^1)}
(\|\epsilon\mathbf{u}_t\|_{L^\infty_t(H^1)}\\
&\quad+\epsilon\|\mathbf{u}\|^2_{L^\infty_t(H^2)})
+\|\rho\|_{L^\infty_{x,t}}\|\epsilon\mathbf{u}_t\|_{L^\infty_t(H^1)}\|\mathbf{u}\|_{L^\infty_t(H^2)}
+\|\mathbf{H}\|_{L^\infty_t(H^2)}\|\epsilon\mathbf{H}_t\|_{L^\infty_t(H^1)}]\\
&\quad+\|\epsilon\theta_t\|_{L^\infty_t(H^1)}[\|\rho\|_{L^\infty_t(H^2)}
\|\mathbf{u}\|_{L^\infty_t(H^2)}\|\epsilon\theta_t\|_{L^2_t(H^2)}
+\|\epsilon\mathbf{u}\|_{L^\infty_t(H^3)}\|\epsilon\mathbf{u}_t\|_{L^2_t(H^2)}\\
&\quad +\|\mathbf{u}\|_{L^2_t(H^3)}({\bl\|\epsilon\sigma_t\|_{L^\infty_t(H^1)}
\|\theta\|_{L^\infty_t(H^1)}}
+\|\rho\|_{L^\infty_t(H^2)}\|\epsilon\theta_t\|_{L^\infty_t(H^1)}\\
&\quad+\|\epsilon\sigma_t\|_{L^\infty_t(H^1)})
{\bl+\|\epsilon\sigma_t\|_{L^\infty_t(H^1)}
\|\epsilon\theta_t\|_{L^2_t(H^2)}}]\}\\
&\quad+Ct^{\frac{1}{4}}\{\|\epsilon\sigma_t\|^{\frac{1}{2}}_{L^\infty_t(H^1)}
\|\epsilon\sigma_t\|^{\frac{1}{2}}_{L^2_t(H^2)}(\|\sigma\|_{L^\infty_t(H^1)}\|\epsilon\mathbf{u}_t\|_{L^2_t(H^2)}
+\|\sigma\|_{L^2_t(H^2)}\|\epsilon\mathbf{u}_t\|_{L^\infty_t(H^1)})\\
&\quad+\|\epsilon\mathbf{u}_t\|_{L^2_t(H^2)}(\|\sigma\|_{L^\infty_t(H^1)}\|\epsilon^2\mathbf{u}_{tt}\|^{\frac{1}{2}}_{L^2_t(H^1)}
\|\epsilon^2\mathbf{u}_{tt}\|^{\frac{1}{2}}_{L^\infty_t(L^2)}\\
&\quad
+\|\epsilon\sigma_t\|_{L^\infty_t(H^1)}\|\theta\|^{\frac{1}{2}}_{L^2_t(H^2)}\|\theta\|^{\frac{1}{2}}_{L^\infty_t(H^1)}
+\|\theta\|_{L^\infty_t(H^1)}\|\epsilon\sigma_t\|^{\frac{1}{2}}_{L^2_t(H^2)}
\|\epsilon\sigma_t\|^{\frac{1}{2}}_{L^\infty_t(H^1)}\\
&\quad
+\|\epsilon\theta_t\|_{L^\infty_t(H^1)}\|\sigma\|^{\frac{1}{2}}_{L^2_t(H^2)}\|\sigma\|^{\frac{1}{2}}_{L^\infty_t(H^1)}
+\|\sigma\|_{L^\infty_t(H^1)}\|\epsilon\theta_t\|^{\frac{1}{2}}_{L^2_t(H^2)}
\|\epsilon\theta_t\|^{\frac{1}{2}}_{L^\infty_t(H^1)})\\
&\quad+\|\epsilon\theta_t\|^{\frac{1}{2}}_{L^2_t(H^2)}\|\epsilon\theta_t\|^{\frac{1}{2}}_{L^\infty_t(H^1)}
[\|\sigma\|_{L^\infty_t(H^1)}\|\epsilon^2\theta_{tt}\|_{{L^2_t(H^1)}}\\
&\quad+\epsilon\|\mathbf{u}\|_{L^2_t(H^3)}
\|\epsilon\sigma_t\|_{L^\infty_t(H^1)}\|\theta\|_{L^\infty_t(H^1)}+\|\epsilon\mathbf{u}_t\|_{L^2_t(H^2)}
(\|\rho\|_{L^\infty_{x,t}}\|\theta\|_{L^\infty_t(H^1)}
+\|\sigma\|_{L^\infty_t(H^1)})\\
&\quad+\|\theta\|_{L^2_t(H^2)}
(\|\epsilon\sigma_t\|_{L^\infty_t(H^1)}\|\mathbf{u}\|_{L^\infty_t(H^2)}
+\|\rho\|_{L^\infty_t(H^2)}\|\epsilon\mathbf{u}_t\|_{L^\infty_t(H^1)})]\}\\
&\leq C_0(M_0) {\rm exp}\{t^{\frac{1}{4}} C_1(M(t))\},\ t\in[0,T],\ \epsilon\in(0,1].
\end{split}\end{equation*}

\hfill$\Box$

\subsection{Estimates of ${\rm curl}\mathbf{u}$}

Next, we give the estimate for ${\rm curl}\mathbf{u}$ and the weighted estimates of its spatial and temporal derivatives.

\begin{lem}Let the assumptions in Lemma \ref{lemrho} be satisfied. Then we have
\begin{equation}\label{curlu1}\begin{split}
\|{\rm curl}\mathbf{u}\|^2_{H^1}(t)+\|{\rm curl}\mathbf{u}\|^2_{L^2_t({H^2})}\leq C_0(M_0){\rm exp}\{t^{\frac{1}{8}}C_1(M(t)) \},
\end{split}\end{equation}
\begin{equation}\label{curlu2}\begin{split}
\|\epsilon{\rm curl}\mathbf{u}_t\|^2_{L^2}(t)
+\|\epsilon{\rm curl}\mathbf{u}_t\|^2_{L^2_t({H^1})}\leq  C_0(M_0){\rm exp}\{t^{\frac{1}{4}}C_1(M(t)) \},
\end{split}\end{equation}
and
\begin{equation}\label{curlu3}\begin{split}
\|\epsilon {\rm curl}\mathbf{u}\|^2_{H^2}(t)+\|\epsilon^2 {\rm curl}\mathbf{u}\|^2_{L^2_t(H^3)}\leq C_0(M_0){\rm exp}\{(t^{\frac{1}{4}}+\epsilon^2)C_1(M(t)) \},
\end{split}\end{equation}
where $\epsilon\in (0,1]$ and $t\in[0,T]$.
\end{lem}

\textbf{Proof}. {\bl\it Step 1. Proof of \ef{curlu1}.} We first adapt the momentum equation of ${\eqref{ns0}}_2$ into the following form
\begin{equation}\label{curl0}\begin{split}
&\mathbf{u}_t+\mathbf{u}\cdot\nabla\mathbf{u}+\frac{1}{\epsilon}\rho^{-1}(\epsilon \sigma)(\nabla\sigma+\nabla\theta)+\rho^{-1}(\epsilon \sigma)\nabla(\sigma\theta)\\
&\quad=\rho^{-1}(\epsilon \sigma)(-\mu{\rm curl}{\rm curl}\mathbf{u}+(\mu+\nu)\nabla{\rm div}\mathbf{u}+{\rm curl}\mathbf{H}\times\mathbf{H}).
\end{split}\end{equation}

Let $\omega={\rm curl}\mathbf{u}=(\partial_{x_2}u_3-\partial_{x_3}u_2, \partial_{x_3}u_1-\partial_{x_1}u_3,\partial_{x_1}u_2-\partial_{x_2}u_1)^t$.
By virtue of the identity
\begin{equation}\label{ident1}\begin{split}
{\rm curl}({\rm curl}\mathbf{H}\times \mathbf{H})&=(\mathbf{H}\cdot\nabla){\rm curl}\mathbf{H}-({\rm curl}\mathbf{H}\cdot\nabla)\mathbf{H}
+{\rm curl}\mathbf{H}\dv \mathbf{H}-\mathbf{H}\dv ({\rm curl}\mathbf{H})
\\&=(\mathbf{H}\cdot\nabla){\rm curl}\mathbf{H}-({\rm curl}\mathbf{H}\cdot\nabla) \mathbf{H},
\end{split}\end{equation}
applying \lq\lq curl''  to  ${\eqref{curl0}}$ and combining ${\eqref{initial}}$ and ${\eqref{boundary}}$, we obtain
\begin{equation}\label{curl1}\left\{\begin{split}
&\omega_t+\mathbf{u}\cdot\nabla\omega-\mu \rho^{-1}(\epsilon \sigma)\Delta\omega=f_1,\quad &(t,x)\in(0,T]\times\Omega,\\
&\omega\times\mathbf{n}=0,\quad&{\rm on} \quad(0,T]\times\partial\Omega,\\
&\omega(0,x)=\omega_0:={\rm curl}\mathbf{u}_0,\quad &x\in\Omega,
\end{split}\right.\end{equation}
where $f_1:=[\mathbf{u}\cdot \nabla,{\rm curl}]\mathbf{u}
-\epsilon\rho^{-2}(\epsilon \sigma)\nabla \sigma\times(\mu\Delta \mathbf{u} + \nu \nabla \dv\mathbf{u}+{\rm curl}\mathbf{H}\times\mathbf{H})+\rho^{-2}{(\epsilon \sigma)}\nabla\sigma\times\nabla\theta+\epsilon\rho^{-2}{(\epsilon \sigma)}\nabla\sigma\times\nabla(\sigma\theta)+\rho^{-1}(\epsilon \sigma)[(\mathbf{H}\cdot\nabla){\rm curl}\mathbf{H}-({\rm curl}\mathbf{H}\cdot\nabla)\mathbf{H}]$. Here the notation $[a,b]:=ab-ba$. Note that $\rho^{-1}(\epsilon\sigma)\nabla\sigma=\nabla g(\epsilon,\sigma)$ for some function $g$ and the facts that
 \be\label{3curl}
\Delta {\rm curl}  =-{\rm curl}{\rm curl}{\rm curl},
\ee
and
 ${\rm curl}\ \nabla=0$ are used in deriving $\ef{curl1}_1$.
Multiplication of ${\eqref{curl1}_1}$ on $\omega$ and integration of the resulting equality on $L^2((0,t)\times\Omega)$ yield the following estimate
\begin{equation}\label{curl1a}\begin{split}
\|\omega & \|^2_{L^2}(t)+\mu\int_0^t\int_\Omega\rho^{-1}(\epsilon \sigma)|{\rm curl}\omega|^2dxds\\
&=\|\omega_0\|^2_{L^2}+{\int_0^t\int_\Omega(\dv\mathbf{u} \omega +f_1)\cdot\omega dxds}\\
&\leq
\|\omega_0\|^2_{L^2}+C\int_0^t\int_\Omega[|\nabla \mathbf{u}|^2+\epsilon|\rho^{-2}|
 |\nabla \sigma| (|\nabla^2\mathbf{u}|+|\nabla\mathbf{H}\|\mathbf{H}|)
 +|\rho^{-2}| |\nabla\sigma| |\nabla\theta|\\
&\quad+|\rho^{-1}|(|\nabla^2\mathbf{H}\|\mathbf{H}|+|\nabla\mathbf{H}|^2)+\epsilon|\rho^{-2}| |\nabla\sigma| (|\nabla\sigma| |\theta|
 +|\nabla\theta| |\sigma|)] |\nabla\mathbf{u}|dx ds\\
&\leq
\|\omega_0\|^2_{L^2}+C t^{\frac{1}{2}}\{\|\mathbf{u}\|^{\frac{3}{2}}_{L^\infty_t(H^1)}
\|\mathbf{u}\|^{\frac{3}{2}}_{L^\infty_t(H^2)}+\|\rho^{-1}\|^2_{L^\infty_{x,t}}
[\|\epsilon\sigma\|_{L^\infty_t(H^2)}
(\|\mathbf{u}\|^2_{L^\infty_t(H^2)}+\|\mathbf{H}\|^2_{L^\infty_t(H^1)})\\
&\quad+\|\sigma\|_{L^\infty_t(H^1)}\|\theta\|_{L^\infty_t(H^1)}
\|\mathbf{u}\|_{L^2_t(H^3)}+\|\sigma\|_{L^\infty_t(H^1)}\|\mathbf{u}\|_{L^\infty_t(H^2)}
(\|\sigma\|_{L^\infty_t(H^1)}\|\epsilon\theta\|_{L^\infty_t(H^2)}\\
&\quad+\|\epsilon\sigma\|_{L^\infty_t(H^2)}\|\theta\|_{L^\infty_t(H^1)})
+\|\rho^{-1}\|_{L^\infty_{x,t}}\|\mathbf{H}\|^2_{L^\infty_t(H^1)}\|\mathbf{u}\|_{L^\infty_t(H^2)}]\}\\
&\leq C_0(M_0){\rm exp}\{{ t^\frac{1}{2}}C_1(M(t))\},\ t\in[0,T].
\end{split}\end{equation}
Next, we apply \lq\lq curl'' to ${\eqref{ns0}}_2$ to derive that
\begin{equation}\label{curlns02}
\rho(\omega_t
+\mathbf{u}\cdot\nabla\omega)={\bl\epsilon\nabla\sigma\times(\mathbf{u}_t+\mathbf{u}\cdot\nabla\mathbf{u})
+\rho\omega\cdot\nabla\mathbf{u}}+\mu\Delta\omega+{\rm curl}({\rm curl}\mathbf{H}\times\mathbf{H}).
\end{equation}
Then, taking the inner product of $\eqref{curlns02}$ and $-\Delta\omega$ and integrating the resulting equality on $L^2((0,t)\times\Omega)$, we obtain
\begin{equation}\label{curl1b}\begin{split}
\| &\sqrt{\rho}{\rm curl} \omega\|^2_{L^2}(t)+\mu\int_0^t\int_\Omega|\Delta\omega|^2dxds\\
&=\|\sqrt{\rho_0}{\rm curl} \omega_0\|^2_{L^2}+\int_0^t\int_\Omega
[ \epsilon\nabla\sigma\times(\mathbf{u}_t+\mathbf{u}\cdot\nabla\mathbf{u})+\rho(\omega\cdot\nabla\mathbf{u}
+\mathbf{u}\cdot\nabla\omega)]\cdot\Delta\omega dxds\\
&\quad-\int_0^t\int_\Omega\epsilon\nabla\sigma\omega_t\times{\rm curl}\omega dxds+\frac{1}{2}\int_0^t\int_\Omega\rho_t({\rm curl}\omega)^2dxds-\int_0^t\int_\Omega{\rm curl}({\rm curl}\mathbf{H}\times\mathbf{H})\cdot\Delta\omega dxds\\
&\leq C_0(M_0)+Ct^{\frac{1}{8}}\|\mathbf{u}\|_{L^2_t(H^3)}[
\|\sigma\|_{L^\infty_t(H^1)}(\|\epsilon\mathbf{u}_t\|^\frac{1}{4}_{L^\infty_t(H^1)}
\|\epsilon\mathbf{u}_t\|^\frac{3}{4}_{L^2_t(H^2)}
+\|\mathbf{u}\|_{L^\infty_t(H^2)}\|\epsilon\mathbf{u}\|_{L^\infty_t(H^3)}\\
&\quad
+\|\epsilon\mathbf{u}_t\|^\frac{1}{2}_{L^\infty_t(H^1)}
\|\epsilon\mathbf{u}_t\|^\frac{1}{2}_{L^2_t(H^2)})
+\|\rho\|_{L^\infty_{x,t}}\|\mathbf{u}\|^2_{L^\infty_t(H^2)}
+\|\epsilon\sigma_t\|_{L^\infty_t(H^1)}\|\mathbf{u}\|_{L^\infty_t(H^2)}
+\|\mathbf{H}\|^2_{L^\infty_t(H^2)}]\\
&\leq C_0(M_0){\rm exp}\{t^{\frac{1}{8}} C_1(M(t))\},\ t\in[0,T].
\end{split}\end{equation}
Moreover, due to Lemmas \ref{hs2} and \ref{lemma1}, \ef{3curl}, \ef{curl1a}, \ef{curl1b}  and the boundary condition ${\rm curl}\mathbf{u}\times\mathbf{n}|_{\partial\Omega}=0$,
we have
\begin{equation}\label{curl1d}\begin{split}
\|{\rm curl}\mathbf{u}\|^2_{L^2_t(H^1)}&\lesssim\|{\rm curl}{\rm curl}\mathbf{u}\|^2_{L^2_t({L^2})}
+\|{\rm curl}\mathbf{u}\|^2_{L^2_t({L^2})}\\
&\leq C_0(M_0){\rm exp} \{t^{\frac{1}{2}} C_1(M(t))\},\ t\in[0,T],
\end{split}\end{equation}
and
\begin{equation}\label{curlh1}\begin{split}
\|{\rm curl}\mathbf{u}\|^2_{L^\infty_t(H^1)}&\lesssim\|{\rm curl}{\rm curl}\mathbf{u}\|^2_{L^\infty_t({L^2})}
+\|{\rm curl}\mathbf{u}\|^2_{L^\infty_t({L^2})}\\
&\leq  C_0(M_0){\rm exp}\{t^{\frac{1}{8}} C_1(M(t))\},  \ t\in[0,T].
\end{split}\end{equation}
Then by virtue of Lemmas \ref{hs1}, \ref{hs2} and \ref{lembdry}, \ef{curl1b} and \ef{curl1d}, we derive
\begin{equation}\label{curl1c}\begin{split}
\|{\rm curl}\mathbf{u}\|^2_{L^2_t(H^2)}&\lesssim\|{\rm curl}{\rm curl}\mathbf{u}\|^2_{L^2_t(H^1)}+\|{\rm curl}\mathbf{u}\|^2_{L^2_t(H^1)}\\
&\lesssim\|\Delta{\rm curl}\mathbf{u}\|^2_{L^2_t(L^2)}+\|{\rm curl}{\rm curl}\mathbf{u}\cdot\mathbf{n}\|^2_{L^2_t( H^\frac 1 2 (\partial \Omega))}+\|{\rm curl}\mathbf{u}\|^2_{L^2_t(H^1)}\\
&\lesssim\|\Delta{\rm curl}\mathbf{u}\|^2_{L^2_t(L^2)} +\|\mathbf{u}\|^2_{L^2_t(H^2)}\\
&\leq  C_0(M_0){\rm exp} \{t^{\frac{1}{8}} C_1(M(t))\},\ t\in[0,T],
\end{split}\end{equation}
where the estimate $\|\mathbf{u}\|^2_{L^2_t(H^2)}\leq t\|\mathbf{u}\|^2_{L^\infty_t(H^2)}\leq tM(t)$ has been used.
Thus \ef{curlu1} is shown by combining \ef{curlh1} and \ef{curl1c}.

{\bl\it Step 2. Proof of \ef{curlu2}.} Next, we apply ``${\rm curl}\partial_t$'' to ${\eqref{curl0}}$   to derive
\begin{equation}\label{curl2}\left\{\begin{split}
&\omega_{tt}+\mathbf{u}\cdot\nabla\omega_t
 +\mu {\rm curl} (\rho^{-1}(\epsilon \sigma){\rm curl}\omega_t) =f_{2}
,\quad &(t,x)\in(0,T]\times\Omega,\\
&\omega_t\times\mathbf{n}=0,\quad&{\rm on} \quad(0,T]\times\partial\Omega,\\
&\omega_t(0,x):={\rm curl}\mathbf{u}_t(0),\quad &x\in\Omega,
\end{split}\right.\end{equation}
where  $\omega_t(0):={\rm curl}[\rho_0^{-1}(\mu\triangle\mathbf{u}_0+(\mu+\nu)\nabla{\rm div}\mathbf{u}_0)]
-{\rm curl}(\mathbf{u}_0\cdot\nabla\mathbf{u}_0)$ and $f_{2}:= -\mathbf{u}_t\cdot\nabla\omega+ \{[\mathbf{u}\cdot \nabla,{\rm curl}]\mathbf{u}\}_t
-  \epsilon (\mu+\nu) (\rho^{-2}(\epsilon \sigma) \nabla\sigma\times\nabla {\rm div} \mathbf{u})_t
+\mu {\rm curl} (\epsilon \sigma_t\rho^{-2}(\epsilon \sigma){\rm curl}\omega)+(\rho^{-2}(\epsilon \sigma)\nabla\sigma\times\nabla\theta+\epsilon\rho^{-2}(\epsilon \sigma)\nabla\sigma\times\nabla(\sigma\theta))_t-\epsilon( \rho^{-2}(\epsilon \sigma) \nabla\sigma\times{\rm curl}\mathbf{H}\times\mathbf{H})_t+(\rho^{-1}(\epsilon \sigma){\rm curl}({\rm curl}\mathbf{H}\times\mathbf{H}))_t:=-\mathbf{u}_t\cdot\nabla\omega+ \{[\mathbf{u}\cdot \nabla,{\rm curl}]\mathbf{u}\}_t+f_3.$\\
We multiply ${\eqref{curl2}}_1$ by $\epsilon^2\omega_t$ and integrate the resulting equality in $L^2((0,t)\times\Omega)$ to deduce that
\begin{equation}\label{curl2a}\begin{split}
& \|\epsilon   \omega_t  \|^2_{L^2}(t)
+\mu \|\epsilon\sqrt{\rho^{-1}(\epsilon\sigma)}{\rm curl}\omega_t\|^2_{L^2_t({L^2})}\\
&=\|\epsilon\omega_t (0)\|^2_{L^2}
+\int_0^t\int_\Omega({\rm div}\mathbf{u}|\e\omega_t|^2+\e f_3  \cdot \e\omega_t ) dxds\\
&\leq C_0(M_0)+Ct^{\frac{1}{2}}\|\mathbf{u}\|_{L^\infty_t(H^2)}
\|\epsilon\mathbf{u}_t\|_{L^2_t(H^2)}
\|\epsilon\mathbf{u}_t\|_{L^\infty_t(H^1)}
+{\int_0^t\int_\Omega\e f_3  \cdot \e\omega_t  dxds}\\
&\leq C_0(M_0) {\rm exp}\{t^{\frac{1}{4}} C_1(M(t))\},\ t\in[0,T],
\end{split}\end{equation}
where
\begin{equation*} \begin{split}
\int_0^t\int_\Omega\e f_3\cdot\e\omega_tdxds
&=-(\mu+\nu)\int_0^t\int_\Omega\epsilon^2(\rho^{-2}(\epsilon \sigma) \nabla\sigma\times\nabla {\rm div} \mathbf{u})_t\cdot\e\omega_tdxds\\
&\quad\quad+\mu\int_0^t\int_\Omega\epsilon{\rm curl} (\epsilon \sigma_t\rho^{-2}(\epsilon \sigma){\rm curl}\omega)\cdot\e\omega_tdxds\\
&\quad\quad+\int_0^t\int_\Omega\epsilon(\rho^{-2}(\epsilon \sigma)\nabla\sigma\times\nabla\theta+\epsilon\rho^{-2}(\epsilon \sigma)\nabla\sigma\times\nabla(\sigma\theta))_t\cdot\e\omega_tdxds\\
&\quad\quad-\int_0^t\int_\Omega\epsilon^2(\rho^{-2}(\epsilon \sigma) \nabla\sigma\times{\rm curl}\mathbf{H}\times\mathbf{H})_t\cdot\e\omega_tdxds\\
&\quad\quad+\int_0^t\int_\Omega(\rho^{-1}(\epsilon \sigma){\rm curl}({\rm curl}\mathbf{H}\times\mathbf{H}))_t\cdot\e\omega_tdxds\\
&=:A_1+A_2+A_3+A_4+A_5.
\end{split}\end{equation*}
By direct calculations and using Corollary \ref{inter1}, we have
\begin{equation}\label{a1}\begin{split}
A_1&= (\mu+\nu)\int_0^t\int_\Omega\epsilon^2(2\rho^{-3}\epsilon\sigma_t\nabla\sigma\times\nabla{\rm div} \mathbf{u}
-\rho^{-2}\nabla\sigma_t\times\nabla{\rm div} \mathbf{u}-\rho^{-2}\nabla\sigma\times\nabla{\rm div} \mathbf{u}_t)\cdot\e\omega_tdxds\\
&\lesssim t^{\frac{1}{4}}\|\epsilon\mathbf{u}_t\|^\frac{1}{2}_{L^2_t(H^2)}
\|\epsilon\mathbf{u}_t\|^\frac{1}{2}_{L^\infty_t(H^1)}
[\epsilon\|\rho^{-1}\|^3_{L^\infty_{x,t}}\|\epsilon\sigma_t\|_{L^\infty_t(H^1)}
\|\epsilon\sigma\|_{L^\infty_t(H^2)}\|\mathbf{u}\|_{L^2_t(H^3)}\\
&\quad+\|\rho^{-1}\|^2_{L^\infty_{x,t}}(\epsilon\|\epsilon\sigma_t\|_{L^\infty_t(H^1)}
\|\mathbf{u}\|_{L^2_t(H^3)}+\|\epsilon\sigma\|_{L^\infty_t(H^2)}
\|\epsilon\mathbf{u}_t\|_{L^2_t(H^2)})],
\end{split}\end{equation}
\begin{equation}\label{a2}\begin{split}
A_2&=\mu\int_0^t\int_\Omega\epsilon(\rho^{-2}\epsilon\nabla\sigma_t\times{\rm curl}\omega-2\epsilon\rho^{-3}\epsilon\sigma_t\nabla\sigma\times{\rm curl}\omega+\rho^{-2}\epsilon\sigma_t {\rm curl}{\rm curl}\omega)\cdot\e\omega_tdxds\\
&\lesssim t^{\frac{1}{4}}\|\epsilon\mathbf{u}_t\|^\frac{1}{2}_{L^2_t(H^2)}
\|\epsilon\mathbf{u}_t\|^\frac{1}{2}_{L^\infty_t(H^1)}\|\mathbf{u}\|_{L^2_t(H^3)}
(\epsilon\|\rho^{-1}\|^2_{L^\infty_{x,t}}\|\epsilon\sigma_t\|_{L^\infty_t(H^1)}\\
&\quad+\epsilon\|\rho^{-1}\|^3_{L^\infty_{x,t}}\|\epsilon\sigma_t\|_{L^\infty_t(H^1)}
\|\epsilon\sigma\|_{L^\infty_t(H^2)}),
\end{split}\end{equation}
{\begin{equation}\label{a3}\begin{split}
A_3&=\int_0^t\int_\Omega\epsilon[-2\rho^{-3}\epsilon\sigma_t(\nabla\sigma\times\nabla\theta+\epsilon\nabla\sigma\times
(\sigma\nabla\theta+\theta\nabla\sigma))+\rho^{-2}(\nabla\sigma_t\times\nabla\theta+\nabla\sigma\times\nabla\theta_t\\
&\quad+\epsilon\nabla\sigma_t\times
(\sigma\nabla\theta+\theta\nabla\sigma)+\epsilon\nabla\sigma\times(\sigma_t\nabla\theta+\sigma\nabla\theta_t
+\theta_t\nabla\sigma+\theta\nabla\sigma_t
))]\cdot\e\omega_tdxds\\
&{ \lesssim t^{\frac{1}{4}}\|\epsilon\mathbf{u}_t\|^\frac{1}{2}_{L^2_t(H^2)}
\|\epsilon\mathbf{u}_t\|^\frac{1}{2}_{L^\infty_t(H^1)}[\|\rho^{-1}\|^3_{L^\infty_{x,t}}
\|\epsilon\sigma_t\|_{L^\infty_t(H^1)}
(\|\epsilon\sigma\|_{L^\infty_t(H^2)}\|\theta\|_{L^2_t(H^2)}}\\
&\quad+\epsilon\|\epsilon\sigma\|_{L^\infty_t(H^2)}
(\|\sigma\|_{L^\infty_t(H^1)}\|\theta\|_{L^2_t(H^2)}
+\|\theta\|_{L^\infty_t(H^1)}\|\sigma\|_{L^2_t(H^2)})\\
&\quad+\|\rho^{-1}\|^2_{L^\infty_{x,t}}(\|\epsilon\sigma_t\|_{L^\infty_t(H^1)}
\|\theta\|_{L^2_t(H^2)}+
\|\epsilon\theta_t\|_{L^\infty_t(H^1)}
\|\sigma\|_{L^2_t(H^2)}\\
&\quad+\|\epsilon\sigma_t\|_{L^\infty_t(H^1)}
(\|\sigma\|_{L^2_t(H^2)}\|\epsilon\theta\|_{L^\infty_t(H^2)}
+\|\theta\|_{L^2_t(H^2)}\|\epsilon\sigma\|_{L^\infty_t(H^2)})\\
&\quad+\|\epsilon\sigma\|_{L^\infty_t(H^2)}
(\|\epsilon\sigma_t\|_{L^\infty_t(H^1)}\|\theta\|_{L^2_t(H^2)}
+\|\epsilon\theta_t\|_{L^\infty_t(H^1)}\|\sigma\|_{L^2_t(H^2)}))],{\bl (modified)}
\end{split}\end{equation}}
\begin{equation}\label{a4}\begin{split}
A_4&= \int_0^t\int_\Omega\epsilon^2[2\rho^{-3}\epsilon\sigma_t\nabla\sigma\times({\rm curl}\mathbf{H}\times\mathbf{H})
-\rho^{-2}\nabla\sigma_t\times({\rm curl}\mathbf{H}\times\mathbf{H})\\
&\quad-\rho^{-2}\nabla\sigma\times({\rm curl}\mathbf{H}_t\times\mathbf{H}+{\rm curl}\mathbf{H}\times\mathbf{H}_t)]\cdot\e\omega_tdxds\\
&\lesssim t^{\frac{1}{2}}\|\epsilon\mathbf{u}_t\|_{L^2_t(H^2)}
[\epsilon^2\|\rho^{-1}\|^3_{L^\infty_{x,t}}\|\epsilon\sigma_t\|_{L^\infty_t(H^1)}
\|\sigma\|_{L^\infty_t(H^1)}\|\mathbf{H}\|^2_{L^\infty_t(H^2)}\\
&\quad+\|\rho^{-1}\|^2_{L^\infty_{x,t}}(\epsilon\|\epsilon\sigma_t\|_{L^\infty_t(H^1)}
\|\mathbf{H}\|^2_{L^\infty_t(H^2)}+\|\epsilon\sigma\|_{L^\infty_t(H^2)}
{\bl\|\mathbf{H}\|_{L^\infty_t(H^2)}\|\epsilon\mathbf{H}_t\|_{L^\infty_t(H^1)}})],
\end{split}\end{equation}
\begin{equation}\label{a5}\begin{split}
A_5&= \int_0^t\int_\Omega\epsilon[-\rho^{-2}\epsilon\sigma_t{\rm curl}({\rm curl}\mathbf{H}\times\mathbf{H})
+\rho^{-1}{\rm curl}({\rm curl}\mathbf{H}\times\mathbf{H})_t]\cdot\e\omega_tdxds\\
&=-\int_0^t\int_\Omega\epsilon^2\rho^{-2}\sigma_t{\rm curl}({\rm curl}\mathbf{H}\times\mathbf{H})\cdot\e\omega_tdxds
-\int_0^t\int_\Omega({\rm curl}\mathbf{H}\times\mathbf{H})_t\cdot{\rm curl}(\rho^{-1}\e\omega_t)dxds\\
&\lesssim t^{\frac{1}{2}}\|\epsilon\mathbf{u}_t\|_{L^2_t(H^2)}
[\epsilon\|\rho^{-1}\|^2_{L^\infty_{x,t}}\|\epsilon\sigma_t\|_{L^\infty_t(H^1)}
\|\mathbf{H}\|^2_{L^\infty_t(H^2)}\\
&\quad+{\bl\|\mathbf{H}\|_{L^\infty_t(H^2)}\|\epsilon\mathbf{H}_t\|_{L^\infty_t(H^1)}}
(\|\rho^{-1}\|_{L^\infty_{x,t}}+\|\rho^{-1}\|^2_{L^\infty_{x,t}}\|\epsilon\sigma\|_{L^\infty_t(H^2)})],
\end{split}\end{equation}
Then \ef{curlu2} follows by applying Lemma \ref{hs2} and $\ef{curl2}_2$ to \ef{curl2a}.

{\it Step 3. Proof of \ef{curlu3}.}
By virtue of \ef{lemma1}, \ef{firstorder1},  \ef{firstorder2}, \ef{curlu2}, ${\eqref{curlns02}}$, and Lemma \ref{hs1}, we obtain
 \begin{equation*} \begin{split}
\|\epsilon\Delta{\rm curl}\mathbf{u}\|^2_{ L^2}(t)&\lesssim (\epsilon\|\epsilon\sigma\|_{L^\infty_t(H^2)}
\|\sigma\|_{L^\infty_t(H^1)}+\|\rho\|_{L^\infty_{x,t}}^2)
(\|\epsilon\mathbf{u}_t\|^2_{L^\infty_t(H^1)}+\epsilon^2\|\mathbf{u}\|^4_{L^\infty_t(H^2)})\\
&\quad
{\bl+\epsilon^2\|\mathbf{H}\|^4_{L^\infty_t(H^2)}}\\
&\leq C_0(M_0) {\rm exp} \{(t^{\frac{1}{4}}+\epsilon^2) C_1(M(t))\},\ t\in[0,T],
\end{split}\end{equation*}
which gives that
\begin{equation}\label{curl2b}\begin{split}
\|\epsilon{\rm curl}{\rm curl}\mathbf{u}\|^2_{ H^1}(t)&\lesssim
\|\epsilon\Delta{\rm curl}\mathbf{u}\|^2_{L^\infty_t( L^2)} + \|\epsilon{\rm curl}{\rm curl}\mathbf{u}\cdot \mathbf{n}\|^2_{L^\infty_t( H^\frac 1 2 (\partial \Omega))} +
\|\epsilon{\rm curl}{\rm curl}\mathbf{u} \|^2_{L^\infty_t( L^2)}
  \\
  &\lesssim
\|\epsilon\Delta{\rm curl}\mathbf{u}\|^2_{L^\infty_t( L^2)} + \epsilon^2 \| \mathbf{u} \|^2_{L^\infty_t( H^2)}\\
&\leq C_0(M_0) {\rm exp} \{(t^{\frac{1}{4}}+\epsilon^2) C_1(M(t))\},\ t\in[0,T],
\end{split}\end{equation}
by applying \ef{3curl}, and Lemmas \ref{hs1} and \ref{lembdry} and the fact that $\Delta=\nabla\dv -{\rm curl}{\rm curl}$. Thus from \ef{curlu1}, \ef{curl2b} and Lemma \ref{hs2}, we can obtain
that
\begin{equation}\label{curl2c}\begin{split}
\|\epsilon{\rm curl}\mathbf{u}\|^2_{L^\infty_t(H^2)}&\lesssim
\|\epsilon{\rm curl}{\rm curl}\mathbf{u}\|^2_{L^\infty_t(H^1)} + \|\epsilon{\rm curl}\mathbf{u}\|^2_{L^\infty_t(H^1)}
\leq C_0(M_0) {\rm exp} \{(t^{\frac{1}{4}}+\epsilon^2) C_1(M(t))\}.
\end{split}\end{equation}
Then we take the inner product of $\epsilon^2\nabla{\eqref{curlns02}}$ and $\epsilon^2\nabla\Delta{\rm curl}\mathbf{u}$ in $L^2((0,t)\times\Omega)$ and use \ef{3curl} to derive
\begin{equation}\label{curl1e0}\begin{split}
 \mu &\int_0^t\int_\Omega |\epsilon^2\nabla\Delta{\rm curl}\mathbf{u}|^2dxds\\
&=\int_0^t\int_\Omega\epsilon^2[\epsilon\nabla^2\sigma\times(\mathbf{u}_t+\mathbf{u}\cdot\nabla\mathbf{u})
+\epsilon\nabla\sigma\times(\nabla\mathbf{u}_t+\nabla\mathbf{u}\cdot\nabla\mathbf{u}
+\mathbf{u}\cdot\nabla^2\mathbf{u})\\
&\quad
+(\mathbf{H}\cdot\nabla){\rm curl}\mathbf{H}-({\rm curl}\mathbf{H}\cdot\nabla)\mathbf{H}
+\nabla\rho(\omega_t+\omega\cdot\nabla\mathbf{u}+\mathbf{u}\cdot\nabla\omega)\\
&\quad
+\rho(\nabla\omega_t+2\nabla\omega\cdot\nabla\mathbf{u}+\omega\cdot\nabla^2\mathbf{u}+\mathbf{u}\cdot\nabla^2\omega)]\cdot \e^2\nabla\Delta{\rm curl}\mathbf{u} dxds\\
&\lesssim\epsilon  \|\e^2\nabla^2{\rm curl}{\rm curl}\mathbf{u}\|_{L^2_t(L^2)}[\|\e \nabla \sigma\|_{L^2_t(H^2)}
(\|\epsilon\mathbf{u}_t\|_{L^\infty_t(H^1)}+\epsilon\|\mathbf{u}\|^2_{L^\infty_t(H^2)})\\
&\quad+\epsilon\|\mathbf{H}\|^2_{L^\infty_t(H^2)}+\|\rho\|_{L^\infty_{x,t}}^2(\|\epsilon\mathbf{u}_t\|_{L^2_t(H^2)}+\epsilon\|\mathbf{u}\|_{L^\infty_t(H^2)}
\|\mathbf{u}\|_{L^2_t(H^3)})].
\end{split}\end{equation}
Thus, due to \ef{lemma1}, \eqref{curl2a} and \eqref{curl1e0} we obtain
\begin{equation}\label{curl1e}\begin{split}
\|\epsilon^2\nabla\Delta{\rm curl}\mathbf{u}\|_{L^2_t(L^2)}
&\leq C_0(M_0) {\rm exp}\{(t^{\frac{1}{4}}+\epsilon^2) C_1(M(t))\},\ t\in[0,T],\ \epsilon\in(0,1].
\end{split}\end{equation}
Next, by use of  \ef{lemma1}, \ef{3curl}, ${\eqref{curl1e}}$, Lemmas \ref{hs1} and \ref{lembdry} we get
\begin{equation}\label{curl1f}\begin{split}
& \|\epsilon^2  {\rm curl}{\rm curl}\mathbf{u}\|^2_{L^2_t(H^2)} \\
 \lesssim & \|\epsilon^2\Delta{\rm curl} \mathbf{u}\|^2_{L^2_t(H^1)}
 +\epsilon^2\|{\rm curl}{\rm curl}\mathbf{u}\cdot\mathbf{n}\|^2_{L^2_t(H^{\frac{3}{2}}(\partial\Omega))}
 +\epsilon^2\| {\rm curl}{\rm curl}\mathbf{u}\|^2_{L^2_t(H^1)} \\
 \lesssim &\|\epsilon^2\nabla\Delta{\rm curl}\mathbf{u}\|^2_{L^2_t(L^2)}+\epsilon^2\|\mathbf{u}\|^2_{L^2_t(H^3)}\\
 \leq & C_0(M_0) {\rm exp} \{(t^{\frac{1}{4}}+\epsilon^2) C_1(M(t))\},\ t\in[0,T],\ \epsilon\in(0,1].
\end{split}\end{equation}
Thus from \ef{curl1c}, \ef{curl1f} and Lemma \ref{hs2}, we can obtain
that
\begin{equation}\label{curl1g}\begin{split}
\|\epsilon{\rm curl}\mathbf{u}\|^2_{L^2_t(H^3)}&\lesssim
\|\epsilon{\rm curl}{\rm curl}\mathbf{u}\|^2_{L^2_t(H^2)} + \|\epsilon{\rm curl}\mathbf{u}\|^2_{L^2_t(H^2)}
\leq C_0(M_0) {\rm exp} \{(t^{\frac{1}{4}}+\epsilon^2) C_1(M(t))\}.
\end{split}\end{equation}
Finally, from  ${\eqref{curl2c}}$ and ${\eqref{curl1g}}$ we conclude  ${\eqref{curlu3}}$ immediately.
\hfill$\Box$

\subsection{Estimates of low-order derivatives.}

\begin{lem} Let the assumptions in Lemma \ref{lemrho} be satisfied. Then the solution to the system ${\eqref{ns0}}$ satisfies the following estimate
\begin{equation}\label{lemma8}\begin{split}
 &\|(\sigma,\theta,\mathbf{H}_t)\|^2_{L^\infty_t(H^1)}+
 \|(\mathbf{u},{\bl\mathbf{H}},\epsilon\sigma,\epsilon\theta)\|^2_{L^\infty_t(H^2)} + \|(\e \sigma_t, \e \theta_t, \e\mathbf{u}_t)\|^2_{L^\infty_t(H^1)} + \|(\theta,\e \mathbf{u}_t, \e \theta_t)\|^2_{L^2_t(H^2)}\\
&\quad \leq  C_0(M_0){\rm exp}\{ ( t^{\frac{1}{8}}+\epsilon^\frac 1 2 )C_1(M(t))\},\q \e\in (0,1],\ t\in [0,1].
\end{split}\end{equation}
\end{lem}
\textbf{Proof}.  {\it Step 1.} It follows from \ef{basicestimate}, \ef{basicestimate1} and \ef{curlu1} that
\begin{equation}\label{h1bound}\begin{split}
\|(\sigma,\mathbf{u},\theta)\|^2_{L^\infty_t(H^1)}
{\bl+\|\mathbf{H}\|^2_{L^\infty_t(H^2)}}+\|\theta\|^2_{L^2_t(H^2)}
\leq C_0(M_0){\rm exp}\{ t^{\frac{1}{8}}C_1(M(t))\}.
\end{split}\end{equation}
By virtue of $\epsilon{\eqref{ns0}}_4$, we have
\begin{equation}\label{q5}\begin{split}
\|\epsilon\Delta\theta\|^2_{L^\infty_t(L^2)}
&\lesssim\|{\rm div}\mathbf{u}\|^2_{L^\infty_t(L^2)}
+\|\rho\|^2_{L^\infty_{x,t}}(\|\epsilon\theta_t\|^2_{L^\infty_t(L^2)}
+\epsilon^2\|\mathbf{u}\|^2_{L^\infty_t(H^2)}\|\theta\|^2_{L^\infty_t(H^1)})\\
&\quad+\epsilon^2\|\sigma\|^2_{L^\infty_t(H^1)}\|\mathbf{u}\|^2_{L^\infty_t(H^2)}
+\epsilon^4\|\mathbf{u}\|^2_{L^\infty_t(H^2)}\\
&\leq C_0(M_0){\rm exp}\{ ( t^{\frac{1}{4}}+\epsilon^2)C_1(M(t))\},
\end{split}\end{equation}
which yields
\begin{equation}\label{q5p}\begin{split}
\|\epsilon\theta\|^2_{L^\infty_t(H^2)}
&\lesssim C_0(M_0){\rm exp}\{ ( t^{\frac{1}{4}}+\epsilon^2)C_1(M(t))\},
\end{split}\end{equation}
by ${\eqref{basicestimate1}}$ and Lemma \ref{hs1}.
We can also derive the following estimates directly from $\nabla {\eqref{ns0}}_2$ that
{\begin{equation}\label{q4}\begin{split}
\|\epsilon&(\nabla^2\sigma+\nabla^2\theta)\|^2_{L^\infty_t(L^2)}\\
&\lesssim \epsilon^2 (\|\epsilon\nabla^2{\rm div}\mathbf{u}\|^2_{L^\infty_t(L^2)}
+\|\epsilon\nabla{\rm curl}{\rm curl}\mathbf{u}\|^2_{L^\infty_t(L^2)})
+\|\epsilon^2\nabla^2(\sigma \theta)\|^2_{L^\infty_t(H^2)}\\
&\quad+\epsilon^2\|\rho\|^2_{L^\infty_t(H^2)}
(\|\epsilon\mathbf{u}_t\|^2_{L^\infty_t(H^1)}
+\epsilon^2\|\mathbf{u}\|^4_{L^\infty_t(H^2)})+\epsilon^2\|\mathbf{H}\|^2_{L^\infty_t(H^2)}\\
&\leq C_0(M_0){\rm exp}\{( t^{\frac{1}{8}}+\epsilon^\frac 1 2 )C_1(M(t))\},
\end{split}\end{equation}}
since
\begin{equation}\label{q4p}\begin{split}
\|\epsilon^2\nabla^2(\sigma\theta)\|_{L^\infty_t(L^2)}^2& \lesssim  \|\epsilon^2 \sigma\nabla^2\theta\|_{L^\infty_t(L^2)}^2 +\|\epsilon^2 \nabla \sigma\cdot\nabla\theta\|_{L^\infty_t(L^2)}^2 + \|\epsilon^2 \theta\nabla^2 \sigma\|_{L^\infty_t(L^2)}^2\\
&\lesssim  \epsilon^\frac{1}{2}\|\epsilon  \sigma\|_{L^\infty_t(H^2)}^\frac{3}{2}\|\sigma\|_{L^\infty_t(H^1)}^\frac{1}{2}\| \epsilon\nabla^2\theta\|_{L^\infty_t(L^2)}^2\\
 &\quad + \epsilon\|\epsilon \nabla \sigma\|_{L^\infty_t(H^1)} \| \nabla \sigma\|_{L^\infty_t(L^2)} \| \epsilon\nabla\theta\|_{L^\infty_t(H^1)}^2\\
&\quad + \epsilon^\frac{1}{2}\|\epsilon \nabla^2\sigma\|_{L^\infty_t(L^2)}^2\| \epsilon\theta\|_{L^\infty_t(H^2)}^{\frac 3 2} \|\theta\|_{L^\infty_t(H^1)}^{\frac 1 2}\\
&\lesssim   \epsilon^\frac{1}{2}M^2(t).
\end{split}\end{equation}
Therefore we obtain the estimates for $\|(\epsilon\sigma,\epsilon\theta)\|^2_{L^\infty_t(H^2)}$ in \ef{lemma8} from  $\eqref{basicestimate1}$, ${\eqref{q4}}$ and ${\eqref{q5p}}$ as follows
\begin{equation}\label{q5pp}\begin{split}
\|(\epsilon\sigma,\epsilon\theta)\|^2_{L^\infty_t(H^2)}
&\lesssim C_0(M_0){\rm exp}\{ ( t^{\frac{1}{8}}+\epsilon^{\frac{1}{2}})C_1(M(t))\}.
\end{split}\end{equation}
Then from ${\eqref{ns0}}_3$ and {\bl\ef{h1bound}}, we have
\begin{equation}\label{h1}\begin{split}
\|\mathbf{H}_t\|^2_{L^\infty_t(H^1)}=\|{\rm curl}(\mathbf{u}\times\mathbf{H})\|^2_{L^\infty_t(H^1)}
\lesssim\|\mathbf{u}\|^2_{L^\infty_t(H^1)}\|\mathbf{H}\|^2_{L^\infty_t(H^2)}
\lesssim C_0(M_0){\rm exp}\{ t^{\frac{1}{8}}C_1(M(t))\}.
\end{split}\end{equation}
Next, applying Lemma \ref{hs1} to \ef{firstorder1}, \ef{firstorder2} and \ef{curlu2}, we conclude
\begin{equation}\label{u4p}\begin{split}
   \|(\e \sigma_t, \e \theta_t, \e\mathbf{u}_t)\|^2_{L^\infty_t(H^1)} +\|(\e \mathbf{u}_t, \e \theta_t)\|^2_{L^2_t(H^2)}
 \leq  C_0(M_0){\rm exp}\{t^\frac 1 4 C_1(M(t))\},
\end{split}\end{equation}
where  $\mathbf{u}_t\cdot \mathbf{n}|_{\partial \Omega}=0$ is used.

{\it Step 2.}  By virtue of ${\eqref{ns0}}_1$ and Lemma \ref{hs1}, we get
\begin{equation}\label{u1}\begin{split}
\|{\rm div}\mathbf{u}\|^2_{H^1}&\lesssim\|\epsilon \sigma_t\|^2_{H^1}+\epsilon^2(\|\nabla\sigma\|^2_{L^3}\|{\rm div}\mathbf{u}\|^2_{H^1}+\|\sigma\|^2_{L^\infty}\|\nabla{\rm div}\mathbf{u}\|^2_{L^2}+\|\nabla^2\sigma\|^2_{L^2}
\|\mathbf{u}\|^2_{L^\infty}\\
&\quad+\|\nabla\sigma\|^2_{L^3}
\|\nabla\mathbf{u}\|^2_{H^1})
+\|\epsilon\sigma\|^2_{H^2}\|\mathbf{u}\|^2_{H^1}\\
&\lesssim\|\epsilon \sigma_t\|^2_{H^1}+(\epsilon\|\epsilon\sigma\|_{H^2}
\|\sigma\|_{H^1}+\epsilon^\frac{1}{2}\|\epsilon\sigma\|^\frac{3}{2}_{H^2}
\|\sigma\|^\frac{1}{2}_{H^1})\|{\rm div}\mathbf{u}\|^2_{H^1}\\
&\quad+\|\epsilon\sigma\|^2_{H^2}
\|\mathbf{u}\|^\frac{3}{2}_{H^2}\|\mathbf{u}\|^\frac{1}{2}_{H^1}
+\|\epsilon\sigma\|^2_{H^2}\|\mathbf{u}\|^2_{H^1}\\
&\lesssim\|\epsilon \sigma_t\|^2_{H^1}+(\epsilon\|\epsilon\sigma\|_{H^2}
\|\sigma\|_{H^1}+\epsilon^\frac{1}{2}\|\epsilon\sigma\|^\frac{3}{2}_{H^2}
\|\sigma\|^\frac{1}{2}_{H^1})\|{\rm div}\mathbf{u}\|^2_{H^1}\\
&\quad+\|\epsilon\sigma\|^2_{H^2}
\|\mathbf{u}\|^\frac{1}{2}_{H^1}(\|{\rm div}\mathbf{u}\|^\frac{3}{2}_{H^1}+\|{\rm curl}\mathbf{u}\|^\frac{3}{2}_{H^1}+\|\mathbf{u}\|^\frac{3}{2}_{H^1})
+\|\epsilon\sigma\|^2_{H^2}\|\mathbf{u}\|^2_{H^1}.
\end{split}\end{equation}
Thus, utilizing \ef{curlu1}, \ef{h1bound}, \ef{q5pp} and \ef{u4p}, we obtain
\begin{equation}\label{u3}\begin{split}
\|{\rm div}\mathbf{u}\|^2_{L^\infty_t(H^1)}&\lesssim
\|\epsilon \sigma_t\|^2_{L^\infty_t(H^1)}
+\epsilon^\frac{1}{2}M^2(t)\|{\rm div}\mathbf{u}\|^2_{L^\infty_t(H^1)}
+\|\epsilon\sigma\|^2_{L^\infty_t(H^2)}
\|\mathbf{u}\|^2_{L^\infty_t(H^1)}\|{\rm div}\mathbf{u}\|^\frac{3}{2}_{L^\infty_t(H^1)}\\
&\quad
+\|\epsilon\sigma\|^2_{L^\infty_t(H^2)}
\|\mathbf{u}\|^2_{L^\infty_t(H^1)}(\|{\rm curl}\mathbf{u}\|^\frac{3}{2}_{L^\infty_t(H^1)}
+\|\mathbf{u}\|^\frac{3}{2}_{L^\infty_t(H^1)})+
\|\epsilon\sigma\|^2_{L^\infty_t(H^2)}
\|\mathbf{u}\|^2_{L^\infty_t(H^1)},
\end{split}\end{equation}
which implies that
\begin{equation}\label{u3}\begin{split}
\|{\rm div}\mathbf{u}\|^2_{L^\infty_t(H^1)}\leq C_0(M_0){\rm exp}\{(t^{\frac{1}{8}}+\epsilon^\frac 1 2 ) C_1(M(t))\}.
\end{split}\end{equation}
Then by using  ${\eqref{basicestimate}}$, ${\eqref{curlu1}}$, ${\eqref{u3}}$ and Lemma \ref{hs1}  for $s=1,2$, we  find
\begin{equation}\label{u4}\begin{split}
\|\mathbf{u}\|^2_{L^\infty_t(H^2)}
&\lesssim\|{\rm div}\mathbf{u}\|^2_{L^\infty_t(H^1)}+\|{\rm curl}\mathbf{u}\|^2_{L^\infty_t(H^1)}
+\|\mathbf{u}\|^2_{L^\infty_t(H^1)}\\
&\lesssim \|{\rm div}\mathbf{u}\|^2_{L^\infty_t(H^1)}+\|{\rm curl}\mathbf{u}\|^2_{L^\infty_t(H^1)}+\|\mathbf{u}\|^2_{L^\infty_t(L^2)}\\
&\leq C_0(M_0){\rm exp}\{(t^{\frac{1}{8}}+\epsilon^\frac 1 2 ) C_1(M(t))\}.
\end{split}\end{equation}
Therefore, we conclude \ef{lemma8} immediately by collecting the estimates above.
\hfill$\Box$\\

\subsection{Weighted estimates of high-order  derivatives}

\begin{lem}Let the assumptions in Lemma \ref{lemrho} be satisfied. Then for $\e\in (0,1],\ t\in [0,T]$, we have
\begin{equation}\label{lemma5}\begin{split}
&\|\epsilon^2 (\sigma_{tt},\mathbf{u}_{tt},\mathbf{H}_{tt},\theta_{tt})\|^2_{L^2}(t)
+\|\epsilon^2(\mathbf{u}_{tt},\theta_{tt})\|^2_{L^2_t(H^1)}\leq C_0(M_0){\rm exp}\{t^\frac 1 8 C_1(M(t))\}.
\end{split}\end{equation}
\end{lem}
\textbf{Proof}. Differentiating ${\eqref{nstime}}$ with respect to t, we deduce that
\begin{equation}\label{ddtns1}\left\{\begin{split}
&\sigma_{ttt}
+\frac{1}{\epsilon}{\rm div}\mathbf{u}_{tt}=-(\mathbf{u}_{tt}\cdot\nabla \sigma+2\mathbf{u}_t\cdot\nabla \sigma_t+\mathbf{u}\cdot\nabla \sigma_{tt})\\
&\quad\quad-(\sigma_{tt}{\rm div}\mathbf{u}+2\sigma_t{\rm div}\mathbf{u}_t+\sigma{\rm div}\mathbf{u}_{tt}),\quad &(t,x)\in(0,T]\times\Omega,\\
&\rho(\mathbf{u}_{ttt}+\mathbf{u}\cdot\nabla\mathbf{u}_{tt})
+\frac{1}{\epsilon}\nabla (\sigma_{tt}+\theta_{tt})-(\mu\Delta\mathbf{u}_{tt}+\nu\nabla{\rm div}\mathbf{u}_{tt})\\
&\quad=-\rho_{tt}(\mathbf{u}_t+\mathbf{u}\cdot\nabla\mathbf{u})-
2\rho_t(\mathbf{u}_{tt}+\mathbf{u}_t\cdot\nabla\mathbf{u}
+\mathbf{u}\cdot\nabla\mathbf{u}_t)\\
&\quad\quad-\rho(\mathbf{u}_{tt}\cdot\nabla\mathbf{u}
+2\mathbf{u}_t\cdot\nabla\mathbf{u}_t)+\sigma_{tt}\nabla\theta+2\sigma_t\nabla\theta_t\\
&\quad\quad+\sigma\nabla\theta_{tt}
+\theta_{tt}\nabla\sigma+2\theta_t\nabla\sigma_t+\theta\nabla\sigma_{tt}
+\mathbf{H}_{tt}\cdot\nabla\mathbf{H}\\
&\quad\quad+2\mathbf{H}_t\cdot\nabla\mathbf{H}_t
+\mathbf{H}\cdot\nabla\mathbf{H}_{tt}+\nabla(|\mathbf{H}_t|^2),\quad &(t,x)\in(0,T]\times\Omega,\\
&\mathbf{H}_{ttt}={\bl\mathbf{H}_{tt}\cdot\nabla\mathbf{u}-\mathbf{H}_{tt}{\rm div}\mathbf{u}+2\mathbf{H}_t\cdot\nabla\mathbf{u}_t-2\mathbf{H}_t{\rm div}\mathbf{u}_t+\mathbf{H}\cdot\nabla\mathbf{u}_{tt}}\\
&\quad\quad\quad-\mathbf{H}{\rm div}\mathbf{u}_{tt}-(\mathbf{u}_{tt}\cdot\nabla\mathbf{H}+2\mathbf{u}_t\cdot\nabla\mathbf{H}_t
+\mathbf{u}\cdot\nabla\mathbf{H}_{tt}),\quad &(t,x)\in(0,T]\times\Omega,\\
&\rho(\theta_{ttt}+\mathbf{u}\cdot\nabla\theta_{tt})
+\frac{1}{\epsilon}{\rm div}\mathbf{u}_{tt}-\kappa\Delta\theta_{tt}\\
&\quad=-\rho_t(\theta_{tt}+\mathbf{u}\cdot\nabla\theta_t)-\rho\mathbf{u}_t\cdot\nabla\theta_t+
\epsilon(2\mu|D(\mathbf{u})|^2+\xi({\rm div}\mathbf{u})^2)_{tt}\\
&\quad\quad-((\rho\theta+\sigma){\rm div}\mathbf{u})_{tt}-(\rho_t\theta_t+(\rho\mathbf{u})_t\cdot\nabla\theta)_t,\quad &(t,x)\in(0,T]\times\Omega,\\
&\mathbf{u}_{tt}\cdot\mathbf{n}=0,\ \mathbf{n}\times{\rm curl}\mathbf{u}_{tt}=0,\ \frac{\partial\theta_{tt}}{\partial\mathbf{n}}=0, \quad&{\rm on} \quad(0,T]\times\partial\Omega,\\
&(\sigma_{tt},\mathbf{u}_{tt},\mathbf{H}_{tt},\theta_{tt})(0,x)=(\sigma_{tt}(0),\mathbf{u}_{tt}(0),\mathbf{H}_{tt}(0),\theta_{tt}(0)),\quad &x\in\Omega,
\end{split}\right.\end{equation}
where
\begin{equation*}\begin{split}
\sigma_{tt}(0):=-\frac{1}{\epsilon}{\rm div}\mathbf{u}_t(0)-\mathbf{u}_t(0)\cdot\nabla \sigma_0-\mathbf{u}_0\cdot\nabla \sigma_t(0)-\sigma_t(0){\rm div}\mathbf{u}_0-\sigma_0{\rm div}\mathbf{u}_t(0),
\end{split}\end{equation*}
\begin{equation*}\begin{split}
\mathbf{u}_{tt}(0)&:=[-\frac{1}{\epsilon}(\nabla \sigma_t(0)+\nabla \theta_t(0))+\mu\triangle\mathbf{u}_t(0)+\nu\nabla{\rm div}\mathbf{u}_t(0)-\rho_t(0)(\mathbf{u}_t(0)
+\mathbf{u}_0\cdot\nabla\mathbf{u}_0)\\
&\quad-\rho_0\mathbf{u}_t(0)\cdot\nabla\mathbf{u}_0
-\nabla\sigma_t(0)\theta_0-\sigma_t(0)\nabla\theta_0-\nabla\sigma_0\theta_t(0)
-\sigma_0\nabla\theta_t(0)\\
&\quad+{\rm curl}\mathbf{H}_t(0)\times\mathbf{H}_0+{\rm curl}\mathbf{H}_0\times\mathbf{H}_t(0)]/\rho_0-\mathbf{u}_0\cdot\nabla\mathbf{u}_t(0),
\end{split}\end{equation*}
and
\begin{equation*}\begin{split}
\theta_{tt}(0)&:=[-\frac{1}{\epsilon}{\rm div}\mathbf{u}_t(0)+\kappa\Delta\theta_t(0)+2\epsilon(2\mu|D\mathbf{u}_0\|D\mathbf{u}_t(0)|+\xi{\rm div}\mathbf{u}_0{\rm div}\mathbf{u}_t(0))\\
&\quad-(\rho_t(0)\theta_0+\rho_0\theta_t(0)+\sigma_t(0)){\rm div}\mathbf{u}_0-(\rho_0\theta_0+\sigma_0){\rm div}\mathbf{u}_t(0)-\rho_t(0)\theta_t(0)\\
&\quad
-\rho_t(0)\mathbf{u}_0\cdot\nabla\theta_0
-\rho_0\mathbf{u}_t(0)\cdot\nabla\theta_0]/\rho_0
-\mathbf{u}_0\cdot\nabla\theta_t(0),
\end{split}\end{equation*}
\begin{equation*}\begin{split}
\mathbf{H}_{tt}(0):=\mathbf{H}_t(0)\cdot\nabla\mathbf{u}_0+\mathbf{H}_0\cdot\nabla\mathbf{u}_t(0)
-\mathbf{u}_t(0)\cdot\nabla\mathbf{H}_0-\mathbf{u}_0\cdot\nabla\mathbf{H}_t(0)-\mathbf{H}_t(0){\rm div}\mathbf{u}_0-\mathbf{H}_0\cdot{\rm div}\mathbf{u}_t(0).
\end{split}\end{equation*}
Integrating the product of ${\eqref{ddtns1}_1}$ and $\epsilon^4 \sigma_{tt}$ in $L^2(\Omega\times(0,t))$ yields that
\begin{equation}\label{ddtns1a}\begin{split}
&\frac{1}{2}\|\epsilon^2 \sigma_{tt}\|^2_{L^2}(t)
+\epsilon^3\int_0^t\int_\Omega{\rm div}\mathbf{u}_{tt}\cdot\sigma_{tt}dxds\\
&=\frac{1}{2}\|\epsilon^2 \sigma_{tt}(0)| |^2_{L^2}-\int_0^t\int_\Omega\epsilon^4(\mathbf{u}_{tt}\cdot\nabla \sigma+2\mathbf{u}_t\cdot\nabla \sigma_t-\frac 1 2{\rm div}\mathbf{u} \sigma_{tt})\sigma_{tt}dxds\\
&\quad-\int_0^t\int_\Omega\epsilon^4(\sigma_{tt}{\rm div}\mathbf{u}+2\sigma_t{\rm div}\mathbf{u}_t+\sigma{\rm div}\mathbf{u}_{tt})\sigma_{tt}dxds\\
&\lesssim C_0(M_0)+Ct^\frac 1 8\|\epsilon^2\sigma_{tt}\|_{L^\infty_t(L^2)}
(\|\epsilon^2\mathbf{u}_{tt}\|_{L^2_t(H^1)}\|\sigma\|^\frac 1 2_{L^2_t(H^2)}\|\sigma\|^\frac 1 2_{L^\infty_t(H^1)}
+\|\epsilon\mathbf{u}_t\|_{L^2_t(H^2)}\|\epsilon \sigma_t\|_{L^\infty_t(H^1)}\\
&\quad+\|\epsilon^2\sigma_{tt}\|_{L^\infty_t(L^2)}\|\mathbf{u}\|_{L^2_t(H^3)}
+\|\sigma\|^\frac 3 4_{L^2_t(H^2)}\|\sigma\|^\frac 1 4_{L^\infty_t(H^1)}\|\epsilon^2\mathbf{u}_{tt}\|_{L^2_t(H^1)}),
\end{split}\end{equation}
where the boundary condition $\mathbf{u}\cdot\mathbf{n}|_{\partial\Omega}=0$ and  {\bl integration} by parts have been used in the first integral on the right-hand side of the {\bl above inequality.}
Then we deduce the following estimate by multiplying  $\epsilon^4{\eqref{ddtns1}_2}$ on $\mathbf{u}_{tt}$ and integrating in $L^2(\Omega\times(0,t))$:
\begin{equation}\label{ddtns1b}\begin{split}
&\frac{1}{2}\|\epsilon^2 \sqrt{\rho}\mathbf{u}_{tt}\|^2_{L^2}(t)+(\mu+\nu)\|\epsilon^2{\rm div}\mathbf{u}_{tt}\|^2_{L^2_t(L^2)}+\mu\|\epsilon^2{\rm curl}\mathbf{u}_{tt}\|^2_{L^2_t(L^2)}
+\epsilon^3\int_0^t\int_\Omega(\nabla \sigma_{tt}+\nabla \theta_{tt})\cdot\mathbf{u}_{tt}dxds\\
&=\frac{1}{2}\|\epsilon^2 \sqrt{\rho}_0\mathbf{u}_{tt}(0)\|^2_{L^2}-\int_0^t\int_\Omega\epsilon^4[\rho_{tt}(\mathbf{u}_t+\mathbf{u}\cdot\nabla\mathbf{u})+
2\rho_t(\mathbf{u}_{tt}+\mathbf{u}_t\cdot\nabla\mathbf{u}
+\mathbf{u}\cdot\nabla\mathbf{u}_t)]\cdot\mathbf{u}_{tt}dxds\\
&\quad-\int_0^t\int_\Omega\epsilon^4\rho(\mathbf{u}_{tt}\cdot\nabla\mathbf{u}
+2\mathbf{u}_t\cdot\nabla\mathbf{u}_t)\cdot\mathbf{u}_{tt}dxds
-\int_0^t\int_\Omega\epsilon^4{\bl\mathbf{H}_{tt}\cdot\mathbf{H}{\rm div}\mathbf{u}_{tt}}dxds\\
&\quad-\int_0^t\int_\Omega\epsilon^4(\sigma_{tt}\nabla\theta+2\sigma_t\nabla\theta_t+\sigma\nabla\theta_{tt}
+\theta_{tt}\nabla\sigma+2\theta_t\nabla\sigma_t+\theta\nabla\sigma_{tt})\cdot\mathbf{u}_{tt}dxds\\
&\quad+\int_0^t\int_\Omega\epsilon^4[\mathbf{H}_{tt}\cdot\nabla\mathbf{H}
+2\mathbf{H}_t\cdot\nabla\mathbf{H}_t+\nabla(|\mathbf{H}_t|^2)]\cdot\mathbf{u}_{tt}dxds\\
&\lesssim C_0(M_0)+Ct^\frac 1 8\{\|\epsilon^2\mathbf{u}_{tt}\|_{L^2_t(H^1)}
[\|\epsilon^2\sigma_{tt}\|_{L^\infty_t(L^2)}
(\|\epsilon\mathbf{u}_t\|_{L^\infty_t(H^1)}+\epsilon\|\mathbf{u}\|^2_{L^\infty_t(H^2)}
+\|\theta\|^\frac 1 2_{L^2_t(H^2)}\|\theta\|^\frac 1 2_{L^\infty_t(H^1)}\\
&\quad+\|\theta\|^\frac 3 4_{L^2_t(H^2)}\|\theta\|^\frac 1 4_{L^\infty_t(H^1)})+\|\epsilon\sigma_t\|_{L^\infty_t(H^1)}(\|\epsilon^2\mathbf{u}_{tt}\|_{L^\infty_t(L^2)}
+\epsilon\|\epsilon\mathbf{u}_t\|_{L^\infty_t(H^1)}\|\mathbf{u}\|_{L^\infty_t(H^2)})\\
&\quad
+\|\rho\|_{L^\infty_{x,t}}(\|\epsilon^2\mathbf{u}_{tt}\|_{L^\infty_t(L^2)}
\|\mathbf{u}\|_{L^\infty_t(H^2)}+\|\epsilon\mathbf{u}_t\|^2_{L^\infty_t(H^1)})
+\|\epsilon\sigma_t\|_{L^\infty_t(H^1)}\|\epsilon\theta_t\|_{L^\infty_t(H^1)}\\
&\quad+\|\epsilon\mathbf{H}_t\|^2_{L^\infty_t(H^1)}+\|\epsilon^2\mathbf{H}_{tt}\|_{L^\infty_t(L^2)}
\|\mathbf{H}\|_{L^\infty_t(H^2)}]\\
&\quad
+\|\epsilon^2\mathbf{u}_{tt}\|^\frac 1 2_{L^2_t(H^1)}\|\epsilon^2\mathbf{u}_{tt}\|^\frac 1 2_{L^\infty_t(L^2)}
\|\epsilon^2\theta_{tt}\|_{L^2_t(H^1)}\|\sigma\|_{L^\infty_t(H^1)}\},
\end{split}\end{equation}
{ where we have used that
\begin{equation*}\begin{split}
\int_0^t\int_\Omega\theta\nabla\sigma_{tt}\cdot\mathbf{u}_{tt}dxds
&=-\int_0^t\int_\Omega\sigma_{tt}\cdot(\theta{\rm div}\mathbf{u}_{tt}+{\bl\nabla\theta\cdot\mathbf{u}_{tt}})dxds\\
&\lesssim\|\epsilon^2\sigma_{tt}\|_{L^\infty_t(L^2)}\|\epsilon^2\mathbf{u}_{tt}\|_{L^2_t(H^1)}
(\|\nabla\theta\|_{L^2_t(L^3)}+\|\theta\|_{L^2_t(L^\infty)}).
\end{split}\end{equation*}}
Next by taking the inner product of ${\eqref{ddtns1}_3}$ and $\epsilon^4\mathbf{H}_{tt}$ in $L^2(\Omega\times(0,t))$ we deduce that
\begin{equation}\label{ddtns1b0}\begin{split}
\frac{1}{2}\|\epsilon^2 \mathbf{H}_{tt}\|^2_{L^2}(t)
&=\frac{1}{2}\|\epsilon^2\mathbf{H}_{tt}(0)\|^2_{L^2}+\epsilon^4\int_0^t\int_\Omega{\rm div}\mathbf{u}|\mathbf{H}_{tt}|^2dxds+\epsilon^4\int_0^t\int_\Omega({\bl\mathbf{H}_{tt}\cdot\nabla\mathbf{u}}
-\mathbf{H}_{tt}{\rm div}\mathbf{u}\\
&\quad+2\mathbf{H}_t\cdot\nabla\mathbf{u}_t-2\mathbf{H}_t{\rm div}\mathbf{u}_t+\mathbf{H}\cdot\nabla\mathbf{u}_{tt}-\mathbf{H}{\rm div}\mathbf{u}_{tt}
-\mathbf{u}_{tt}\cdot\nabla\mathbf{H}-2\mathbf{u}_t\cdot\nabla\mathbf{H}_t)\cdot\mathbf{H}_{tt}dxds\\
&\leq C_0(M_0)+Ct^\frac{1}{2}\|\epsilon^2 \mathbf{H}_{tt}\|_{L^\infty_t(L^2)}(\|\epsilon^2 \mathbf{H}_{tt}\|_{L^\infty_t(L^2)}\|\mathbf{u}\|_{L^2_t(H^3)}
+\|\epsilon\mathbf{H}_t\|_{L^\infty_t(H^1)}\|\epsilon\mathbf{u}_t\|_{L^2_t(H^2)}\\
&\quad+\|\mathbf{H}\|_{L^\infty_t(H^2)}\|\epsilon^2\mathbf{u}_{tt}\|_{L^2_t(H^1)}).
\end{split}\end{equation}
Then we integrate the product of ${\eqref{ddtns1}_4}$ and $\epsilon^4\theta_{tt}$ in $L^2(\Omega\times(0,t))$ to obtain that
\begin{equation}\label{ddtns1c}\begin{split}
\frac{1}{2}&\|\epsilon^2 \sqrt{\rho}\theta_{tt}\|^2_{L^2}(t)+\kappa\|\epsilon^2\nabla\theta_{tt}\|^2_{L^2_t(L^2)}
+\epsilon^3\int_0^t\int_\Omega{\rm div}\mathbf{u}_{tt}\cdot\theta_{tt}dxds\\
&=\frac{1}{2}\|\epsilon^2 \sqrt{\rho}_0\theta_{tt}(0)\|^2_{L^2}-\int_0^t\int_\Omega\epsilon^4[\rho_t(\theta_{tt}+\mathbf{u}\cdot\nabla\theta_t)+\rho\mathbf{u}_t\cdot\nabla\theta_t+\epsilon(2\mu|D(\mathbf{u})|^2+\xi({\rm div}\mathbf{u})^2)_{tt}\\
&\quad-((\rho\theta+\sigma){\rm div}\mathbf{u})_{tt}-(\rho_t\theta_t+(\rho\mathbf{u})_t\cdot\nabla\theta)_t]\cdot\theta_{tt}dxds\\
&\lesssim  C_0(M_0)+ Ct^\frac 1 4 \{\|\epsilon^2\theta_{tt}\|_{L^2_t(H^1)}[\|\epsilon\sigma_t\|_{L^\infty_t(H^1)}
(\|\epsilon^2\theta_{tt}\|_{L^\infty_t(L^2)}
+\epsilon\|\mathbf{u}\|_{L^\infty_t(H^2)}\|\epsilon\theta_t\|_{L^\infty_t(H^1)})\\
&\quad
+\|\rho\|_{L^\infty_{x,t}}\|\epsilon\mathbf{u}_t\|_{L^\infty_t(H^1)}
\|\epsilon\theta_t\|_{L^\infty_t(H^1)}+\|\mathbf{u}\|_{L^\infty_t(H^2)}
(\epsilon\|\epsilon^2\sigma_{tt}\|_{L^\infty_t(L^2)}\|\theta\|_{L^\infty_t(H^1)}
\\
&\quad+\epsilon\|\epsilon\sigma_t\|_{L^\infty_t(H^1)}\|\epsilon\theta_t\|_{L^\infty_t(H^1)}
+\|\rho\|_{L^\infty_{x,t}}\|\epsilon^2\theta_{tt}\|_{L^\infty_t(L^2)}
+\|\epsilon^2\sigma_{tt}\|_{L^\infty_t(L^2)})\\
&\quad+\|\epsilon^2\sigma_{tt}\|_{L^\infty_t(L^2)}
\|\epsilon\theta_t\|_{L^\infty_t(H^1)}
+\|\epsilon^2\theta_{tt}\|_{L^\infty_t(L^2)}\|\epsilon\sigma_t\|_{L^\infty_t(H^1)}\\
&\quad+\|\epsilon\mathbf{u}_t\|_{L^\infty_t(H^1)}
(\epsilon\|\epsilon\sigma_t\|_{L^\infty_t(H^1)}\|\theta\|_{L^\infty_t(H^1)}
+\|\rho\|_{L^\infty_{x,t}}\|\epsilon\theta_t\|_{L^\infty_t(H^1)}
+\|\epsilon\sigma_t\|_{L^\infty_t(H^1)})\\
&\quad
+\|\epsilon\theta_t\|_{L^\infty_t(H^1)}
(\|\rho\|_{L^\infty_{x,t}}\|\epsilon\mathbf{u}_t\|_{L^\infty_t(H^1)}
+\epsilon\|\epsilon\sigma_t\|_{L^\infty_t(H^1)}\|\mathbf{u}\|_{L^\infty_t(H^2)})]\\
&\quad+\|\epsilon^2\theta_{tt}\|^\frac{1}{2}_{L^\infty_t(L^2)}
\|\epsilon^2\theta_{tt}\|^\frac{1}{2}_{L^2_t(H^1)}
[\|\epsilon^2\mathbf{u}_{tt}\|_{L^2_t(H^1)}
(\|\rho\|_{L^\infty_{x,t}}\|\theta\|_{L^\infty_t(H^1)}
+\|\sigma\|_{L^\infty_t(H^1)})\\
&\quad+\epsilon(\|\mathbf{u}\|_{L^\infty_t(H^2)}
\|\epsilon^2\mathbf{u}_{tt}\|_{L^2_t(H^1)}+\|\epsilon\mathbf{u}_t\|_{L^\infty_t(H^1)}
\|\epsilon\mathbf{u}_t\|_{L^2_t(H^2)})\\
&\quad+\|\theta\|_{L^2_t(H^2)}(\|\epsilon^2\sigma_{tt}\|_{L^\infty_t(L^2)}
\|\mathbf{u}\|_{L^\infty_t(H^2)}+\epsilon\|\epsilon\sigma_t\|_{L^\infty_t(H^1)}
\|\epsilon\mathbf{u}_t\|_{L^\infty_t(H^1)}
+\|\rho\|_{L^\infty_{x,t}}\|\epsilon^2\mathbf{u}_{tt}\|_{L^\infty_t(L^2)})]\}.
\end{split}\end{equation}
Hence by summarizing {\bl${\eqref{ddtns1a}}$-${\eqref{ddtns1c}}$} and using Lemma \ref{hs1}, we conclude \ef{lemma5} immediately.
\hfill$\Box$

Next, we carry out the weighted  estimates of highest order, which are crucial to close the uniform estimates.
\begin{lem}Let the assumptions in Lemma \ref{lemrho} be satisfied. Then the solution to the system ${\eqref{ns0}}$ satisfies the following estimate
\begin{equation}\label{lemma6}\begin{split}
&\|(  \epsilon  \nabla^2 \sigma_t,\epsilon \nabla^2 \theta_t, \epsilon \nabla^3 \sigma,\epsilon \nabla^3 \theta, \epsilon^2 \nabla^3{\rm div}\mathbf{u}, \epsilon^2\nabla^4\theta )\|^2_{L^2_t(L^2)}+\|\mathbf{u}\|^2_{L^2_t(H^3)}\\
&\quad+\|(\epsilon \nabla^2{\rm div}\mathbf{u},\epsilon^2 \nabla^3 \sigma,\epsilon^2 \nabla^3 \theta)\|^2_{L^2}(t) \\
&\leq  C_0(M_0){\rm exp}\{( t^{\frac{1}{8}}+\epsilon^\frac{1}{2})C_1(M(t))\},
\end{split}\end{equation}
for each $\e\in (0,1],\ t\in [0,T].$
\end{lem}
\textbf{Proof}. {\it Step 1.}
We first derive  the estimates  of $\|\epsilon^2\nabla^3{\rm div}\mathbf{u}\|^2_{L^2_t(L^2)}$, $\|\epsilon^2\nabla^3 \sigma\|^2_{L^\infty_t(L^2)}$ and $\|\epsilon^2\nabla^3 \theta\|^2_{L^\infty_t(L^2)}$. To this end, we take the inner product of $\epsilon^2\partial_i\nabla{\eqref{ns0}}_2$ and $\epsilon^2\partial_i\nabla^2{\rm div}\mathbf{u}$ in $L^2(\Omega)$ to obtain that
\begin{equation}\label{dxns1u}\begin{split}
(\mu &+\nu)\|\epsilon^2\partial_i\nabla^2{\rm div}\mathbf{u}\|^2_{L^2}-\epsilon^3\int_\Omega(\partial_i\nabla^2 \sigma+\partial_i\nabla^2 \theta)\cdot\partial_i\nabla^2{\rm div}\mathbf{u}dx\\
&=\int_\Omega\epsilon^2\partial_i\nabla^2{\rm div}\mathbf{u}\cdot\epsilon^2\partial_i\nabla
[\rho(\mathbf{u}_t+\mathbf{u}\cdot\nabla\mathbf{u})+\nabla(\sigma\theta) +\mu{\rm curl}{\rm curl}\mathbf{u}-{\rm curl}\mathbf{H}\times\mathbf{H}]dx\\
&\lesssim\|\epsilon^2\partial_i\nabla^2{\rm div}\mathbf{u}\|_{L^2}[\|\rho\|_{H^2}(\epsilon\|\epsilon\mathbf{u}_t\|_{H^2}
+\epsilon^2\|\mathbf{u}\|_{H^2}\|\mathbf{u}\|_{H^3})+\|\epsilon^2\nabla\sigma\|_{L^\infty}
\|\theta\|_{H^2}+\|\epsilon^2\nabla\theta\|_{L^\infty}
\|\sigma\|_{H^2}\\
&\quad+\|\sigma\|_{L^\infty}\|\epsilon^2\partial_i\nabla^2 \theta\|_{L^2}+\|\theta\|_{L^\infty}\|\epsilon^2\partial_i\nabla^2 \sigma\|_{L^2}+ \|\epsilon^2{\rm curl}{\rm curl}\mathbf{u}\|_{H^2}+\|\epsilon^2\mathbf{H}\|_{H^3}\|\mathbf{H}\|_{H^2}]\\
&\lesssim \delta\|\epsilon^2\partial_i\nabla^2{\rm div}\mathbf{u}\|^2_{L^2}+C_\delta[\epsilon^2\|\rho\|^2_{H^2}
(\|\epsilon\mathbf{u}_t\|^2_{H^2}
+ \|\mathbf{u}\|^2_{H^2}\|\mathbf{u}\|^2_{H^3})
+\|\epsilon^2\nabla\sigma\|^2_{L^\infty}
\|\theta\|^2_{H^2}\\
&\quad+\|\epsilon^2\nabla\theta\|^2_{L^\infty}
\|\sigma\|^2_{H^2}+\|\sigma\|^2_{L^\infty}\|\epsilon^2\partial_i\nabla^2 \theta\|^2_{L^2}+\|\theta\|^2_{L^\infty}\|\epsilon^2\partial_i\nabla^2 \sigma\|^2_{L^2}\\
&\quad+\|\epsilon^2{\rm curl}{\rm curl}\mathbf{u}\|^2_{H^2}+\|\epsilon^2\mathbf{H}\|^2_{H^3}\|\mathbf{H}\|^2_{H^2}].
\end{split}\end{equation}
Next we multiply $\epsilon^2\partial_i\nabla^2{\eqref{ns0}}_1$ by $\epsilon^2\partial_i\nabla^2 \sigma$ in $L^2(\Omega)$ and use the boundary condition $\mathbf{u}\cdot\mathbf{n}|_{\partial\Omega}=0$ to obtain that
\begin{equation}\label{dxns1s}\begin{split}
\frac{1}{2}&\frac{d}{dt}\|\epsilon^2\partial_i\nabla^2 \sigma\|^2_{L^2}+\epsilon^3\int_\Omega\partial_i\nabla^2{\rm div}\mathbf{u}\cdot\partial_i\nabla^2 \sigma dx\\
&=-\int_\Omega\epsilon^2\partial_i\nabla^2 \sigma\cdot(\epsilon^2\partial_i\nabla^2(\mathbf{u}\cdot\nabla \sigma)+\epsilon^2\partial_i\nabla^2(\sigma{\rm div}\mathbf{u}))dx\\
&=\frac{1}{2}\int_\Omega{\rm div}\mathbf{u}|\epsilon^2\partial_i\nabla^2 \sigma|^2dx-\int_\Omega\epsilon^2\partial_i\nabla^2 \sigma\cdot\epsilon^2[\partial_i\nabla^2,\mathbf{u}]\cdot\nabla \sigma dx-\int_\Omega\epsilon^2\partial_i\nabla^2
\sigma\cdot\epsilon^2\partial_i\nabla^2(\sigma{\rm div}\mathbf{u})dx\\
&:=I_1+I_2+I_3,
\end{split}\end{equation}
where the operator $[a,b]=ab-ba$ and
\begin{equation}\label{dxns1d1}\begin{split}
I_1\lesssim\|\mathbf{u}\|_{H^3}\|\epsilon^2\partial_i\nabla^2\sigma\|^2_{L^2},
\end{split}\end{equation}
\begin{equation}\label{dxns1d2}\begin{split}
I_2&\leq\int_\Omega|\epsilon^2\partial_i\nabla^2 \sigma|\epsilon^2(|\partial_i\nabla^2\mathbf{u}\|\nabla \sigma|+|\partial_i\nabla\mathbf{u}\|\nabla^2\sigma|+|\nabla\mathbf{u}\|\partial_i\nabla^2 \sigma|)dx\lesssim\|\epsilon^2\nabla \sigma\|^2_{H^2}\|\mathbf{u}\|_{H^3},
\end{split}\end{equation}
\begin{equation}\label{dxns1d3}\begin{split}
I_3&\leq\int_\Omega|\epsilon^2\partial_i\nabla^2\sigma|
\epsilon^2({ |\nabla^3 \sigma|}|{\rm div}\mathbf{u}|+|\nabla^2 \sigma|{ |\nabla{\rm div}\mathbf{u}|}+|\nabla \sigma|{ |\nabla^2{\rm div}\mathbf{u}|}+|\sigma\|\partial_i\nabla^2{\rm div}\mathbf{u}|) dx\\
&\lesssim\|\epsilon^2\partial_i\nabla^2\sigma\|_{L^2}(\|\epsilon^2\nabla \sigma\|_{H^2}\|\mathbf{u}\|_{H^3}+\|\sigma\|_{L^\infty}\|\epsilon^2\partial_i\nabla^2{\rm div}\mathbf{u}\|_{L^2}).
\end{split}\end{equation}
Similarly, multiplying $\epsilon^2\partial_i\nabla^2{\eqref{ns0}}_3$ on $\epsilon^2\partial_i\nabla^2 \theta$ in $L^2(\Omega)$ yields that
\begin{equation}\label{dxns1t}\begin{split}
\frac{1}{2}&\frac{d}{dt}\|\epsilon^2\sqrt{\rho}\partial_i\nabla^2 \theta\|^2_{L^2}+\kappa\|\epsilon^2\partial_i\nabla^3 \theta\|^2_{L^2}
+\epsilon^3\int_\Omega\partial_i\nabla^2{\rm div}\mathbf{u}\cdot\partial_i\nabla^2 \theta dx\\
&{=\int_\Omega\epsilon^2\{\partial_i\nabla^2\rho(\theta_t+\mathbf{u}\cdot\nabla\theta)
+\partial_i\nabla\rho\cdot\nabla(\theta_t+\mathbf{u}\cdot\nabla\theta)
+\partial_i\rho\nabla^2(\theta_t+\mathbf{u}\cdot\nabla\theta)}\\
&\quad
+\rho(\partial_i\nabla^2\mathbf{u}\cdot\nabla\theta+\partial_i\nabla\mathbf{u}\cdot\nabla^2\theta
+\partial_i\mathbf{u}\cdot\nabla^3\theta) +\partial_i\nabla^2[(\rho\theta+\sigma){\rm div}\mathbf{u}]\\ &\quad+\partial_i\nabla^2[\epsilon(2\mu|D(\mathbf{u})|^2+\xi({\rm div}\mathbf{u})^2)]\}\cdot\partial_i\nabla^2 \theta dx\\
&\lesssim\|\epsilon^2\partial_i\nabla^2 \theta\|_{L^2}[\|\epsilon^2\partial_i\nabla^2 \sigma\|_{L^2}
(\|\epsilon\theta_t\|_{H^2}+\|\mathbf{u}\|_{H^2}\|\epsilon\theta\|_{H^3})+\|\rho\|_{L^\infty}
\|\epsilon^2\partial_i\nabla^2 \theta\|_{L^2}\|\mathbf{u}\|_{H^3}\\
&\quad+(\|\epsilon^2\partial_i\nabla^2 \sigma\|_{L^2}\|\epsilon\theta\|_{H^2}+\|\epsilon^2\partial_i\nabla^2 \sigma\|_{L^2})\|\mathbf{u}\|_{H^3}+\|\epsilon\mathbf{u}\|_{H^3}\|\epsilon^2\partial_i\nabla^2{\rm div}\mathbf{u}\|_{L^2}].
\end{split}\end{equation}
By use of Corollary \ref{inter}, we combine {\bl${\eqref{dxns1u}}$-${\eqref{dxns1t}}$} and
integrate the resulting inequality on $[0,t]$ to derive that
\begin{equation}\label{dxns116}\begin{split}
 &\|\epsilon^2\partial_i\nabla^2{\rm div}\mathbf{u}\|^2_{L^2_t(L^2)}
+\|\epsilon^2\partial_i\nabla^2 \sigma\|^2_{L^2}(t)+\|\epsilon^2\sqrt{\rho}\partial_i\nabla^2 \theta\|^2_{L^2}(t)
+\kappa\|\epsilon^2\partial_i\nabla^3 \theta\|^2_{L^2_t(L^2)}\\
&\lesssim\|\epsilon^2\partial_i\nabla^2 \sigma_0\|^2_{L^2}
+\|\epsilon^2\partial_i\nabla^2 \theta_0\|^2_{L^2}+ \epsilon^2 \|\rho\|^2_{L^\infty_t(H^2)}
(\|\epsilon\mathbf{u}_t\|^2_{L^2_t(H^2)}
+ \|\mathbf{u}\|^2_{L^\infty_t(H^2)}\|\mathbf{u}\|^2_{L^2_t(H^3)})\\
&\quad+\epsilon^\frac 1 2(\|\epsilon^2\partial_i\nabla^2 \sigma\|^\frac 3 2_{L^\infty_t(L^2)}\|\epsilon\sigma\|^\frac 1 2_{L^\infty_t(H^2)}\|\theta\|^2_{L^2_t(H^2)}
+\|\epsilon^2\partial_i\nabla^2 \theta\|^\frac 3 2_{L^\infty_t(L^2)}\|\epsilon\theta\|^\frac 1 2_{L^\infty_t(H^2)}\|\sigma\|^2_{L^2_t(H^2)})\\
&\quad+t^\frac 1 4(\|\sigma\|^\frac 3 2_{L^2_t(H^2)}\|\sigma\|^\frac 1 2_{L^\infty_t(H^1)}\|\epsilon^2\partial_i\nabla^2 \theta\|^2_{L^\infty_t(L^2)}+\|\theta\|^\frac 3 2_{L^2_t(H^2)}\|\theta\|^\frac 1 2_{L^\infty_t(H^1)}\|\epsilon^2\partial_i\nabla^2 \sigma\|^2_{L^\infty_t(L^2)})\\
&\quad
+\|\epsilon^2{\rm curl}{\rm curl}\mathbf{u}\|^2_{L^2_t(H^2)}+t^\frac 1 8 \{\|\epsilon^2\nabla \sigma\|^2_{L^\infty_t(H^2)}\|\mathbf{u}\|_{L^2_t(H^3)}+\|\rho\|_{L^\infty_{x,t}}
\|\epsilon^2\partial_i\nabla^2 \theta\|^2_{L^\infty_t(L^2)}\|\mathbf{u}\|_{L^2_t(H^3)}\\
&\quad+\|\epsilon^2\partial_i\nabla^2 \theta\|_{L^\infty_t(L^2)}[\|\epsilon^2\partial_i\nabla^2 \sigma\|_{L^\infty_t(L^2)}
(\|\epsilon\theta_t\|_{L^2_t(H^2)}+\|\mathbf{u}\|_{L^\infty_t(H^2)}\|\epsilon\theta\|_{L^2_t(H^3)})\\
&\quad+(\|\epsilon^2\partial_i\nabla^2 \sigma\|_{L^\infty_t(L^2)}\|\epsilon\theta\|_{L^\infty_t(H^2)}+\|\epsilon^2\partial_i\nabla^2 \sigma\|_{L^\infty_t(L^2)})\|\mathbf{u}\|_{L^2_t(H^3)}+\|\epsilon\mathbf{u}\|_{L^\infty_t(H^3)}\|\epsilon^2\partial_i\nabla^2{\rm div}\mathbf{u}\|_{L^2_t(L^2)}]\\
&\quad+\|\epsilon^2\partial_i\nabla^2\sigma\|_{L^\infty_t(L^2)}
\|\sigma\|^\frac 1 4_{L^\infty_t(H^1)}\|\sigma\|^\frac 3 4_{L^2_t(H^2)}\|\epsilon^2\partial_i\nabla^2{\rm div}\mathbf{u}\|_{L^2_t(L^2)}
+t\|\epsilon^2\mathbf{H}\|^2_{L^\infty_t(H^3)}\|\mathbf{H}\|^2_{L^\infty_t(H^2)}\} \\
&\leq  C_0(M_0){\rm exp}\{( t^{\frac{1}{8}}+\epsilon^\frac{1}{2})C_1(M(t))\}.
\end{split}\end{equation}

{\it Step 2.} Next, we derive the bounds of $\|\epsilon\nabla^3\theta\|^2_{L^2_t(L^2)}$ and $\|\epsilon\nabla^3\sigma\|^2_{L^2_t(L^2)}$.
Applying $\epsilon^2\partial_i\nabla$ to ${\eqref{ns0}}_2$ and using $\Delta\mathbf{u}=\nabla{\rm div}\mathbf{u}-{\rm curl}{\rm curl}\mathbf{u}$, we derive
\begin{equation}\label{dxns1b}\begin{split}
&\|\epsilon\partial_i\nabla^2 \sigma+\epsilon\partial_i\nabla^2 \theta\|_{L^2_t(L^2)}^2\\
&\quad\lesssim \|\epsilon^2\partial_i\nabla^2 {\rm div}\mathbf{u}\|_{L^2_t(L^2)}^2+\|\epsilon^2\partial_i\nabla{\rm curl}{\rm curl}\mathbf{u}\|_{L^2_t(L^2)}^2+t\|\epsilon^2\mathbf{H}\|^2_{L^\infty_t(H^3)}\|\mathbf{H}\|^2_{L^\infty_t(H^2)}\\
&\quad\quad+\epsilon^2 \|\rho\|_{L^\infty_t(H^2)}^2
(\|\epsilon\mathbf{u}_t\|_{L^2_t(H^2)}^2+\epsilon^2\|\mathbf{u}\|_{L^\infty_t(H^2)}^2
\|\mathbf{u}\|_{L^2_t(H^3)}^2)+\|\sigma\|^2_{L^2_t(L^\infty)}\|\epsilon^2\partial_i\nabla^2 \theta\|^2_{L^\infty_t(L^2)}\\
&\quad\quad+\|\theta\|^2_{L^2_t(L^\infty)}\|\epsilon^2\partial_i\nabla^2 \sigma\|^2_{L^\infty_t(L^2)}+\|\nabla\sigma\|^2_{L^2_t(L^3)}
\|\epsilon^2\nabla^2 \theta\|^2_{L^\infty_t(H^1)}
+\|\nabla\theta\|^2_{L^2_t(L^3)}\|\epsilon^2\nabla^2 \sigma\|^2_{L^\infty_t(H^1)}\\
&\quad\lesssim \|\epsilon^2\partial_i\nabla^2 {\rm div}\mathbf{u}\|_{L^2_t(L^2)}^2+\|\epsilon^2\partial_i\nabla{\rm curl}{\rm curl}\mathbf{u}\|_{L^2_t(L^2)}^2+t\|\epsilon^2\mathbf{H}\|^2_{L^\infty_t(H^3)}\|\mathbf{H}\|^2_{L^\infty_t(H^2)}\\
&\quad\quad+\epsilon^2 \|\rho\|_{L^\infty_t(H^2)}^2
(\|\epsilon\mathbf{u}_t\|_{L^2_t(H^2)}^2+\epsilon^2\|\mathbf{u}\|_{L^\infty_t(H^2)}^2
\|\mathbf{u}\|_{L^2_t(H^3)}^2)\\
&\quad\quad+t^\frac 1 4(\|\sigma\|^\frac 3 2_{L^2_t(H^2)}\|\sigma\|^\frac 1 2_{L^\infty_t(H^1)}\|\epsilon^2\partial_i\nabla^2 \theta\|^2_{L^\infty_t(L^2)}+\|\theta\|^\frac 3 2_{L^2_t(H^2)}\|\theta\|^\frac 1 2_{L^\infty_t(H^1)}\|\epsilon^2\partial_i\nabla^2 \sigma\|^2_{L^\infty_t(L^2)})\\
&\quad\quad+{\bl t^\frac 1 2}(\|\sigma\|_{L^2_t(H^2)}\|\sigma\|_{L^\infty_t(H^1)}\|\epsilon^2\nabla^2 \theta\|^2_{L^\infty_t(H^1)}+\|\theta\|_{L^2_t(H^2)}\|\theta\|_{L^\infty_t(H^1)}
\|\epsilon^2\nabla^2 \sigma\|^2_{L^\infty_t(H^1)})\\
&\leq C_0(M_0){\rm exp}\{( t^{\frac{1}{8}}+\epsilon^\frac{1}{2})C_1(M(t))\},
\end{split}\end{equation}
where the previous estimates  \eqref{curl1f} and \ef{dxns116}, and Corollary \ref{inter} are utilized.
We can also deduce the following estimate resulting from $\epsilon\nabla{\eqref{ns0}}_3$ that
\begin{equation}\label{dxns1c}\begin{split}
\|\epsilon\partial_i\nabla^2 \theta\|_{L^2_t(L^2)}^2 &\lesssim \epsilon^4\|\mathbf{u}\|_{L^\infty_t(H^2)}^2
\|\mathbf{u}\|_{L^2_t(H^3)}^2+ t\|\mathbf{u}\|_{L^\infty_t(H^2)}^2
+\epsilon^2\|\mathbf{u}\|_{L^\infty_t(H^2)}^2\|\sigma\|_{L^2_t(H^2)}^2\\
&\quad+\|\rho\|_{L^\infty_t(H^2)}^2
({t\|\epsilon\theta_t\|_{L^\infty_t(H^1)}^2}+\epsilon^2\|\mathbf{u}\|_{L^\infty_t(H^2)}^2
\|\theta\|_{L^2_t(H^2)}^2)\\
&\leq C_0(M_0){\rm exp}\{( t^{\frac{1}{4}}+\epsilon^{\frac{1}{2}})C_1(M(t))\}.
\end{split}\end{equation}
Thus, by virtue of ${\eqref{dxns1b}}$ and ${\eqref{dxns1c}}$, we can obtain
\begin{equation}\label{dxns1d}\begin{split}
\|\epsilon\partial_i\nabla^2 \sigma\|_{L^2_t(L^2)}^2\leq C_0(M_0){\rm exp}\{(t^{\frac{1}{8}}+\epsilon^\frac{1}{2})C_1(M(t))\}.
\end{split}\end{equation}

{\it Step 3.} Then we estimate $\|(\epsilon\nabla^2\sigma_t,\epsilon\nabla^2\theta_t)\|^2_{L^2_t(L^2)}$, $\|\epsilon\nabla^2{\rm div}\mathbf{u}\|^2_{L^\infty_t(L^2)}$ and $\|\epsilon^2{\rm curl}{\rm curl}\mathbf{u}_t\|^2_{L^\infty_t(H^1)}$.
We apply $\partial_i\nabla$ to ${\eqref{ns0}_1}$ and ${\eqref{ns0}_3}$ and multiply the resulting equality by $\epsilon^2\partial_i\nabla \sigma_t$ and $\epsilon^2\partial_i\nabla \theta_t$, respectively, in $L^2((0,t)\times\Omega)$ to get
\begin{equation}\label{dxns1a}\begin{split}
&\|\epsilon\partial_i\nabla \sigma_t\|^2_{L^2_t(L^2)}+\epsilon\int_0^t\int_\Omega\partial_i\nabla{\rm div}\mathbf{u}\cdot\partial_i\nabla \sigma_tdxds\\
&\quad=-\int_0^t\int_\Omega\epsilon^2
(\partial_i\nabla\mathbf{u}\cdot\nabla \sigma+\nabla\mathbf{u}\cdot\partial_i\nabla \sigma+\partial_i\mathbf{u}\cdot\nabla^2 \sigma+\mathbf{u}\cdot\partial_i\nabla^2 \sigma)\cdot\partial_i\nabla \sigma_tdxds\\
&\quad\quad-\int_0^t\int_\Omega\epsilon^2(\partial_i\nabla \sigma{\rm div}\mathbf{u}+\nabla \sigma\partial_i{\rm div}\mathbf{u}+\partial_i \sigma\nabla{\rm div}\mathbf{u}+\sigma\partial_i\nabla{\rm div}\mathbf{u})\cdot\partial_i\nabla \sigma_tdxds\\
&\quad\lesssim\|\epsilon\partial_i\nabla \sigma_t\|_{L^2_t(L^2)}(\|\mathbf{u}\|_{L^2_t(H^3)}
\|\epsilon\nabla\sigma\|_{L^\infty_t(L^3)}+\|\mathbf{u}\|_{L^\infty_t(H^2)}
\|\epsilon\nabla^2\sigma\|_{L^2_t(L^3)}\\
&\quad\quad+\|\mathbf{u}\|_{L^2_t(H^3)}
\|\epsilon\sigma\|_{L^\infty_t(L^\infty)}
+\|\mathbf{u}\|_{L^\infty_t(H^2)}\|\epsilon\partial_i\nabla^2 \sigma\|_{L^2_t(L^2)})\\
&\quad \lesssim\delta\|\epsilon\partial_i\nabla \sigma_t\|^2_{L^2_t(L^2)}+C_\delta
(\epsilon\|\mathbf{u}\|^2_{L^2_t(H^3)}
\|\epsilon\sigma\|_{L^\infty_t(H^2)}
\|\sigma\|_{L^\infty_t(H^1)}\\
&\quad\quad+t^\frac{1}{2}\|\mathbf{u}\|^2_{L^\infty_t(H^2)}
\|\epsilon\nabla^2\sigma\|_{L^\infty_t(L^2)}
\|\epsilon\nabla^2\sigma\|_{L^2_t(H^1)})\\
&\quad\quad+\epsilon^{\frac{1}{2}}\|\mathbf{u}\|^2_{L^2_t(H^3)}
\|\epsilon\sigma\|^{\frac{3}{2}}_{L^\infty_t(H^2)}
\|\sigma\|^{\frac{1}{2}}_{L^\infty_t(H^1)}
+\|\mathbf{u}\|^2_{L^\infty_t(H^2)}\|\epsilon\partial_i\nabla^2 \sigma\|^2_{L^2_t(L^2)})\\
&\quad \lesssim\delta\|\epsilon\partial_i\nabla \sigma_t\|^2_{L^2_t(L^2)}+C_\delta C_0(M_0){\rm exp}\{ ( t^{\frac{1}{8}}+\epsilon^\frac{1}{2})C_1(M(t))\},
\end{split}\end{equation}
 and
\begin{equation}\label{dxns1f}\begin{split}
&\|\epsilon\sqrt{\rho}\partial_i\nabla \theta_t\|^2_{L^2_t(L^2)}+\epsilon\int_0^t\int_\Omega\partial_i\nabla{\rm div}\mathbf{u}\cdot\partial_i\nabla \theta_tdxds\\
&\quad=-\int_0^t\int_\Omega\epsilon^2
\{\partial_i\nabla\rho(\theta_t+\mathbf{u}\cdot\nabla\theta)+\nabla\rho(\partial_i\theta_t
+\partial_i\mathbf{u}\cdot\nabla\theta+\mathbf{u}\cdot\partial_i\nabla\theta)\\
&\quad\quad+\partial_i\rho(\nabla\theta_t
+\nabla\mathbf{u}\cdot\nabla\theta+\mathbf{u}\cdot\nabla^2\theta)
+\rho(\partial_i\nabla\mathbf{u}\cdot\nabla\theta+\nabla\mathbf{u}\cdot\partial_i\nabla\theta
+\partial_i\mathbf{u}\cdot\nabla^2\theta+\mathbf{u}\cdot\partial_i\nabla^2\theta)\\
&\quad\quad+(\partial_i\rho\theta
+\rho\partial_i\theta+\partial_i\sigma)\nabla{\rm div}\mathbf{u}+(\rho\theta+\sigma)\partial_i\nabla{\rm div}\mathbf{u}+\kappa\partial_i\nabla\Delta\theta\\
&\quad\quad
+\partial_i\nabla[\epsilon(2\mu|D(\mathbf{u})|^2+\xi({\rm div}\mathbf{u})^2)]\}\cdot\partial_i\nabla \theta_tdxds\\
&\quad\lesssim\|\epsilon\partial_i\nabla \theta_t\|_{L^2_t(L^2)}[\|\rho\|_{L^\infty_t(H^2)}(\|\epsilon\theta_t\|_{L^2_t(L^\infty)}
+\|\epsilon\nabla\theta_t\|_{L^2_t(L^3)}+\|\mathbf{u}\|_{L^\infty_t(H^2)}
\|\epsilon\theta\|_{L^2_t(H^3)})\\
&\quad\quad+\|\mathbf{u}\|_{L^2_t(H^3)}(\|\rho\|_{L^\infty_t(H^2)}\|\epsilon\nabla\theta\|_{L^\infty_t(L^3)}
+\|\epsilon\nabla\sigma\|_{L^\infty_t(L^3)}+\|\rho\|_{L^\infty_t(H^2)}\|\epsilon\theta\|_{L^\infty_t(L^\infty)}\\
&\quad\quad
+\|\epsilon\sigma\|_{L^\infty_t(L^\infty)})+\|\epsilon^2\partial_i\nabla\Delta\theta\|_{L^2_t(L^2)}
+\epsilon\|\epsilon\mathbf{u}\|_{L^\infty_t(H^3)}
\|\mathbf{u}\|_{L^2_t(H^3)}]\\
&\quad \lesssim\delta\|\epsilon\sqrt{\rho}\partial_i\nabla \theta_t\|^2_{L^2_t(L^2)}+C_\delta\|\rho^{-1}\|_{L^\infty_{x,t}}\{\|\rho\|^2_{L^\infty_t(H^2)}
(t^{\frac{1}{4}}\|\epsilon\theta_t\|^{\frac{3}{2}}_{L^2_t(H^2)}\|\epsilon\theta_t\|^{\frac{1}{2}}_{L^\infty_t(H^1)}\\
&\quad\quad{+t^{\frac{1}{2}}\|\epsilon\theta_t\|_{L^2_t(H^2)}\|\epsilon\theta_t\|_{L^\infty_t(H^1)}}+\|\mathbf{u}\|^2_{L^\infty_t(H^2)}\|\epsilon\theta\|^2_{L^2_t(H^3)})\\
&\quad\quad
+\|\mathbf{u}\|^2_{L^2_t(H^3)}[\epsilon(\|\rho\|^2_{L^\infty_t(H^2)}
\|\epsilon\theta\|_{L^\infty_t(H^2)}\|\theta\|_{L^\infty_t(H^1)}
+\|\epsilon\sigma\|_{L^\infty_t(H^2)}\|\sigma\|_{L^\infty_t(H^1)})\\
&\quad\quad+\epsilon^{\frac{1}{2}}(\|\rho\|^2_{L^\infty_t(H^2)}
\|\epsilon\theta\|^{\frac{3}{2}}_{L^\infty_t(H^2)}\|\theta\|^{\frac{1}{2}}_{L^\infty_t(H^1)}
+\|\epsilon\sigma\|^{\frac{3}{2}}_{L^\infty_t(H^2)}\|\sigma\|^{\frac{1}{2}}_{L^\infty_t(H^1)})]\\
&\quad\quad+\|\epsilon^2\partial_i\nabla\Delta\theta\|^2_{L^2_t(L^2)}
+\epsilon\|\epsilon\mathbf{u}\|^2_{L^\infty_t(H^3)}
\|\mathbf{u}\|^2_{L^2_t(H^3)}\}\\
&\quad \lesssim\delta\|\epsilon\sqrt{\rho}\partial_i\nabla \theta_t\|^2_{L^2_t(L^2)}+C_\delta C_0(M_0){\rm exp}\{( t^{\frac{1}{8}}+\epsilon^\frac{1}{2})C_1(M(t))\},
\end{split}\end{equation}
by employing \eqref{lemma8}, \ef{dxns116}, \ef{dxns1c} and  \ef{dxns1d}.\\
We can obtain the estimate for $\|\e^2 {\rm curl} {\rm curl} \mathbf{u}_t\|_{L^2_t(H^1)}$ directly from ${\eqref{curl2}_1}$ with the help of Lemmas  \ref{hs1} and \ref{lembdry} as follows
\begin{equation}\label{curlut}\begin{split}
 \|\e^2 & {\rm curl} {\rm curl} \mathbf{u}_t\|^2_{L^2_t(H^1)}\equiv\|\e^2 {\rm curl} \omega_t\|^2_{L^2_t(H^1)}\\
&\lesssim\|\e^2 {\rm curl}{\rm curl} \omega_t\|^2_{L^2_t(L^2)}+\|\e^2 {\rm curl}\omega_t\cdot\mathbf{n}\|^2_{L^2_t( H^\frac 1 2 (\partial \Omega))}+\|\e^2 {\rm curl} \omega_t\|^2_{L^2_t(L^2)}\\
&\lesssim  \|\e^2 {\rm curl}{\rm curl} \omega_t\|^2_{L^2_t(L^2)}+ \e^2 \|\e \mathbf{u}_t\|^2_{L^2_t(H^2)}\\
&\leq C_0(M_0) {\rm exp}\{(t^{\frac{1}{8}} +\epsilon^\frac{1}{2}) C_1(M(t))\},
\end{split}\end{equation}
where $\omega_t= {\rm curl} \mathbf{u}_t$. Note that in the last inequality of \ef{curlut},  we use $\ef{curl2}_1$, \ef{lemma1}, \ef{lemma8}, \ef{dxns116}, \ef{dxns1b}, \ef{dxns1d} to derive that
\begin{equation*}\begin{split}
&\|\e^2 {\rm curl}  {\rm curl} \omega_t\|^2_{L^2_t(L^2)}\\
&\lesssim\|\rho\|^2_{L^\infty_{x,t}}[\|\epsilon^2\omega_{tt}\|^2_{L^2_t(L^2)}
+\epsilon^2\|\mathbf{u}\|^2_{L^\infty_t(H^2)}\|\epsilon\mathbf{u}_t\|^2_{L^2_t(H^2)}\\
&\quad+\epsilon^2\|\rho^{-1}\|^6_{L^\infty_{x,t}}\|\epsilon \sigma_t\|^2_{L^\infty_{t}(H^1)}\|\epsilon\sigma\|^2_{L^\infty_t(H^2)}(\|\mathbf{u}\|^2_{L^2_t(H^3)}
+\|\theta\|^2_{L^2_t(H^2)}+\| \epsilon\nabla \sigma\|^2_{L^\infty_{t}(H^1)}\|\theta\|^2_{L^2_t(H^2)})\\
&\quad
+\|\rho^{-1}\|^4_{L^\infty_{x,t}}(\|\epsilon\sigma_t\|^2_{L^\infty_{t}(H^1)}\|\epsilon^2\nabla{\rm div} \mathbf{u}\|^2_{L^2_t(H^2)}
+\|\epsilon^2\nabla \sigma\|^2_{L^\infty_{t}(H^2)}\|\epsilon\mathbf{u}_t\|^2_{L^2_t(H^2)}\\
&\quad+t\|\epsilon\sigma_t\|^2_{L^\infty_{t}(H^1)}\|\epsilon^2\mathbf{H}\|^2_{L^\infty_t(H^3)}
\|\mathbf{H}\|^2_{L^\infty_t(H^2)}+t\|\epsilon^2\nabla\sigma\|^2_{L^\infty_{t}(H^2)}\|\epsilon\mathbf{H}_t\|^2_{L^\infty_t(H^1)}
\|\mathbf{H}\|^2_{L^\infty_t(H^2)}\\
&\quad
+\|\epsilon\sigma_t\|^2_{L^\infty_{t}(H^1)}\|\epsilon^2{\rm curl}{\rm curl} \mathbf{u}\|^2_{L^2_t(H^2)}
+\|\epsilon\sigma_t\|^2_{L^\infty_{t}(H^1)}\|\epsilon\nabla \theta\|^2_{L^2_{t}(H^2)}
+\|\epsilon\theta_t\|^2_{L^\infty_{t}(H^1)}\|\epsilon\nabla \sigma\|^2_{L^2_{t}(H^2)}\\
&\quad
+\|\epsilon\sigma_t\|^2_{L^\infty_{t}(H^1)}(\|\epsilon \sigma\|^2_{L^\infty_{t}(H^2)}\|\epsilon\nabla\theta\|^2_{L^2_t(H^2)}
+\|\epsilon^2\nabla \sigma\|^2_{L^\infty_{t}(H^2)}\|\theta\|^2_{L^2_t(H^2)})\\
&\quad+\|\epsilon^2\nabla \sigma\|^2_{L^\infty_{t}(H^2)}\|\epsilon\sigma_t\|^2_{L^\infty_{t}(H^1)}
\|\theta\|^2_{L^2_t(H^2)}+\|\epsilon\nabla \sigma\|^2_{L^2_{t}(H^2)}\|\epsilon\theta_t\|^2_{L^\infty_{t}(H^1)}
\|\epsilon\sigma\|^2_{L^\infty_t(H^2)})\\
&\quad+t\|\epsilon^2\mathbf{H}_t\|^2_{L^\infty_t(H^2)}
\|\mathbf{H}\|^2_{L^\infty_t(H^2)}]\\
&\leq C_0(M_0) {\rm exp}\{(t^{\frac{1}{8}} +\epsilon^\frac{1}{2}) C_1(M(t))\},
\end{split}\end{equation*}
and  in the second inequality of \ef{curlut}, we apply Lemma \ref{lembdry} to ${\rm curl} {\rm curl} \mathbf{u}_t$ to give  that
\begin{equation*}\begin{split}
\|\e^2 {\rm curl}{\rm curl} \mathbf{u}_t\cdot\mathbf{n}\|^2_{L^2_t( H^\frac 1 2 (\partial \Omega))}\lesssim\epsilon^2\|\epsilon\mathbf{u}_t\|^2_{L^2_t(H^2)}.
\end{split}\end{equation*}
Then we differentiate $ {\eqref{ns0}_2}$ with respect to $x_i$ and $t$, and multiply the resulting equality by $\epsilon^2\partial_i\nabla{\rm div}\mathbf{u}$ in $L^2((0,t)\times\Omega)$ to obtain that
\begin{equation}\label{dxns1e}\begin{split}
(\mu &+\nu)\|\epsilon\partial_i\nabla{\rm div}\mathbf{u}\|^2_{L^2}(t)-\epsilon\int_0^t\int_\Omega(\partial_i\nabla \sigma_t+\partial_i\nabla\theta_t)\cdot\partial_i\nabla{\rm div}\mathbf{u}dxds\\
&=(\mu+\nu)\|\epsilon\partial_i\nabla{\rm div}\mathbf{u}_0\|^2_{L^2}+\epsilon^2\int_0^t\int_\Omega
[\partial_i\rho_t(\mathbf{u}_t
+\mathbf{u}\cdot\nabla\mathbf{u})+\rho_t(\partial_i\mathbf{u}_t
+\partial_i\mathbf{u}\cdot\nabla\mathbf{u}+\mathbf{u}\cdot\partial_i\nabla\mathbf{u})\\
&\quad
+\p_i\rho( \mathbf{u}_{tt}
+  \mathbf{u}_t\cdot\nabla\mathbf{u}+
 \mathbf{u} \cdot\nabla\mathbf{u}_t
 )
+\rho(\partial_i\mathbf{u}_{tt}
+ \partial_i\mathbf{u}_t\cdot\nabla\mathbf{u}
+\mathbf{u}_t\cdot\partial_i\nabla\mathbf{u}
+ \partial_i\mathbf{u} \cdot\nabla\mathbf{u}_t
+\mathbf{u}\cdot\partial_i\nabla\mathbf{u}_t)\\
&\quad+\partial_i\sigma_t\nabla\theta+\partial_i\sigma\nabla\theta_t+
\sigma_t\partial_i\nabla\theta+\sigma\partial_i\nabla\theta_t+\theta_t\partial_i\nabla\theta
+\theta\partial_i\nabla\sigma_t+\nabla\sigma_t\partial_i\theta+\nabla\sigma\partial_i\theta_t]
\cdot\partial_i\nabla{\rm div}\mathbf{u}dxds\\
&\quad +\epsilon^2\int_0^t\int_\Omega\partial_i{\rm curl}{\rm curl}\mathbf{u}_t\cdot\partial_i\nabla{\rm div}\mathbf{u}dxds\\
& \lesssim \|\mathbf{u}_0\|^2_{H^3}+ \|{\rm div}\mathbf{u}\|_{L^2_t(H^2)} [\|\epsilon\sigma_t\|_{L^\infty_t(H^1)} (\|\epsilon\mathbf{u}_t\|_{L^2_t(H^2)}
+t^\frac 1 4 \|\mathbf{u}\|_{L^\infty_t(H^2)}^\frac 3 2 \|\mathbf{u}\|_{L^2_t(H^3)}^\frac 1 2\\
&\quad
+t^\frac 1 8 \|\mathbf{u}\|_{L^\infty_t(H^2)}^\frac 5 4 \|\mathbf{u}\|_{L^2_t(H^3)}^\frac 3 4)
 +\|\rho\|_{L^\infty_{t}(H^2)}(\|\epsilon^2\mathbf{u}_{tt}\|_{L^2_t(H^1)}
+\|\epsilon\mathbf{u}_t\|_{L^2_t(H^2)}\|\mathbf{u}\|_{L^\infty_t(H^2)})\\
&\quad+\|\epsilon\sigma_t\|_{L^2_t(H^2)}\|\epsilon\theta\|_{L^\infty_t(H^2)}
+\|\epsilon\theta_t\|_{L^2_t(H^2)}\|\epsilon\sigma\|_{L^\infty_t(H^2)}+ \|\e^2 {\rm curl} {\rm curl} \mathbf{u}_t\|_{L^2_t(H^1)} ]\\
&\lesssim\delta\|{\rm div}\mathbf{u}\|^2_{L^2_t(H^2)}+C_\delta C_0(M_0){\rm exp}\{(t^\frac 1 8+\epsilon^\frac{1}{2}) C_1(M(t))\},
\end{split}\end{equation}
with the help of  \eqref{lemma1}, \eqref{lemma8}, \eqref{lemma5} and { \eqref{curlut}}.
By virtue of $\ef{ns0}_1$, \eqref{lemma8} and \eqref{dxns1d}, we get
\begin{equation}\label{u5}\begin{split}
\|{\rm div}\mathbf{u}\|^2_{L^2_t(H^2)}
&\lesssim \|\epsilon \sigma_t\|^2_{L^2_t(H^2)}
+\|\mathbf{u}\|^2_{L^\infty_t(H^2)}\|\epsilon\nabla \sigma\|^2_{L^2_t(H^2)}
+\epsilon^\frac{1}{2}\|\epsilon\sigma\|^\frac{3}{2}_{L^\infty_t(H^2)}
\|\sigma\|^\frac{1}{2}_{L^\infty_t(H^1)}
\|\mathbf{u}\|^2_{L^2_t(H^3)}\\
& \lesssim \|\epsilon \nabla^2 \sigma_t\|^2_{L^2_t(L^2)} + C_0(M_0){\rm exp}\{(t^\frac 1 8+\epsilon^\frac{1}{2}) C_1(M(t))\}.
\end{split}\end{equation}
Summarizing ${\eqref{dxns1a}}$ and ${\eqref{dxns1e}}$, then plugging \ef{u5} into the resulting inequality and choosing $\delta$ to be small, one  derives
\begin{equation}\label{dxns1f}\begin{split}
&\|\epsilon\partial_i\nabla \sigma_t\|^2_{L^2_t(L^2)}+\|\epsilon\partial_i\nabla \theta_t\|^2_{L^2_t(L^2)}+\|\epsilon\partial_i\nabla{\rm div}\mathbf{u}\|^2_{L^2}(t) \leq C_0(M_0){\rm exp}\{ (t^\frac 1 8+\epsilon^\frac{1}{2})C_1(M(t))\}.
\end{split}\end{equation}
By combining the estimates ${\eqref{curlu1}}$, \eqref{lemma8}, \ef{dxns1f} and ${\eqref{u5}}$, we find
\begin{equation}\label{u6}\begin{split}
\|\mathbf{u}\|^2_{L^2_t({H^3})}
&\lesssim\|{\rm div}\mathbf{u}\|^2_{L^2_t({H^2})}+\|{\rm curl}\mathbf{u}\|^2_{L^2_t({H^2})}
+\|\mathbf{u}\cdot\mathbf{n}\|^2_{L^2_t(H^{\frac{5}{2}}(\partial\Omega))}
+\|\mathbf{u}\|^2_{L^2_t({H^2})}\\
&\leq C_0(M_0){\rm exp}\{( t^{\frac{1}{8}}+\epsilon^\frac{1}{2})C_1(M(t))\},
\end{split}\end{equation}
where we have used Lemma \ref{hs1} and the boundary condition $\mathbf{u}\cdot\mathbf{n}|_{\partial\Omega}=0$.\\
From  ${\eqref{dxns116}}$, \eqref{dxns1b}, ${\eqref{dxns1f}}$ and \ef{u6} we conclude ${\eqref{lemma6}}$ immediately.

\hfill$\Box$

\begin{lem} For the solution to the system ${\eqref{ns0}}$, we have
\begin{equation}\label{lemma71}\begin{split}
\|\sigma\|^2_{L^2_t(H^2)}
\leq C_0(M_0){\rm exp}\{ ( t^{\frac{1}{8}}+\epsilon^\frac 1 2 )C_1(M(t))\},\q \e\in (0,1],\ t\in [0,1],
\end{split}\end{equation}
\begin{equation}\label{lemma72}\begin{split}
\|\epsilon^2\mathbf{H}\|^2_{L^\infty_t(H^3)}+\|\epsilon^2\mathbf{H}_t\|^2_{L^\infty_t(H^2)}
\leq C_0(M_0){\rm exp}\{ ( t^{\frac{1}{8}}+\epsilon^\frac 1 2 )C_1(M(t))\},\q \e\in (0,1],\ t\in [0,1].
\end{split}\end{equation}
\end{lem}
\textbf{Proof}. Observe that
\begin{equation}\label{q21p}\begin{split}
\|\epsilon\nabla^2(\sigma\theta)\|_{L^2_t(L^2)}^2& \lesssim  \|\epsilon \sigma\nabla^2\theta\|_{L^2_t(L^2)}^2 +\|\epsilon\nabla \sigma\cdot\nabla\theta\|_{L^2_t(L^2)}^2 + \|\epsilon\theta\nabla^2 \sigma\|_{L^2_t(L^2)}^2\\
&\lesssim  \epsilon^\frac{1}{2}\|\epsilon  \sigma\|_{L^\infty_t(H^2)}^\frac{3}{2}\|\sigma\|_{L^\infty_t(H^1)}^\frac{1}{2}\| \nabla^2\theta\|_{L^2_t(L^2)}^2\\
 &\quad + \epsilon\|\epsilon \nabla \sigma\|_{L^\infty_t(H^1)} \| \nabla \sigma\|_{L^\infty_t(L^2)} \| \nabla\theta\|_{L^2_t(H^1)}^2\\
&\quad + t^\frac{1}{4}\|\epsilon \nabla^2\sigma\|_{L^\infty_t(L^2)}^2\| \theta\|_{L^2_t(H^2)}^{\frac 3 2} \|\theta\|_{L^\infty_t(H^1)}^{\frac 1 2}\\
&\lesssim   (\epsilon^\frac{1}{2}+ t^\frac{1}{4})M^2(t).
\end{split}\end{equation}
By virtue of \eqref{lemma1}, \eqref{curlu3}, \eqref{lemma8} and \ef{lemma6}, we apply \lq\lq$\epsilon\nabla$" to ${\eqref{ns0}}_2$ to get
\begin{equation}\label{q21}\begin{split}
\|\nabla^2\sigma+\nabla^2\theta\|^2_{L^2_t(L^2)}
&\lesssim \|\epsilon\nabla^2{\rm div}\mathbf{u}\|^2_{L^2_t(L^2)}
+ \|\epsilon\nabla{\rm curl}{\rm curl}\mathbf{u}\|^2_{L^2_t(L^2)}+{\|\epsilon\nabla^2(\sigma\theta)\|^2_{L^2_t(L^2)}}\\
&\quad+\|\rho\|^2_{L^\infty_t(H^2)}
(\|\epsilon\mathbf{u}_t\|^2_{L^2_t(H^1)}
+\epsilon^2\|\mathbf{u}\|^4_{L^\infty_t(H^2)})\\
&\leq C_0(M_0){\rm exp}\{( t^{\frac{1}{8}}+\epsilon^\frac 1 2 )C_1(M(t))\}.
\end{split}\end{equation}
Therefore, \eqref{lemma71} follows from  ${\eqref{lemma8}}$ and ${\eqref{q21}}$.

Let $\alpha$ be a multi-index of spatial derivatives with $0\leq|\alpha|\leq3$. Applying $\partial^\alpha$ to ${\eqref{ns0}}_3$ and multiplying the resulting equality by $\epsilon^4\partial^\alpha\mathbf{H}$ in $L^2((0,t)\times\Omega)$ lead to
\begin{equation}\label{high1}\begin{split}
\frac{1}{2}&\|\epsilon^2\partial^\alpha\mathbf{H}\|^2_{L^2}(t)\\
&=\frac{1}{2}\|\epsilon^2\partial^\alpha\mathbf{H}(0)\|^2_{L^2}
+\int_0^t\int_\Omega{\rm div}\mathbf{u}|\epsilon^2\partial^\alpha\mathbf{H}|^2dxds\\
&\quad+\epsilon^4\int_0^t\int_\Omega\{\partial^\alpha[(\mathbf{H}\cdot\nabla)\mathbf{u}]
-[\partial^\alpha,\mathbf{u}\cdot\nabla]\mathbf{H}-\partial^\alpha(\mathbf{H}{\rm div}\mathbf{u})\}\cdot\partial^\alpha\mathbf{H}dxds\\
&\leq C_0(M_0)+Ct^\frac{1}{2}\|\epsilon^2\partial^\alpha\mathbf{H}\|_{L^\infty_t(L^2)}
(\|\mathbf{u}\|_{L^2_t(H^3)}\|\epsilon^2\mathbf{H}\|_{L^\infty_t(H^3)}
+\|\epsilon^2\mathbf{u}\|_{L^2_t(H^4)}\|\mathbf{H}\|_{L^\infty_t(H^2)}).
\end{split}\end{equation}
Thus we sum over $\alpha$ to obtain that
\begin{equation}\label{high2}\begin{split}
\|\epsilon^2\mathbf{H}\|^2_{L^\infty_t(H^3)}\leq C_0(M_0){\rm exp}\{( t^{\frac{1}{8}}+\epsilon^\frac 1 2 )C_1(M(t))\}.
\end{split}\end{equation}
Then from ${\eqref{ns0}}_3$ we directly derive
\begin{equation}\label{high3}\begin{split}
\|\epsilon^2\mathbf{H}_t\|^2_{L^\infty_t(H^2)}&\leq\epsilon^2
\|\epsilon\mathbf{u}\|^2_{L^\infty_t(H^3)}\|\mathbf{u}\|^2_{L^\infty_t(H^2)}
+\|\mathbf{u}\|^2_{L^\infty_t(H^2)}\|\epsilon^2\mathbf{H}\|^2_{L^\infty_t(H^3)}\\
&\leq C_0(M_0){\rm exp}\{( t^{\frac{1}{8}}+\epsilon^\frac 1 2 )C_1(M(t))\}.
\end{split}\end{equation}
\hfill$\Box$

Collecting the estimates in this section, we obtain Proposition \ref{prop} immediately.

\section{Proof of Main Theorems}

\subsection{Proof of Theorem \ref{unithm}}

The proof is in the spirit of Alazard \cite{Alazard} and M\'etivier and Schochet \cite{Metivier}. Assume $T^\epsilon<+\infty$ is the maximal life time of existence for the solution obtained in Theorem \ref{main1}. Then by the uniform estimates in Proposition \ref{prop} for any $0\leq t\leq \min\{T^\epsilon,1\}$, we have
\begin{equation}\label{pf1}
M^\epsilon(t)\leq C_0(M^\epsilon_0){\rm exp}\{ (t^\frac{1}{8}+\epsilon^\frac 1 2 )C_1(M^\epsilon(t))\},
\end{equation}
where $M^\epsilon_0\leq D_0$ for $0<\epsilon\leq1$. In the sequel, we choose $C(D_0)<D<+\infty$ then $\epsilon_1\in(0,1]$ and $T_1\in(0,1]$ such that
\begin{equation}\label{pf2}
C_0(D_0){\rm exp}\{ (T_1^\frac{1}{8}+\epsilon_1^\frac 1 2 )C_1(D)\}<D.
\end{equation}

Let $t<\min\{T^\epsilon,T_1\}$. By combining the inequalities ${\eqref{pf1}}$ and ${\eqref{pf2}}$ with the hypothesis
$M^\epsilon(0)=M^\epsilon_0$, we have that $M^\epsilon(t)\neq D$. Otherwise, since $C_0(\cdot)$ and $C_1(\cdot)$ are nondecreasing, by the inequalities ${\eqref{pf1}}$ and ${\eqref{pf2}}$, we have
$$
D=M^\epsilon(t)\leq C_0(M^\epsilon_0){\rm exp}\{ (t^\frac{1}{8}+\epsilon^\frac 1 2 )C_1(M^\epsilon(t))\}<D,
$$
which leads to a contradiction.

Besides, we can assume without restriction that $D_0<D$, so that $M^\epsilon_0<D$. Since the function $M(t)$ is continuous, we obtain
\begin{equation}\label{pf3}
M^\epsilon(t)<D\ {\rm for}\ t<\min\{T^\epsilon,T_1\},\ 0<\epsilon_1\leq1.
\end{equation}

Consequently, by the standard continuation principle, we have $T^\epsilon\geq  T_1$ for $0<\epsilon\leq\epsilon_1$. Therefore, $M^\epsilon(t)<D$ for any $t\in[0,T_1]$ where $T_1$ is independent of $\epsilon\in(0,\epsilon_1]$. \hfill$\Box$

\subsection{Proof of weak convergence in Theorem \ref{convergence}}

In the spirit of Secchi's argument \cite{Secchi}, we assume that $\{(\sigma^\epsilon,\mathbf{u}^\epsilon,\theta^\epsilon)\}_{\{\epsilon>0\}}$ is a sequence of local strong solutions to   \eqref{ns0} where the initial data $\{(\sigma^\epsilon_0,\mathbf{u}^\epsilon_0,\theta^\epsilon_0)\}_{\{\epsilon>0\}}$ satisfy the assumptions in Theorem \ref{unithm}.  Then, due to the uniform estimates \eqref{uniform},  we have the following convergence results up to a subsequence as $\e \to 0$:
\be\label{con3}
\mathbf{u}^\epsilon\rightharpoonup\mathbf{v}\quad weakly\ {\rm in}\ L^2(0,T;H^3(\Omega))\ {\rm and}\  weakly-\ast\ {\rm in}\ L^\infty(0,T;H^2(\Omega)),
\ee
\be \label{con4}
\sigma^\e  \rightharpoonup \sigma \q {weakly\ {\rm in}\ L^2(0,T;H^2(\Omega))}\ {\rm and}\  weakly-\ast\ {\rm in}\ L^\infty(0,T;{H^1}(\Omega)),
\ee
\be\label{con41}
\theta^\epsilon\rightharpoonup\theta\quad weakly\ {\rm in}\ L^2(0,T;H^3(\Omega))\ {\rm and}\  weakly-\ast\ {\rm in}\ L^\infty(0,T;{H^1}(\Omega)),
\ee
\be\label{con42}
\mathbf{H}^\epsilon\rightharpoonup\mathbf{h}\quad  weakly-\ast\ {\rm in}\ L^\infty(0,T;{H^2}(\Omega)),
\ee
and
\be\label{con5}
\rho^\e=1+\e \sigma^\e  \to 1 \q   {\rm in}\ L^2(0,T;{H^2}(\Omega))\cap L^\infty(0,T;{H^1}(\Omega)),
\ee
\be\label{con6}
\mathcal{T}^\e=1+\e \theta^\e  \to 1 \q   {\rm in}\ L^2(0,T;{H^2}(\Omega))\cap L^\infty(0,T;{H^1}(\Omega)).
\ee

Consider the orthonormal decomposition $L^2(\Omega)=H\oplus G$ where
$$
H=\{u\in L^2(\Omega)|\nabla\cdot u=0\ {\rm in}\ \Omega,\ u\cdot n=0\ {\rm on}\ \partial\Omega\},\quad G=\{\nabla\psi|\psi\in H^1(\Omega)\}.
$$
Let $P$ be the projection onto $H$ and $Q=I-P$. Below we will use the fact that
\begin{equation}\label{fact}
\|Pu\|_s\leq\|u\|_s
\end{equation}
for any number $s\geq0$ (see Hu and Wang \cite{HW2009}; Temam \cite{Temam}). Indeed, from the results of Galdi \cite{Galdi} (see also Hu and Wang \cite{HW2009}), we know that the operators $P$ and $Q$ are linear bounded operators in $W^{s,p}(\Omega)$ for all $s\ge 0$ and $1<p<\infty$ in a bounded domain $\Omega$ with smooth boundary $\partial \Omega$. In particular, $P\in\mathcal{L}(H^s(\Omega),H^s(\Omega))$ for $s\ge 0.$

Next, observing the equivalence of $\eqref{ns}$ and $\eqref{ns0}$, we show the converge of $\eqref{ns}$ for convenience. By virtue of
$\eqref{ns}_1$, we obtain
\be
\begin{split}
\rho^\e {\rm div}\mathbf{u}^\epsilon
&=- \rho^\e_t - \mathbf{u}^\epsilon\cdot\nabla \rho^\e.
\end{split}
\ee
Due to \ef{con3} and \ef{con5}, we have
$\rho^\e {\rm div}\mathbf{u}^\epsilon\rightarrow {\rm div}\mathbf{v}$ and $-\rho^\e_t\to 0$ in $\mathcal{D}'((0,T)\times\Omega)$ as $\epsilon\rightarrow0$. Moreover, we derive from \ef{uniform} that $\mathbf{u}^\epsilon\cdot\nabla \rho^\e \to 0$  in $L^\infty(0,T;H^1(\Omega))$ as $\epsilon\rightarrow0$.
With the uniqueness of limit, we find
\be\label{div}
{\rm div}\mathbf{v}=0\q a.e.\q (0,T)\times\Omega,
 \ee
 which yields  $\mathbf{v}=P\mathbf{v}$ and $Q\mathbf{v}=0$ $a.e.$  in $(0,T)\times\Omega$. Besides, substituting \eqref{nonisentropic} into $\eqref{ns}_3$, we get
\be\label{temp0}
\begin{split}
p^\e{\rm div}\mathbf{u}^\epsilon&=\epsilon^2(2\mu|D(\mathbf{u}^\epsilon)|^2+\xi({\rm div}\mathbf{u}^\epsilon)^2)+{\rm div}(\kappa\nabla\mathcal{T}^\e)-\rho^\e\mathcal{T}^\e_t-\rho^\e\mathbf{u}^\e\cdot\nabla\mathcal{T}^\e\\
&=-\rho^\e\mathcal{T}^\e_t+\epsilon[\epsilon(2\mu|D(\mathbf{u}^\epsilon)|^2+\xi({\rm div}\mathbf{u}^\epsilon)^2)+\kappa\Delta\theta^\e-\rho^\e\mathbf{u}^\e\cdot\nabla\theta^\e]\\
&=:-\rho^\e\mathcal{T}^\e_t+\epsilon(\ast).
\end{split}
\ee
By utilizing the convergence results \ef{con3}, \ef{con5} and \ef{con6}, we have $p^\e {\rm div}\mathbf{u}^\epsilon\rightarrow {\rm div}\mathbf{v}$ and $-\rho^\e\mathcal{T}^\e_t=-\mathcal{T}^\e_t-\epsilon\sigma^\e\mathcal{T}^\e_t\to 0$ in $\mathcal{D}'((0,T)\times\Omega)$ as $\epsilon\rightarrow0$. Moreover, from the uniform estimates in \ef{uniform}, $(\ast)$ is bounded in $L^2(0,T;L^2(\Omega))$. Thus, as $\epsilon\rightarrow0$, we obtain
\be\label{div0}
{\rm div}\mathbf{v}=0\q a.e.\q (0,T)\times\Omega,
 \ee
which can be covered by \ef{div}.

We project ${\eqref{ns}}_2$  onto $H$ by the operator $P$ to obtain that
\be\label{peqn}
\partial_tP(\rho^\epsilon\mathbf{u}^\epsilon)=-P({\rm div}(\rho^\epsilon\mathbf{u}^\epsilon\otimes
\mathbf{u}^\epsilon))+\mu P\Delta\mathbf{u}^\epsilon+P({\rm curl}\mathbf{H}^\epsilon\times\mathbf{H}^\epsilon)
\ee
is bounded in $L^\infty(0,T;L^2(\Omega))\cap L^2(0,T;H^1(\Omega))$.
Then by Lemma \ref{als}, $P(\rho^\epsilon\mathbf{u}^\epsilon)$ converges to $\mathbf{v}$ strongly in $L^\infty(0,T;H^2(\Omega))\cap L^2(0,T;H^3(\Omega))$.

Note that
\begin{equation}\label{pu}\begin{split}
\partial_tP\mathbf{u}^\epsilon=\partial_tP(-\rho^\epsilon\mathbf{u}^\epsilon
+\mathbf{u}^\epsilon)+\partial_tP(\rho^\epsilon\mathbf{u}^\epsilon)=P(-\epsilon\sigma^\epsilon\partial_t\mathbf{u}^\epsilon
-\epsilon\sigma^\epsilon_t\mathbf{u}^\epsilon)
+\partial_tP(\rho^\epsilon\mathbf{u}^\epsilon)
\end{split}\end{equation}
is bounded in $L^\infty(0,T;L^2(\Omega))\cap L^2(0,T;H^1(\Omega))$.
Thus, by virtue of Lemma \ref{als}, we have up to a subsequence that for any $\delta\in(0,1)$
\begin{equation}\label{con1}\begin{split}
P\mathbf{u}^\epsilon\rightarrow\mathbf{v}\q  {\rm in}\ L^2(0,T;H^{3-\delta}(\Omega))\cap {\bl C([0,T],H^{2-\delta}(\Omega))},\q \e \to 0,
\end{split}\end{equation}
\begin{equation}\label{con1}\begin{split}
P\mathbf{u}^\epsilon\rightharpoonup\mathbf{v}\quad weakly\ {\rm in}\ L^2(0,T;H^3(\Omega))\ {\rm and}\  weakly-\ast\ {\rm in}\ L^\infty(0,T;H^2(\Omega)), \q \e \to 0,
\end{split}\end{equation}
and
\begin{equation}\label{con2}\begin{split}
Q\mathbf{u}^\epsilon\rightharpoonup0\quad weakly\ {\rm in}\ L^2(0,T;H^3(\Omega))\ {\rm and}\  weakly-\ast\ {\rm in}\ L^\infty(0,T;H^2(\Omega))
, \q \e \to 0.
\end{split}\end{equation}

Integrating the equality \eqref{peqn} with respect to $t$,  one achieves
\begin{equation}\label{q5}\begin{split}
P(\rho^\epsilon\mathbf{u}^\epsilon)(t)=P(\rho^\epsilon_0\mathbf{u}^\epsilon_0)
-\int_0^tP({\rm div}(\rho^\epsilon\mathbf{u}^\epsilon\otimes
\mathbf{u}^\epsilon))ds
+\int_0^tP(\mu\Delta\mathbf{u}^\epsilon +{\rm curl}\mathbf{H}^\epsilon\times\mathbf{H}^\epsilon)ds.
\end{split}\end{equation}

{
Note that
\begin{equation}\label{project2}\begin{split}
{\rm div}(\rho^\epsilon\mathbf{u}^\epsilon\otimes
\mathbf{u}^\epsilon)={\rm div}(\rho^\epsilon\mathbf{u}^\epsilon\otimes
P\mathbf{u}^\epsilon)+{\rm div}(P(\rho^\epsilon\mathbf{u}^\epsilon)\otimes
Q\mathbf{u}^\epsilon)+{\rm div}(Q(\rho^\epsilon\mathbf{u}^\epsilon)\otimes
Q\mathbf{u}^\epsilon).
\end{split}\end{equation}

Combining the weak and strong convergence in \eqref{con1} and \eqref{con2} we find
\be
{\rm div}(\rho^\epsilon\mathbf{u}^\epsilon\otimes
P\mathbf{u}^\epsilon)\rightharpoonup {\rm div}(\mathbf{v}\otimes\mathbf{v})\ {\rm in}\  L^2(0,T;H^1(\Omega)),\; \e\to 0,
\ee
since $\rho^\epsilon\mathbf{u}^\epsilon$ converges  to $\mathbf{v}$ in $L^2(0,T;H^3(\Omega))$,
and
\be
{\rm div}(P(\rho^\epsilon\mathbf{u}^\epsilon)\otimes
Q\mathbf{u}^\epsilon)\rightharpoonup 0\ {\rm in}\  L^2(0,T;H^1(\Omega)),\; \e\to 0.
\ee

Next, we need to estimate the remaining part ${\rm div}(Q(\rho^\epsilon\mathbf{u}^\epsilon)\otimes
Q\mathbf{u}^\epsilon)$ in the spirit of {\bl\cite{LM}}. By virtue of \ef{con3} and \ef{con5}, we can deduce that
$$
\rho^\epsilon\mathbf{u}^\epsilon-\mathbf{u}^\epsilon
=(\rho^\epsilon-1)\mathbf{u}^\epsilon\rightarrow0\ {\rm in}\  L^2(0,T;H^2(\Omega)),\ \e\to 0,
$$
which implies that
$$
Q(\rho^\epsilon\mathbf{u}^\epsilon)\rightarrow Q\mathbf{u}^\epsilon\ {\rm in}\  L^2(0,T;H^2(\Omega)),\ \e\to 0.
$$
Then we arrive at
$$
{\rm div}((Q(\rho^\epsilon\mathbf{u}^\epsilon)-Q\mathbf{u}^\epsilon)\otimes
Q\mathbf{u}^\epsilon)\rightharpoonup0\ {\rm in}\  L^2(0,T;H^1(\Omega)),\ \e\to 0.
$$
Then we rewrite ${\rm div}(Q\mathbf{u}^\epsilon\otimes
Q\mathbf{u}^\epsilon)$  as
$(Q\mathbf{u}^\epsilon\cdot\nabla)Q\mathbf{u}^\epsilon$
and set $Q\mathbf{u}^\epsilon=\nabla\psi^\epsilon$ for some function $\psi^\epsilon$ such that $(Q\mathbf{u}^\epsilon\cdot\nabla)Q\mathbf{u}^\epsilon
=\frac{1}{2}\nabla|\nabla\psi^\epsilon|^2$. Thus, the term $(Q\mathbf{u}^\epsilon\cdot\nabla)Q\mathbf{u}^\epsilon$ can be incorporated in the pressure $\pi$. Noticing that
$$
{\rm div}(Q(\rho^\epsilon\mathbf{u}^\epsilon)\otimes
Q\mathbf{u}^\epsilon)={\rm div}((Q(\rho^\epsilon\mathbf{u}^\epsilon)-Q\mathbf{u}^\epsilon)\otimes
Q\mathbf{u}^\epsilon)+{\rm div}(Q\mathbf{u}^\epsilon\otimes
Q\mathbf{u}^\epsilon),
$$
the term ${\rm div}(Q(\rho^\epsilon\mathbf{u}^\epsilon)\otimes
Q\mathbf{u}^\epsilon)$ can be incorporated in the pressure $\pi$.
On the other hand, owing to \ef{con3}, and \ef{div}, we obtain
\begin{equation}\label{project3}
\mu\Delta\mathbf{u}^\epsilon
+\nu\nabla{\rm div}\mathbf{u}^\epsilon\rightharpoonup\mu\Delta\mathbf{v}\quad weakly\ {\rm in}\ L^2(0,T;H^1(\Omega))\ {\rm and}\ \ weakly-\ast\ in\ L^\infty(0,T;L^2),
\end{equation}
as $\e \to 0$. {\bl From the bound in \ef{uniform},} we deduce that $\mathbf{H}^\epsilon$ is bounded in $ L^\infty(0,T;H^2(\Omega))$ and
\begin{equation}\label{conh1}
\partial_t\mathbf{H}^\epsilon={\rm curl}(\mathbf{u}^\epsilon\times\mathbf{H}^\epsilon)
\end{equation}
is bounded in $L^\infty(0,T;L^2(\Omega))$. Thus, by Lemma \ref{als}, we have up to a subsequence that
\begin{equation}\label{conh2}
\mathbf{H}^\epsilon\rightarrow\mathbf{h}\  {\rm in}\  C([0,T],H^{2-\delta}(\Omega)),\quad \delta\in(0,1),
\end{equation}
as $\e \to 0$, which gives that
\begin{equation}\label{conh3}
{\rm curl}\mathbf{H}^\epsilon\times\mathbf{H}^\epsilon{\bl\rightharpoonup}{\rm curl}\mathbf{h}\times\mathbf{h}\  weakly-\ast\  {\rm in}\  L^\infty(0,T;H^1(\Omega)),\quad \epsilon\rightarrow0.
\end{equation}
Then we pass to the limit in \eqref{peqn} to get
\begin{equation}\label{q6}
\mathbf{v}(t)=P\mathbf{w}_0-\int_0^tP({\rm div}(\mathbf{v}\otimes\mathbf{v}))ds+\int_0^tP(\mu\Delta\mathbf{v}(s)
+{\rm curl}\mathbf{h}(s)\times\mathbf{h}(s))ds,
\end{equation}
which is the abstract form of $\eqref{ins}_1$ for the hydrostatic pressure function $\pi(x,t)$ and the initial value $P\mathbf{w}_0$. Integrating \ef{conh1} in $[0,T]$ and passing to the limit indicate the abstract form of $\eqref{ins}_3$ with initial value $\mathbf{h}_0$ that
\begin{equation}\label{conh4}
\mathbf{h}(t)=\mathbf{h}_0-\int_0^t{\rm curl}(\mathbf{u}(s)\times\mathbf{h}(s))ds.
\end{equation}
In addition, we deduce that $\mathbf{v}\in L^2(0,T;H^3(\Omega))\cap {\bl L^\infty(0,T;H^2(\Omega))}$ and $\mathbf{h}\in {\bl L^\infty(0,T;H^2(\Omega))}$.

\hfill$\Box$

\vskip 5mm
\noindent {\bf Acknowledgement.} This work was supported in part by  National Natural Science Foundation of China through grants 11971477, 12131007, and 11761141008.

\vskip 5mm
\noindent {\bf Data availability.} Data sharing not applicable to this article as no datasets were generated or analysed during the current study.

\vskip 5mm
\noindent {\bf Conflict of interest.} On behalf of all authors, the corresponding author states that there is no conflict of interest.

\end{document}